\documentclass[11pt]{amsart}
\usepackage{fullpage}
\usepackage[utf8]{inputenc}
\usepackage{microtype}
\usepackage{graphicx}

\usepackage{amsmath,amssymb,stmaryrd,array,multirow,dsfont,mathrsfs}
\usepackage{cite}

\usepackage{scalerel,stackengine}
\usepackage{mathtools}

\usepackage{amsthm}
\usepackage{enumerate}
\usepackage{xcolor}
\usepackage{color}
\usepackage{hyperref}
\usepackage{indentfirst}
\usepackage{pifont}
\usepackage[T1]{fontenc}
\usepackage{microtype}
\usepackage{graphics}
\usepackage{amsthm}
\usepackage{extarrows}
\usepackage{stmaryrd}
\graphicspath{{Images/}}
\xdefinecolor{darkgreen}{rgb}{0,0.4,0}
\usepackage{tikz}

\newcommand{\be}{\begin{equation} }
\newcommand{\ee}{\end{equation}}
\newcommand{\bse}{\begin{subequations}}
\newcommand{\ese}{\end{subequations}}

\newcommand\bpm{\begin{pmatrix}}
\newcommand\epm{\end{pmatrix}}
\newtheorem{thm}{Theorem}[section]
\newtheorem{lem}[thm]{Lemma}
\newtheorem{prop}[thm]{Proposition}
\newtheorem{rmk}[thm]{Remark}

\newcommand{\beq}{\begin{equation}}
\newcommand{\eeq}{\end{equation}}
\newcommand{\beqs}{\begin{equation*}}
\newcommand{\eeqs}{\end{equation*}}
\newcommand{\ben}{\begin{eqnarray}}
\newcommand{\een}{\end{eqnarray}}
\newcommand{\beno}{\begin{eqnarray*}}
\newcommand{\eeno}{\end{eqnarray*}}
\newcommand{\mr}{\mathbb{R}}
\newcommand{\m}{\mathbb}
\newcommand{\p}{\partial}

\newcommand{\ang}[1]{\langle #1\rangle}
\newcommand {\Div}{\text{div}}

\def \d {\mathrm {d}}

\def\cC{{\mathcal C}}
\def\cD{{\mathcal D}}

\def\cI{{\mathcal I}}

\def\cK{{\mathcal K}}

\def\cN{{\mathcal N}}
\def\cO{{\mathcal O}}

\def\cS{{\mathcal S}}

\def\cU{{\mathcal U}}
\def\cV{{\mathcal V}}

\def\cY{{\mathcal Y}}

\let\f=\frac
\def \p {\partial}
\def\mR {\mathbb{R}}

\def \vep {\varepsilon}
\def \diag {\text{diag}}
\def\vp{\varphi}
\def \kpa {\kappa}
\def \T {\mathbb{T}}
\def \m {\mathfrak m}
\def \rM {\mathrm{M}}
\def \rD {\mathrm{D}}
\def \rA {\mathrm{A}}
\def \rF {\mathrm{F}}

\def \rG {\mathrm{G}}
\def \rN {\mathrm{N}}
\def \rQ {\mathrm{Q}}
\def \rV {\mathrm{V}}
\def \rB  {\mathrm{B}}
\def \rU {\mathrm{U}}
\def \Id {\textrm{Id}}

\def \tah {\mathrm{th}}

\def \per{\mathrm{per}}
\def \i {\mathrm{i}}

\def \sech {\mathrm{sech}}
\def \sh {\mathrm {sh}}
\def \ch {\mathrm {ch}}

\newtheorem{theorem}{Theorem}[section]
\newtheorem{corollary}[theorem]{Corollary}

\newtheorem{lemma}{Lemma}[section]
\newtheorem{proposition}{Proposition}[section]

\allowdisplaybreaks[4]

\numberwithin{equation}{section}
\title[Modulational stability of finite-depth Stokes waves with constant vorticity]{Modulational stability of finite-depth Stokes waves with constant vorticity}
\begin{document}

\author{Lili Fan}
\address[Lili Fan]{School of Mathematics and Statistics,
Henan Normal University, Xinxiang 453007, China}
\email{fanlili@htu.edu.cn}
\author{Changzhen Sun$^*$}
\address[Changzhen Sun]{Universit\'e  Marie et Louis-Pasteur, CNRS, Laboratoire de Math\'ematiques de Besan\c{c}on (UMR CNRS 6623), F-25000 Besançon, France.}
\email{changzhen.sun@cnrs.fr, (corresponding author)}
\date{}

\begin{abstract}
We investigate the effect of constant vorticity on the modulational spectral stability of small-amplitude periodic gravity waves for the two-dimensional finite-depth water-wave system. We derive an explicit stability criterion that agrees with the prediction obtained by Thomas, Kharif and Manna \cite{TKM2012} through a nonlinear Schr\"odinger approximation. 
For waves propagating to the right, we identify a sharp range of dimensionless positive (counterclockwise) vorticity, for which modulational instability is suppressed at every depth, while at any fixed depth sufficiently strong negative (clockwise) vorticity also restores modulational stability.
We give a complete description of the spectrum near the origin and identify a figure-eight spectral curve in the unstable regime. A key distinction from the irrotational case is that the usual Hamiltonian factorization of the  linearized operator fails at zero Floquet exponent on the full periodic Sobolev space. We therefore carry out a spectral reduction based on reversibility, without relying on Hamiltonian factorization.

\vspace{0.5em}


\vspace{0.5em}


\end{abstract}

\maketitle

\setcounter{tocdepth}{1}
\tableofcontents

\section{Introduction}
In this paper, we study the spectral stability of periodic waves of the two-dimensional water-wave
system with vorticity in a channel with finite depth. The equations of motion are given by
\begin{equation}\label{1.1}
\left\{\begin{array}{l}
\partial_t U+U \cdot \nabla U+\nabla P+g e_y=0, \\
\operatorname{div} U=0,\quad \operatorname{curl} U=\gamma.
\end{array} \quad( x, y) \in \Omega_t\,.\right.
\end{equation}
Here $\gamma\in\mr$ is the vorticity, and the fluid domain $\Omega_t$, which is also unknown, is defined as
\begin{equation*}
\Omega_t:=\left\{(x, y) \in \mathbb{R}^2: -h<y<\eta(t, x)\right\},
\end{equation*}
where $h>0$ is a fixed constant and $\eta(t, x)>-h$ is the free surface, which satisfies the kinematic condition
\begin{equation}\label{1.2}
\partial_t \eta(t, x)+\partial_x \eta(t, x) \cdot U_1(t, x, \eta(t,x))-U_2(t, x, \eta(t,x))=0.
\end{equation}
The system \eqref{1.1}-\eqref{1.2} is closed with the boundary conditions for $U_2$ at the lower boundary and for the pressure $P$ at the upper boundary:
\begin{equation}\label{1.3}
\left.U_2\right|_{y=-h}=0,  \;\;\text {and}\;\;
\left.P\right|_{y=\eta}=P_{\text{atm}},
\end{equation}
where $P_{\text{atm}}$ is a constant and can be chosen to be zero.

Since $(U_1+\gamma y,U_2)^t$ is curl-free, we can find a velocity potential $\Phi$ such that
\begin{equation}\label{1.5}
\Phi_{x}=U_1+\gamma y \;\; \text { and } \;\; \Phi_{y}=U_2\,.
\end{equation}
Let $\Psi$ be the stream function such that 
\begin{align*}
    \Psi_y=\Phi_x, \quad \Psi_x=-\Phi_y, \quad \Psi|_{y=-h}=0\,.
\end{align*}
Passing to the frame moving with speed $c>0,$  equations \eqref{1.1}-\eqref{1.3} reduce to 
\begin{equation}
\begin{cases}
\left(\partial_x^2+\partial_y^2\right) \Phi=0, & \text { in } \Omega_t, \\
\partial_t \Phi-c\partial_x \Phi+\frac{1}{2}|\nabla\Phi|^2+g y-\gamma y \p_x \Phi+\gamma\Psi={\rm const}, & \text { on } y=\eta, \\
\partial_t \eta-c\partial_x \eta+\partial_x \eta \partial_x \Phi-\gamma \eta\partial_x \eta-\partial_y \Phi=0, & \text { on } y=\eta, \\
\partial_y \Phi=0, & \text { on } y=-h.\end{cases}
\end{equation}
To make the equations non-dimensional, we perform the change of variables
\begin{equation}\label{1.13}
\tilde{t}=c t / h, \quad \tilde{x}=x / h, \quad \tilde{y}=y / h, \quad \tilde{\eta}=\eta / h, \quad \tilde{\Psi}=\Psi /(c h), \quad \tilde{\Phi}=\Phi /(c h), \quad \tilde{\gamma}=\gamma h /c.
\end{equation}
The new equations then read
\begin{equation}\label{eq-1}
\begin{cases}\left(\partial_x^2+\partial_y^2\right) \Phi=0, & \text { in } -1<y<\eta, \\
\partial_t \Phi-\partial_x \Phi+\frac{1}{2}|\nabla\Phi|^2+\alpha\eta-\gamma\eta\p_x\Phi+
\gamma\Psi={\rm const}, & \text { on } y=\eta, \\
\partial_t \eta-\partial_x \eta+\partial_x \eta \partial_x \Phi-\gamma \eta\partial_x \eta-\partial_y \Phi=0, & \text { on } y=\eta, \\
\partial_y \Phi=0, & \text { on } y=-1, \end{cases}
\end{equation}
where $\alpha=\frac{g h}{c^2}$
is the inverse square of the Froude number, and we have dropped the tildes for notational convenience.

Let $\varphi$ be the trace of the velocity potential on the surface $\eta$ and define the Dirichlet-Neumann (DN) operator $G[\eta]$ by
$$
G[\eta] \varphi=\left.\left(\partial_y \Phi-\partial_x \eta \cdot \partial_x \Phi\right)\right|_{y=\eta},
$$
where $\Phi$ is harmonic with $\left.\Phi\right|_{y=\eta}=\varphi$ and $\left.\partial_y \Phi\right|_{y=-1}=0$. 
Moreover, with the stream function normalized to vanish on the bottom boundary, we have 
\begin{align}\label{def-Heta}
 \Psi(x,\eta(x))=\int_{-1}^{\eta} \p_x \Phi (x,y) \,\d y := H[\eta]\varphi.
\end{align}
It is straightforward to check that 
\begin{equation}\label{1.9}
-\partial_x  H[\eta]\varphi=G[\eta]\varphi.
\end{equation}
We remark that the relation \eqref{1.9} may be regarded as determining $H[\eta]\varphi$ only up to an additive constant, which, however, is fixed by the definition \eqref{def-Heta}.

With these preparations, the system \eqref{eq-1} can be rewritten as a closed system for the surface elevation $\eta$ and the trace of the velocity potential $\varphi:$
\begin{equation}\label{1.8}
\left\{\begin{array}{l}
\partial_t \eta=\partial_x (\eta+\frac {\gamma}{2}\eta^2)+G[\eta] \varphi, \\[3pt]
\partial_t \varphi=\partial_x \varphi-\frac{1}{2}\left|\partial_x \varphi\right|^2+\frac{1}{2} \frac{\left(G[\eta] \varphi+\partial_x \varphi \cdot \partial_x \eta\right)^2}{1+\left|\partial_x \eta\right|^2}-\alpha \eta+\gamma \eta \partial_x \varphi-\gamma H[\eta]\varphi+\mathrm{const}.
\end{array}\right.
\end{equation}
For each fixed background wave, the temporal gauge of the potential is chosen so that the Bernoulli term is constant in time, and this constant is held fixed later in the linearization.
When $\gamma=0$, the above system reduces to the pure gravity water-wave problem, which admits a family of periodic stationary solutions for $\alpha>1$, commonly known as Stokes waves. The stability of Stokes waves has attracted significant interest in recent years. When $\gamma\neq 0, \alpha>0,$ such that $ \alpha+\gamma>1,$ there exists a family of periodic small-amplitude waves
 $\left(\tilde{\eta}_{\varepsilon}, \tilde{\varphi}_{\varepsilon}\right)$  parameterized by the amplitude $\varepsilon$,
$$
\left(\tilde{\eta}_{\varepsilon}, \tilde{\varphi}_{\varepsilon}\right)(x)=\left(\eta_{\varepsilon}, \varphi_{\varepsilon}\right)\left(k_{\varepsilon} x\right),
$$
where $\left(\eta_{\varepsilon}, \varphi_{\varepsilon}\right)$ are $2 \pi$-periodic functions, $k_{\varepsilon}$ is the wave number and satisfies the expansion $k_{\vep}=\kappa+\cO(\vep^2)$ with $\frac{\kappa}{\tah (\kappa)}=\alpha+\gamma.$
For details on the expansions of wave profiles $\left(\eta_{\varepsilon}, \varphi_{\varepsilon}\right)$, see Proposition \ref{prop-expan-profile}.  Studying the stability of $\left(\tilde{\eta}_{\varepsilon}, \tilde{\varphi}_{\varepsilon}\right)$ for the system \eqref{1.8} is equivalent to studying the stability of $\left(\eta_{\varepsilon}, \varphi_{\varepsilon}\right)$ for the system obtained by the scaling $x \rightarrow k_{\varepsilon} x$. The scaled DN operator $G\left[\eta_{\varepsilon}\right]$ is defined as $G\left[\eta_{\varepsilon}\right] \varphi:=\left.\left(\partial_y \Phi-k_{\varepsilon}^2 \partial_x \Phi \partial_x \eta_{\varepsilon}\right)\right|_{y=\eta_{\varepsilon}}$, where $\Phi$ solves the elliptic problem
\beqs
\left(k_{\varepsilon}^2 \partial_x^2+\partial_y^2\right) \Phi=0 \quad \text { in } -1<y<\eta_{\varepsilon} ;\left.\quad \Phi\right|_{y=\eta_{\varepsilon}}=\varphi,\left.\quad \partial_y \Phi\right|_{y=-1}=0.
\eeqs
Using the shape derivative of the DN operator \cite{Lannes}: 
\begin{align*}
\frac{\d}{\d s} G[\eta_{\vep}+s {\eta}](\vp_{\vep}+s{\vp})\big|_{s=0}=G[\eta_{\vep}](\vp-Z_{\vep}\eta)-k_{\vep}\p_x (v_{\vep}\eta)
\end{align*}
where
\begin{equation}\label{1.10}
Z_{\varepsilon}:=\frac{G\left[\eta_{\varepsilon}\right] \varphi_{\varepsilon}+k_{\varepsilon}^2 \partial_x \eta_{\varepsilon} \partial_x \varphi_{\varepsilon}}{1+\left|k_{\varepsilon} \partial_x \eta_{\varepsilon}\right|^2}, \quad \quad v_{\varepsilon}:=k_{\varepsilon}\left(\partial_x \varphi_{\varepsilon}-Z_{\varepsilon} \partial_x \eta_{\varepsilon}\right),
\end{equation}
we can compute that 
\begin{align*}
\frac{\d}{\d s} H[\eta_{\vep}+s {\eta}](\vp_{\vep}+s{\vp})\big|_{s=0}=H[\eta_{\vep}](\vp-Z_{\vep}\eta)+v_{\vep}\eta\,.
\end{align*}
Linearizing the scaled system around $\left(\eta_{\varepsilon}, \varphi_{\varepsilon}\right),$ we find the linearized system
\beq
\left\{\begin{array}{l}
\partial_t \eta=k_{\varepsilon} \partial_x \eta+G\left[\eta_{\varepsilon}\right]\left(\varphi-Z_{\varepsilon} \eta\right)-k_{\varepsilon} \partial_x\left(v_{\varepsilon} \eta\right)+\gamma k_{\varepsilon} \partial_x(\eta_{\varepsilon} \eta), \\[3pt]
\partial_t \varphi=\left(1-v_{\varepsilon}\right) k_{\varepsilon} \partial_x \varphi+Z_{\varepsilon} G\left[\eta_{\varepsilon}\right]\left(\varphi-Z_{\varepsilon} \eta\right)-\left(\alpha+Z_{\varepsilon} k_{\varepsilon} \partial_x v_{\varepsilon}\right) \eta\\[3pt]
\qquad\quad+\gamma \eta_{\varepsilon} k_{\varepsilon} \partial_x \varphi+\gamma k_{\varepsilon} \partial_x \varphi_{\varepsilon} \eta-\gamma v_{\varepsilon} \eta -\gamma 
H\left[\eta_{\varepsilon}\right]\left(\varphi-Z_{\varepsilon} \eta\right).
\end{array}\right.
\eeq
Define the new unknowns $V_1=\eta, V_2=\varphi-Z_{\varepsilon} \eta$. Then the linearized system becomes
\begin{equation}\label{1.12}
\partial_t V=L^{\varepsilon} V, \quad L^{\varepsilon}=\left(\begin{array}{cc}
k_{\varepsilon} \partial_x\left(d_{\varepsilon} \cdot\right) & G\left[\eta_{\varepsilon}\right] \cdot \\
-w_{\varepsilon} & k_{\varepsilon} d_{\varepsilon} \partial_x-\gamma  H \left[\eta_{\varepsilon}\right]
\end{array}\right),
\end{equation}
where
$$
d_{\varepsilon}:=1-v_{\varepsilon}+\gamma \eta_{\varepsilon}, \quad w_{\varepsilon}:=\alpha-d_{\varepsilon} k_{\varepsilon} \partial_x Z_{\varepsilon}\,.
$$
Here we again suppress the dependence of $L^{\varepsilon}$ on the parameters $(\alpha, \gamma)$ for simplicity.

We are interested in the modulational stability of the periodic waves for the system \eqref{1.8} under localized long wave perturbations. This is equivalent to studying the spectrum of the linear operator $L^{\varepsilon}$ near the origin of the complex plane in the (localized) space
$$
Y=L^2(\mathbb{R}) \times H_*^{\frac{1}{2}}(\mathbb{R}), \quad\left(H_*^{s+\frac{1}{2}}(\mathbb{R}):=\sqrt{G[0]}^{-1} H^s(\mathbb{R}),\, s\geq 0\right)
$$
with the domain $H^{1}(\mathbb{R}) \times H_*^{\frac{3}{2}}(\mathbb{R})$.
\subsection{Main results} 

The main theoretical result of this work is a criterion characterizing modulational stability and instability of periodic Stokes waves with constant vorticity.  To state our main result, define the Bloch operator $L_{\xi}^{\varepsilon}:=e^{-\i x\xi}L^{\vep}e^{\i x \xi}=L^{\varepsilon}\left(\partial_x+\i \xi\right)$ which acts on periodic function space 
$$Y_{\text {per }}=L^2\left(\mathbb{T}_{2 \pi}\right) \times H^{\frac{1}{2}}\left(\mathbb{T}_{2 \pi}\right)$$ with domain $D_{\per}=H^1(\mathbb{T}_{2 \pi})\times H^{\frac32}\left(\mathbb{T}_{2 \pi}\right)$. The modulational stability concerns the spectrum of $L_{\xi}^{\vep}$ near the origin for sufficiently small $|\xi|.$

\begin{theorem}\label{thm-main}
Let $\alpha>0, \kappa>0$  and $\gamma\in \mR$ be such that 
$\alpha+\gamma=\frac{\kappa}{\tah(\kappa)}.$
There exist $\vep_0,\xi_0,r_0>0$ such that, for
$|\vep|\leq\vep_0$ and $|\xi|\leq\xi_0$, the spectrum of
$L_\xi^\vep$ in $B_{r_0}(0)$ consists of four eigenvalues,
counted with algebraic multiplicity. For $\xi\neq0$, these are
\beqs
\begin{aligned}
\lambda_1^\pm(\vep,\xi)
 =\i\kappa\xi\left(c_*(\vep,\xi)  \pm\frac12\sqrt{\Delta(\vep,\xi)}\right), \qquad
\lambda_2^\pm(\vep,\xi) =\i\kappa\xi\,\widetilde X_2^\pm(\vep,\xi),
\end{aligned}
\eeqs
where $c_*,\Delta,\widetilde X_2^\pm$ are smooth real-valued functions satisfying 
\beqs
\begin{aligned}
c_*(\vep,\xi)&=1-\nu_*+O(\vep^2+\xi^2),\\
\widetilde X_2^\pm(\vep,\xi)
 &=1-\frac\gamma2\pm\frac12\sqrt{\gamma^2+4\alpha}
   +O(|\vep|+\xi^2),\\
\Delta(\vep,\xi)
 &=4\delta_*\left(\delta_*\kappa^2\xi^2
        +\mathcal C_*\vep^2\right)
   +O\big((\xi^2+\vep^2)(\xi^2+|\vep|)\big),
\end{aligned}
\eeqs
and $\nu_*<1, \delta_*>0$, $\mathcal C_*$ is defined in \eqref{def-cC}.
At $\xi=0$, all four eigenvalues are zero. The following alternatives hold.
\begin{itemize}
\item If $\mathcal C_*>0$, then
$\sigma(L_\xi^\vep)\cap B_{r_0}(0)\subset\i\mR$
for $|\vep|\leq\vep_0$, $|\xi|\leq\xi_0$.
\item If $\mathcal C_*<0$, then, for each $0<|\vep|<\vep_0$,
there exists a unique $\xi_*(\vep)=\frac{\sqrt{-\mathcal C_*}}{\kappa\sqrt{\delta_*}}
|\vep|+O(\vep^2)\in(0,\xi_0)$ such that
for $0<|\xi|<\xi_*(\vep)$, two eigenvalues are purely imaginary
and the other two have opposite nonzero real parts.
For $\xi_*(\vep)<|\xi|\leq\xi_0$, all four eigenvalues
are purely imaginary. 
\end{itemize}
\end{theorem}

\begin{rmk}
In the small-amplitude modulational regime, gravity waves with constant vorticity can be formally described by a nonlinear
Schr\"odinger (NLS) equation governing the slowly varying envelope of a periodic carrier wave \cite{TKM2012}.  The criterion in Theorem~\ref{thm-main} agrees with the prediction of this approximation. We refer to Section~\ref{app:TKM-comparison} for the detailed comparison.
\end{rmk}

\begin{rmk}
\begin{figure}[htbp]\label{Figure-8}
    \centering
    \includegraphics[
        width=0.7\linewidth,
        height=0.7\textheight,
        keepaspectratio
    ]{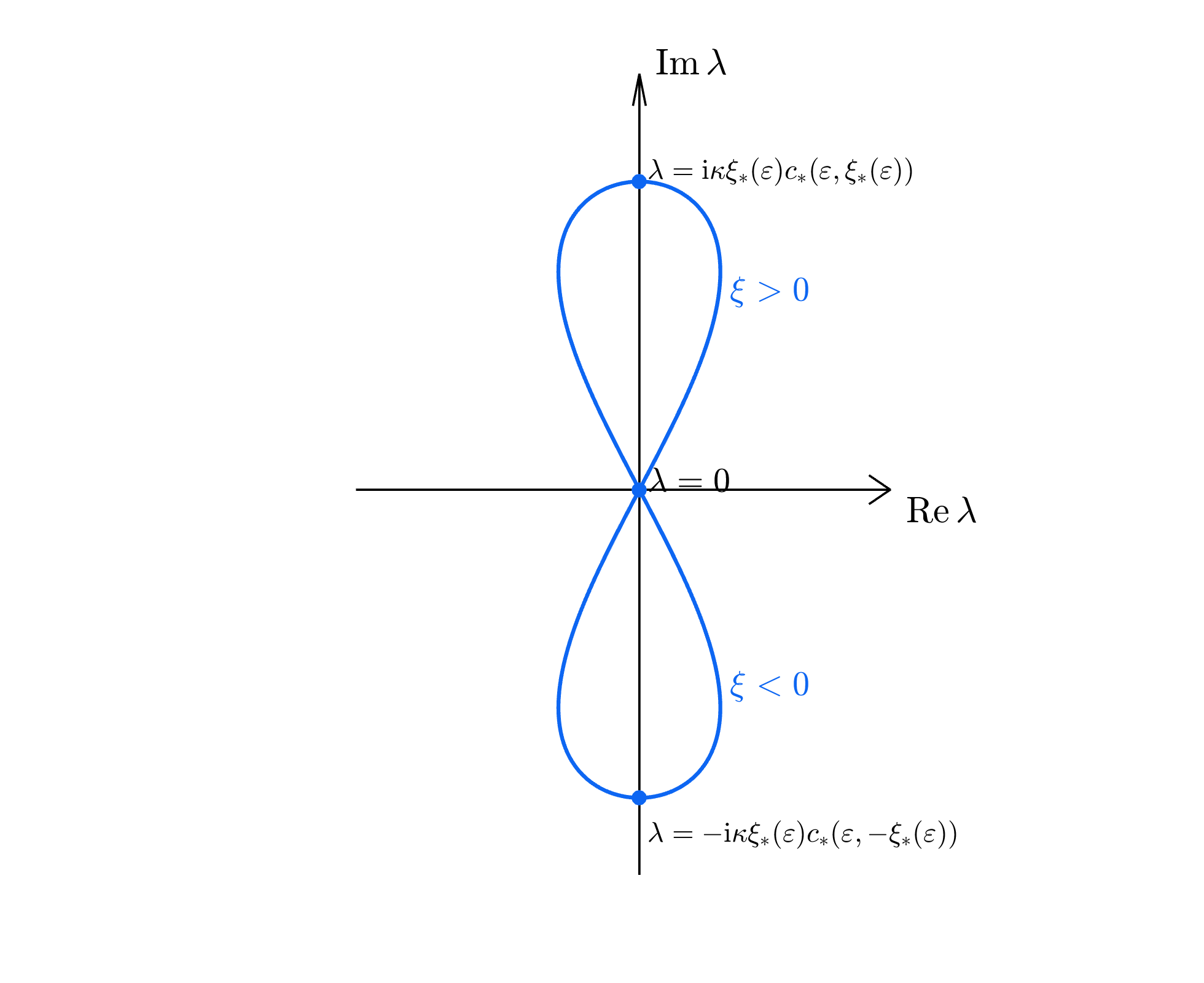}
    \caption{Illustration of the figure-eight spectral curves traced out by $\lambda_1^{\pm}(\vep,\xi)$ for fixed, sufficiently small $\vep$, as $\xi$ varies from $-\xi_{*}(\vep)$ to $\xi_{*}(\vep)$. The upper and lower loops correspond to positive and negative Floquet exponents, respectively. For $\xi>0$, the right and left branches correspond to
$\lambda_1^-(\vep,\xi)$ and $\lambda_1^+(\vep,\xi)$,
respectively; this correspondence is reversed for $\xi<0$.}
    \label{fig:Figure8}
\end{figure}

Theorem \ref{thm-main} provides asymptotic formulas for all four
eigenvalues of $L_\xi^\vep$ near the origin. When $\mathcal C_*<0$, the two eigenvalue branches
$\lambda_1^\pm(\vep,\xi)$ trace out a figure-eight-shaped curve
in the complex plane as $\xi$ varies from $-\xi_*(\vep)$ to $\xi_*(\vep)$, as illustrated in Figure \ref{fig:Figure8}. 
\end{rmk}

\begin{rmk}
In the irrotational case $\gamma=0$, our stability
index recovers the classical finite-depth
Benjamin--Feir criterion. Indeed, direct algebraic
computations show that
\beqs
\mathcal C_*(\kappa,0)
=A(\kappa)\mathcal C_{\rm BF}(\kappa),
\qquad
A(\kappa)=
\frac{\kappa^3}
{32\sigma^3(1-\sigma^2)
(\kappa-\sigma\nu_*^2)}>0,
\eeqs
where $\sigma=\tah\kappa$ and
\beqs
\begin{aligned}
\cC_{\rm BF}(\kappa)
={}&\kappa^2\left(9-28\sigma^2+38\sigma^4-28\sigma^6+9\sigma^8\right)\\
&+2\kappa\sigma\left(3\sigma^6-23\sigma^4+13\sigma^2-9\right)
+\sigma^2\left(9+38\sigma^2-7\sigma^4\right).
\end{aligned}
\eeqs
The expression $\cC_{\rm BF}$ is equivalent to
the classical Benjamin--Feir index
\cite{BM-BFins,BMP-ARMA}, which changes sign
at the critical depth
$\kappa_{\rm BF}\approx1.3627827567$.
Consequently, our criterion recovers the modulational
stability of small-amplitude gravity waves for
$0<\kappa<\kappa_{\rm BF}$ and their instability
for $\kappa>\kappa_{\rm BF}$.
\end{rmk}

\subsubsection{Stability diagram}

For convenience, we introduce the  dimensionless vorticity parameter $X=-\gamma\tah(\kappa)/\kappa>-1$ and write the criterion $\mathcal{C}_{*}=\mathcal C_*(\kappa,X).$   
We are able to prove rigorously the following proposition. 

\begin{prop}\label{prop:stability}
The following properties hold

\rm{(1)\,{(Fix $X,$ varying $\kappa$).} }
 If $X\in (-1,-\frac23],$ then for any $\kappa>0,$ we have
$\mathcal C_*(\kappa,X)>0.$
For every fixed $X>-2/3$, there exists $\kappa_X>0$ such that
$\mathcal C_*(\kappa,X)<0$ for all $\kappa>\kappa_X$.

\rm{(2)\,{(Fix $\kappa,$ varying $X$).} } For any fixed $\kappa>0,$ there exists $X_{\kappa}>0$, such that $\mathcal C_*(\kappa,X)>0$ for any $X>X_{\kappa}.$ 

\end{prop}

 The above results clearly illustrate the effects of constant vorticity on the modulational stability of periodic water waves. For waves traveling to the right, on the one hand, a sufficiently strong positive vorticity, corresponding to $X\in(-1,-\f23]$, suppresses modulational instability of small-amplitude waves at every dimensionless depth. On the other hand, at any fixed depth, a sufficiently strong negative vorticity ($X>0$) can also stabilize waves that may otherwise be modulationally unstable in the absence of shear.

\begin{figure}[htbp]
    \centering
    \includegraphics[
        width=0.7\linewidth,
        height=0.6\textheight,
        keepaspectratio
    ]{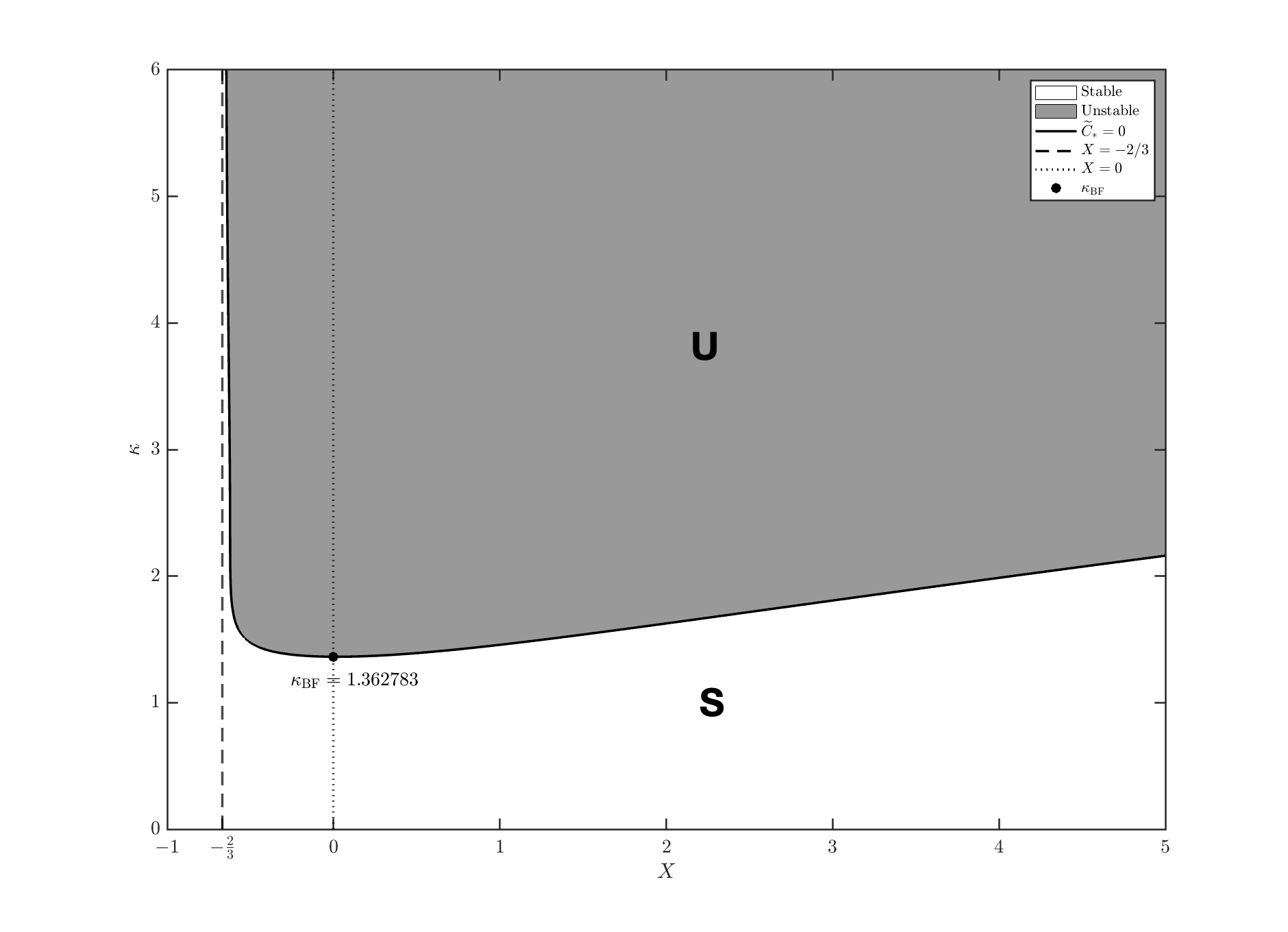}
    \caption{Numerically computed modulational stability diagram in the
    $(X,\kappa)$-plane. The gray region corresponds to modulationally
    unstable waves, while the white region corresponds to
    modulationally stable waves. The stability region when $X\in (-1,-\f23]$ is rigorously proven in Proposition \ref{prop:stability}.}
    \label{fig:stability-diagram}
\end{figure}

\textbf{Numerical evaluations.}
Numerical evaluation of the index $\mathcal C_*(\kappa,X)$
gives the stability diagram displayed in
Figure~\ref{fig:stability-diagram}. 
For  $X>-2/3$, the numerical results in Figure~\ref{fig:stability-diagram} suggest that $\mathcal C_*(\kappa,X)$ vanishes at a unique depth $\kappa=\kappa_{\rm crit}(X)$, with modulationally stable waves below the curve and unstable waves above it. The curve passes through the classical Benjamin--Feir point at $X=0$ and tends to infinity as $X\to-2/3^+$.

An interesting feature is that the critical depth does not vary monotonically with $X$. It decreases to a minimum at $X_{\min}\approx0.0196712$ and then increases. The minimum value, $\kappa_{\rm crit}(X_{\min})\approx1.3627149$, is only about $7\times10^{-5}$ below the irrotational threshold $\kappa_{\rm BF}\approx1.3627828$. Thus, increasing $X$ reduces the range of stable depths when $-2/3<X<X_{\min}$ and enlarges it when $X>X_{\min}$. Compared with the irrotational case, the range is slightly smaller only for $0<X<X_{\rm BF}$, where $X_{\rm BF}\approx0.0396727$. These observations are based on numerical computations.

\subsection{Background on the stability of periodic water waves}

In the 1960s, Benjamin and Feir \cite{benjamin1967disintegration} observed experimentally and through formal perturbation arguments that periodic gravity waves in deep water are unstable with respect to long-wave sideband perturbations. Since then, the stability of periodic traveling solutions of the water-wave system has been an important topic in water-wave theory. A rigorous treatment in finite depth was later given by Bridges and Mielke \cite{BM-BFins}, who used a Hamiltonian spatial-dynamics formulation together with a center-manifold reduction. The deep-water case, where this reduction is no longer available, was addressed more recently by Nguyen and Strauss \cite{Nguyen-Strauss}. Hur and Yang \cite{Hur-Yang-ARMA} subsequently revisited the finite-depth problem by means of an Evans-function analysis. Besides recovering the modulational instability criterion, they also found instability near the next nonzero spectral crossing. A different approach, based on Bloch--Floquet decomposition and Kato's analytic perturbation theory \cite{Kato}, was developed by Berti and collaborators in a series of works \cite{BMP-ARMA,BMP-deep,berti-highfreq}. These papers provide a detailed description of the spectrum near the origin and at higher-frequency crossings. Earlier numerical and formal studies of high-frequency and superharmonic instabilities include \cite{Longuet-Higgins1978,DeconinckOliveras2011,CDT-highfreq}. The nonlinear modulational instability of small-amplitude deep-water Stokes waves was proved in \cite{Chen-Su} using the NLS approximation of the water-wave system. Recent developments also include the stability analysis of gravity-capillary periodic waves \cite{Hur-Yang-capillary,Mariana-Tien-Erik,Hsiao-Maspero}, the study of planar Stokes waves under purely transverse perturbations \cite{CNS-deep,CNS-finite}, and the analysis of perturbations involving both longitudinal and transverse frequencies \cite{JRSY-cmp,BMR-Mclean}.

The results above concern the classical irrotational water-wave problem. In many physical situations, however, the underlying current is vertically sheared and the fluid has nonzero vorticity, for instance due to wind forcing, tidal motion or interaction with the bottom. The case of constant vorticity corresponds to a linear shear current and provides a natural first model for studying how rotational effects modify the dispersion relation, the wave profile and the nonlinear interaction between the carrier wave and the induced mean flow. The existence of periodic traveling waves with nontrivial vorticity has been studied extensively; see, for instance, \cite{SimmenSaffman1985,SilvaPeregrine2018,Constaintin-Strauss-CPAM,CSV-acta,Erik-Weber,Kozlov-Lokharu}.

Modulational instability in the presence of vorticity has also been studied through full-dispersion shallow-water models. Hur and Johnson \cite{HurJohnson2015} considered a Whitham equation whose linear dispersion relation accounts for constant vorticity and surface tension. They derived a modulational stability criterion for small-amplitude periodic waves and investigated the stabilizing effect of vorticity. Related higher-order vortical Whitham models were considered in \cite{CKKA2022}. These works provide complementary descriptions of the effect of constant vorticity on modulational stability through different Whitham-type equations.

In the long-wave modulation regime, the stability of small-amplitude periodic waves is closely related to the focusing or defocusing character of the NLS equation governing their slow modulations. For constant-vorticity flows, Li, Hui and Donelan \cite{LiHuiDonelan1987} derived an NLS equation for deep-water waves (see also the subsequent work of Baumstein \cite{Baumstein_deep}). Thomas, Kharif and Manna \cite{TKM2012} carried out the analysis in finite depth and derived explicit expressions for the NLS coefficients in terms of the vorticity and dimensionless depth. Their criterion shows that vorticity may substantially modify the growth rate and bandwidth of the modulational instability, and that the nonlinear coupling between the modulation-induced mean flow and the background vorticity plays an essential role. In particular, sufficiently strong shear can suppress the Benjamin--Feir instability, while for a suitable adverse current the waves are predicted to be modulationally stable for all depths.

Despite these developments, rigorous results on the spectral stability of periodic water waves with vorticity remain limited. Constantin and Strauss \cite{CS-CPAM} investigated variational stability properties of exact periodic traveling water waves with vorticity, including a formal stability result for small-amplitude waves with wavelengths large relative to the fluid depth. Hur and Lin \cite{HL-CMP2008,HL-CMP2013} established the linear instability of small-amplitude periodic waves bifurcating from certain unstable shear flows with nonconstant vorticity. These works address stability questions distinct from the modulational spectral problem considered here. 

Very recently, Chen and Nguyen \cite{chen2026suppression} studied the Benjamin--Feir instability of small-amplitude Stokes waves in infinite depth.
They proved that the waves are modulationally unstable below a critical vorticity and
 stable above it, confirming the
transition predicted in \cite{TKM2012}. 
Their analysis applies Wahl\'en's change of
variables \cite{Erik-LMP}, which puts the linearized problem in canonical Hamiltonian
form for nonzero Floquet exponents, and handles the resulting singular behavior of the
vorticity term on the zero Fourier mode as the Floquet exponent tends to zero.
 In contrast, the present paper treats finite depth and derives a stability criterion depending on both depth and vorticity through a spectral reduction based on reversibility, without relying on Hamiltonian factorization at zero Floquet exponent. The two works therefore provide complementary results in distinct depth regimes, each involving different challenges and proof strategies.
 
During the preparation of this work, we learned that Hsiao, Hur and Yang \cite{HHY-finite} independently investigated the spectral instability of small-amplitude periodic gravity waves with constant vorticity in finite depth, both near the origin and at the next nonzero spectral crossing, using an Evans-function approach which differs substantially from the spectral reduction developed here. They also investigated the corresponding infinite-depth problem using Kato's perturbation theory \cite{HHY-deep}. 

\subsection{Main contributions and strategy of proof}

The present work has two complementary contributions.
First, we obtain a complete description of the modulational spectrum of small-amplitude periodic gravity waves with
constant vorticity in finite depth. The resulting explicit
stability index depends on both the dimensionless depth
$\kappa$ and the vorticity parameter $X$. In particular, we identify the sharp range $-1<X\leq-2/3$ for which modulational stability holds at every dimensionless depth. We also establish a second
stabilizing regime for sufficiently large positive $X$ at any fixed depth. These results provide a rigorous description of the stabilizing effects of constant
vorticity in the finite-depth setting.

Second, we develop a spectral reduction that applies to the constant-vorticity problem. 
Let us recall that, from the spectral point of view, modulational (in)-stability is determined
by the behavior of the spectrum of the Bloch operators $L_\xi^\vep$
near the origin when both the wave amplitude $\vep$ and the Floquet
parameter $\xi$ are small.
The proof is based on 
Kato's perturbation theory, a framework used in many studies of periodic traveling waves \cite{GH-NLS,Maspero-Kdv,Deconinck-KdV,ARS-KdV,NRS-EP,CFWX-KP}.
For Stokes waves, Berti, Maspero and Ventura \cite{BMP-deep,BMP-ARMA} extend a basis of the unperturbed generalized kernel jointly in the Floquet exponent and the amplitude by Kato's similarity transformation, and analyze the resulting $4\times4$ Hamiltonian and reversible matrix through a sequence of block-decoupling transformations inspired by KAM theory. In \cite{SW-capillary}, the second author and E.~Wahl\'en proposed a different reduction: the basis is obtained by extending, in the Floquet exponent only, a basis of the generalized kernel of $L_0^{\vep}$ built exactly from the family of traveling waves. With this choice, the representation matrix inherits the Jordan structure of $L_0^{\vep}$, so that, after a rescaling, its eigenvalues can be analyzed by a single block reduction. The construction in \cite{SW-capillary} still relies on the Hamiltonian structure, both to construct the dual basis and to derive the expansions of the matrix entries. With constant vorticity, however, the Hamiltonian
factorization is not available for the linearized
operator at zero Floquet exponent (see below
and Appendix~\ref{appen-Hamiltonian-0}).
We overcome this difficulty by constructing
the dual basis directly in the generalized kernel
of $(L_0^\vep)^*$. Reversibility, together with the smooth dependence
on the Floquet exponent, then provides the structure needed for the spectral reduction,
without relying on the Hamiltonian formulation at $\xi=0$. 
This approach may also be useful for other reversible problems where the Hamiltonian structure cannot be directly exploited at zero Floquet exponent.

We now describe two specific difficulties
arising from constant vorticity in the
spectral reduction. The first difficulty concerns the construction of a generalized eigenbasis for the zero eigenvalue of $L_0^{\vep}$. In the irrotational case,
differentiating the traveling-wave equations with respect to the wave parameters directly yields
the desired Jordan relation. In the presence of
constant vorticity, however, the variation of the
Bernoulli constant introduces an additional term.
More precisely, combining derivatives of the
traveling-wave profile with respect to the amplitude
and wave speed, we obtain a vector
$\widehat q_2^\varepsilon$ satisfying
$$
L_0^\varepsilon\widehat q_2^\varepsilon
=a_\varepsilon q_1^\varepsilon+r_\varepsilon e_2,
\qquad e_2=(0,1)^t.
$$
Although both $q_1^\varepsilon$ and $e_2$ belong
to the kernel of $L_0^{\vep}$, the additional term prevents
$\widehat q_2^\varepsilon$ from satisfying the
desired Jordan relation. To eliminate this term,
we construct, by a Lyapunov--Schmidt argument,
a generalized eigenvector $q_4^\varepsilon$
satisfying
$$
L_0^\varepsilon q_4^\varepsilon
=-e_2-\Lambda_\varepsilon q_1^\varepsilon.
$$
The corrected vector
$q_2^\varepsilon=\widehat q_2^\varepsilon
+r_\varepsilon q_4^\varepsilon$
then satisfies
$L_0^\varepsilon q_2^\varepsilon
=d_{12}(\varepsilon,0)q_1^\varepsilon$.
This choice of basis is particularly useful for
the subsequent spectral reduction, since the correction enters the coefficients determining modulational stability.

The second difficulty concerns the Hamiltonian
structure of the linearized operator at zero Floquet exponent.  For $\xi\neq0$, the Bloch operator admits
a Hamiltonian factorization
$L_\xi^\varepsilon=J_\xi^\varepsilon A_\xi^\varepsilon$ where $J_\xi^\vep$ is skew-adjoint and $A_\xi^\vep$ is self-adjoint.
However, this factorization involves the inverse of
$D_\xi=-\i\partial_x+\xi$, which is not defined on
the full periodic space when $\xi=0$.
In fact, as shown in Appendix \ref{appen-Hamiltonian-0}, the physical linearized operator satisfies
$$
L_0^\varepsilon
=J_0^\varepsilon A_0^\varepsilon+R^\varepsilon,
$$
where  $J_0^\vep$ is defined in \eqref{def-J0ep}  and $R^\varepsilon$ is a nontrivial rank-one correction
taking values in the constant mode, which is proportional to the vorticity and
vanishes when $\gamma=0$.  
This distinction is related to the Hamiltonian
formulation of constant-vorticity water waves
in \cite{Erik-LMP}, where the velocity potential
is considered modulo additive constants.
In the irrotational problem \cite{SW-capillary},  the adjoint
generalized kernel can be recovered from suitable linear combinations
of 
$\big\{J^{-1}q_j^\vep\big\}_{j=1}^4$
and since the Kato extension preserves the corresponding
symplectic pairing, the same relation remains valid after the Floquet
parameter is introduced. In the present problem, the identity
$L_0^\vep=J_0^\vep A_0^\vep$ does not hold for the linearized operator
on the full Sobolev spaces, and this symplectic construction is no
longer available. 
We therefore construct the dual vectors
$\{\tilde q_j^\vep\}_{j=1}^4$ directly in the generalized kernel of
$(L_0^\vep)^*$ by
imposing the biorthogonality conditions
$\ang{\tilde q_j^\vep,q_\ell^\vep}
=\delta_{j\ell}.$
The basis and dual basis must then be extended separately to nonzero
Floquet parameters. We therefore construct a pair of compatible Kato extension operators $\cU^\vep(\xi)$ and $\cV^\vep(\xi)$ such that
$\cV^\vep(\xi)^*\cU^\vep(\xi)=\Id$
and define
\beqs
q_j^\vep(\xi):=\cU^\vep(\xi)q_j^\vep,
\qquad \tilde {q}_j^\vep(\xi)
:=\cV^\vep(\xi)\tilde q_j^\vep.
\eeqs
The biorthogonality relation is thereby preserved for nonzero $\xi.$

Having constructed these bases, we define the representation matrix
\beqs
\rD_\xi^\vep:=\left(
\ang{\tilde q_j^\vep(\xi),
L_\xi^\vep q_\ell^\vep(\xi)}
\right)_{1\leq j,\ell\leq4}.
\eeqs
The key observation is that
reversibility determines the real or purely imaginary character of its entries, without relying on the Hamiltonian structure. Together with the Jordan
relations and the smooth dependence on $\xi$,
this allows us to identify their vanishing orders
and perform the rescaling
\beqs
\rM_\xi^\vep
:=\frac{1}{\i\kappa\xi}
\mathcal P_\xi\rD_\xi^\vep\mathcal P_\xi^{-1} =\bpm  \rM_{11} & \vep \rM_{12} \\
    \vep \rM_{21} & \rM_{22}\epm (\vep,\xi),
\qquad
\mathcal P_\xi
:=\diag(\i\kappa\xi,1,\i\kappa\xi,1).
\eeqs
This rescaling reduces the
problem to the study of the real matrix $\rM_\xi^\vep$ which extends
smoothly to $\xi=0$.
The matrix $\rM_\xi^\vep$ separates at leading order into a carrier-wave block $\rM_{11}$ and a mean-flow block $\rM_{22}$. 
As in the  irrotational case \cite{SW-capillary}, the order-$\vep$ coupling
between the carrier-wave and mean-flow blocks produces
an order-$\vep^2$ Schur-complement correction to the effective carrier block. However, the presence of vorticity modifies the coefficients through the corrected
generalized eigenvector, the dual basis and the
additional terms in the linearized operator.
More precisely, the coefficient governing the splitting
of the carrier modes takes the form
\beqs
\mathcal C_*=\Gamma_*+\mathfrak n_{12},
\qquad
\Gamma_*=
\frac{2\kappa k_2}{c\partial_c\kappa}
-2\mathfrak b_2\Lambda_0,
\qquad
\mathfrak n_{12}
=-\frac{\mathfrak m_{13}\mathfrak m_{42}}
{\nu_*^2-\gamma\nu_*-\alpha}.
\eeqs
Here $\Gamma_*$ is the direct contribution to the
carrier block, while $\mathfrak n_{12}$ arises from
the coupling with the mean-flow block.
In particular, the term $-2\mathfrak b_2\Lambda_0$
arises from the variation of the Bernoulli
constant under our normalization of the
stream function \eqref{def-Heta}.
This correction vanishes in the irrotational
case, and we recover the structure found in
\cite{SW-capillary}.

After the block reduction, the discriminant of the carrier block satisfies
\beqs
\Delta=4\delta_*\left(
\delta_*\kappa^2\xi^2+\mathcal C_*\vep^2\right)
+\cO\big((\xi^2+\vep^2)(|\vep|+\xi^2)\big).
\eeqs
Since $\delta_*>0$, the sign of $\mathcal C_*$
determines the modulational stability or instability
of the carrier modes, while the two mean-flow
eigenvalues remain purely imaginary.
The reduction therefore yields the modulational stability criterion together with asymptotic formulas for all
four eigenvalues of $L_\xi^\vep$ near the origin.

\subsection{Notation}
For simplicity of notation, we use the following conventions throughout the paper:
$$
\operatorname{sh}(\cdot):=\sinh (\cdot), \quad \operatorname{ch}(\cdot):=\cosh (\cdot), \quad \operatorname{th}(\cdot):=\tanh (\cdot).
$$
We sometimes denote $S=\sh(\kappa), C=\ch(\kappa)$ for further simplication.
\section{Spectral preliminaries}\label{sec prelim}
\subsection{Bloch transform}
To analyze the spectrum of operators with periodic coefficients, it is expedient to introduce Bloch symbols, associated with the (Floquet-)Bloch transform.

In a similar spirit to the inverse Fourier transform, the inverse Bloch transform provides a way to reconstruct a function from its Bloch components. For any $f \in L^1(\mathbb{R})$ such that $\hat{f} \in L^1(\mathbb{R})$, it holds that
\begin{equation}\label{2.1}
f(x)=\int_{-1 / 2}^{1 / 2} e^{\i x \xi} \check{f}(\xi, x) \mathrm{d} \xi,
\end{equation}
where the direct Bloch transform of $f, \mathcal{B}(f)=\check{f}$, is defined as
\begin{equation}\label{2.2}
\mathcal{B}(f)(\xi, x)=\check{f}(\xi, x):=\sum_{j \in \mathbb{Z}} e^{\mathrm{i} j x} \hat{f}(\xi+j).
\end{equation}
Note that $\xi$ is called Floquet exponent/parameter.
With the help of the inverse Bloch transform \eqref{2.1}, one can turn an operator with periodic coefficients acting on a function over $\mathbb{R}$ into a family of operators acting on periodic functions. Explicitly, it holds that,
\begin{equation}\label{2.3}
L^{\varepsilon}(u)(x)=\int_{-1 / 2}^{1 / 2} e^{\i x \xi} L_{\xi}^{\varepsilon}(\check{u}(\xi, \cdot))(x) \mathrm{d} \xi,
\end{equation}
where each $L_{\xi}^{\varepsilon}:=L^{\varepsilon}\left(\partial_x+\i \xi\right)$ acts on $Y_{\text {per }}=L^2\left(\mathbb{T}_{2 \pi}\right) \times H^{\frac{1}{2}}\left(\mathbb{T}_{2 \pi}\right)$ with domain $D_{\per}=H^1(\mathbb{T}_{2 \pi})\times H^{\frac32}\left(\mathbb{T}_{2 \pi}\right)$. The main gain when replacing the direct analysis of $L^{\varepsilon}$ with that of $L_{\xi}^{\varepsilon}$ is that each $L_{\xi}^{\varepsilon}$ has compact resolvent, and thus its spectrum consists of eigenvalues of finite multiplicity. A key related observation is that
\beq
\sigma\left(\left.L^{\varepsilon}\right|_Y\right)=\bigcup_{\xi \in(-1 / 2,1 / 2]} \sigma\left(\left.L_{\xi}^{\varepsilon}\right|_{Y_{\rm per}}\right)
\eeq
We may also regard $L_\xi^\vep$ as an operator on
$\cY_{\per}:=\big(L^2(\mathbb T_{2\pi})\big)^2$
with domain
$\cD_{\per}:=\big(H^1(\mathbb T_{2\pi})\big)^2$, which has also compact resolvent. 
By elliptic
regularity, the eigenfunctions and generalized
eigenfunctions under both realizations are smooth. Consequently, their spectra,
including algebraic multiplicities, coincide
\[
\sigma\left(\left.L_\xi^\vep\right|_{Y_{\per}}\right)
=
\sigma\left(\left.L_\xi^\vep\right|_{\cY_{\per}}\right).
\]
In the following, we study $L_\xi^\vep$ on $\cY_{\per}$,
equipped with the usual $L^2$ inner product.

\subsection{Symmetry of the spectrum}
An important feature of the operator $L_{\xi}^{\vep}$ is that, for $\xi\neq 0$, it admits the Hamiltonian factorization
\beq\label{def-G-xi-esp}
L_{\xi}^{\vep}=
\begin{pmatrix}
\i k_{\vep}D_{\xi}(d_{\vep}\cdot) & G_{\xi}[\eta_{\vep}]\\
-w_{\vep} & \i k_{\vep}d_{\vep}D_{\xi}-\gamma
H_{\xi}[\eta_{\vep}]
\end{pmatrix}=J_{\xi}^{\vep}A_{\xi}^{\vep},
\eeq
where
\beqs
J_{\xi}^{\vep}
:=
\left(\begin{array}{cc}
0&1\\
-1&\gamma(\i k_{\vep}D_{\xi})^{-1}
\end{array}\right),
\qquad
A_{\xi}^{\vep}
:=
\left(\begin{array}{cc}
w_{\vep}+\gamma d_{\vep}&-\i k_{\vep}d_{\vep}D_{\xi}\\
\i k_{\vep}D_{\xi}(d_{\vep}\cdot)&G_{\xi}[\eta_{\vep}]
\end{array}\right),
\eeqs
with
\begin{equation}\label{2.4}
D_{\xi}:=-\i\p_x+\xi,
\qquad
G_{\xi}[\eta_{\vep}]
:=
e^{-\i x\xi}G[\eta_{\vep}]e^{\i x\xi}.
\end{equation}
For every real $\xi\neq 0$, we have $(J_{\xi}^{\vep})^*=-J_{\xi}^{\vep},(A_{\xi}^{\vep})^*=A_{\xi}^{\vep}.$ 
Indeed, $D_{\xi}$ and $G_{\xi}[\eta_{\vep}]$ are self-adjoint, and
$\left(-\i k_{\vep}d_{\vep}D_{\xi}\right)^*
=\i k_{\vep}D_{\xi}(d_{\vep}\cdot).$
It follows that $(L_{\xi}^{\vep})^*=-(J_{\xi}^{\vep})^{-1}L_{\xi}^{\vep}J_{\xi}^{\vep}$, for $\xi\neq 0.$ Consequently,
\beq \label{spec-symmetry}
\sigma(L_{\xi}^{\vep})
=-\overline{\sigma(L_{\xi}^{\vep})},
\qquad \xi\neq 0,
\eeq
and hence the spectrum of $L_{\xi}^{\vep}$ is symmetric with respect to the imaginary axis.

However, we need some caution when considering the situation when $\xi=0.$
In this case,
$\i k_{\vep}D_0=k_{\vep}\p_x$
is not invertible on the full periodic functions. Using its zero-mean inverse as in Appendix \ref{app B} defines 
$J_{0}^{\vep}$ but the linearized operator
satisfies   $L_0^{\vep}=J_0^{\vep}A_0^{\vep}+R^{\vep}$ (see Appendix \ref{appen-Hamiltonian-0}) with a correction taking values in the constant potential mode. The correction vanishes after passage to the quotient by constant potential.
Nevertheless, although the Hamiltonian factorization is not directly available for the  operator
$L_0^\vep$, its spectrum is still symmetric with respect to the imaginary axis. This
follows by continuity from the Bloch operators $L_\xi^\vep$ with $\xi\neq 0$. 
\footnote{\,It also follows from the reversibility of the operator $L_{\xi}^{\vep}$. Indeed, let $\cI$ denote the anti-linear reversibility operator (see \eqref{def-revisibility}), and
$$\cI L_{\xi}^{\vep}=-L_{\xi}^{\vep}\cI,$$
then $$L_{\xi}^{\vep}q=\lambda q\quad\Longrightarrow\quad L_{\xi}^{\vep}(\cI q)
= -\overline{\lambda}\,\cI q.
$$
It follows that
$$
\lambda\in\sigma(L_{\xi}^{\vep})
\quad\Longrightarrow\quad
-\overline{\lambda}\in\sigma(L_{\xi}^{\vep}), \quad \text{ for }\, \xi \in\left(-{1}/{2}, {1}/{2}\right].
$$
}
Indeed, let
$\lambda_0\in\sigma(L_0^\vep). $
Since $L_0^\vep$ has compact resolvent, $\lambda_0$ is an isolated eigenvalue of
finite algebraic multiplicity. Let $\Gamma$ be a sufficiently small circle enclosing $\lambda_0$ and no other point of
$\sigma(L_0^\vep)$. For $\xi$ sufficiently close to zero, define the Riesz
projection
\[
\Pi_\xi:=\frac{1}{2\pi\i}
\int_{\Gamma}(z-L_\xi^\vep)^{-1}\,\d z.
\]
By the continuity of the resolvent,
$\Pi_\xi\rightarrow\Pi_0$ 
as $\xi \rightarrow0.$
Consequently, $\operatorname{rank}\Pi_\xi=\operatorname{rank}\Pi_0>0$
for all sufficiently small $\xi$. It follows that there exists a sequence $\xi_n\neq 0$ tending to $0$ and 
$\lambda_{\xi_n}\in\sigma(L_{\xi_n}^\vep)$ such that $\lambda_{\xi_n}\rightarrow \lambda_{0}.$
By \eqref{spec-symmetry}
$$
-\overline{\lambda_{\xi_n}}\in\sigma(L_{\xi_n}^\vep),
\quad
-\overline{\lambda_{\xi_n}}\rightarrow-\overline{\lambda_0},
\quad\text{as }n\rightarrow +\infty.
$$
Again by the continuity of the resolvent of $L_{\xi}^{\vep}$ in $\xi,$ it holds that
$-\overline{\lambda_0}\in\sigma(L_0^\vep).$
Consequently,
\[\sigma(L_0^\vep)=-\overline{\sigma(L_0^\vep)}.
\]
\subsection{Study of the constant operator $L_{\xi}^0$}
By definition,
$$
L_{\xi}^{\varepsilon}=\left(\begin{array}{cc}
\i k_{\vep} D_{\xi} \left(d_{\vep}\cdot\right) & G_{\xi}\left[\eta_{\varepsilon}\right] \cdot \\
-w_{\varepsilon} &\i k_{\varepsilon} d_{\varepsilon}D_{\xi}-\gamma H_{\xi}\left[\eta_{\varepsilon}\right]\cdot
\end{array}\right).
$$
Specifically, when $\varepsilon=0$,
\begin{equation}\label{L0}
L_{\xi}^0=\left(\begin{array}{cc}
\i\kappa D_{\xi}
& \kappa D_{\xi} \operatorname{th}\left(\kappa D_{\xi}\right) \\
-\alpha & \i\kappa D_{\xi}
-\i\gamma \operatorname{th}\left(\kappa D_{\xi}\right)
\end{array}\right).
\end{equation}
Since $L_{\xi}^0$ has constant coefficients, its spectrum can be found explicitly using the Fourier transform. Direct computations yield
\begin{equation}\label{eigenrelation}
L_{\xi}^0 v^\pm_{j}(\xi,x)=\lambda_j^{ \pm}(\xi) v^\pm_{j}(\xi,x), \quad j \in \mathbb{Z},
\end{equation}
where
\begin{equation}\label{lam}
\lambda_j^{ \pm}(\xi)=\i\left(\zeta_j-\frac{\gamma\operatorname{th}(\zeta_j)} 2 \mp  w(\zeta_j)\right), \qquad v_j^\pm(\xi, x) = \begin{pmatrix} 1 \\ \mathrm{ i}\dfrac{-\gamma\operatorname{th}(\zeta_j)\mp 
2w(\zeta_j)}
{2\zeta_j \operatorname{th}(\zeta_j)} \end{pmatrix} e^{\mathrm{i} jx},
\end{equation}
with 
\begin{align}\label{def-wzetaj}
\zeta_j:=\kappa(j+\xi), 
 \qquad
 w(\zeta_j):=\frac {\zeta_j} {2} \sqrt{4\alpha \tah(\zeta_j)/\zeta_j+ (\gamma \tah(\zeta_j)/\zeta_j )^2}\,.
\end{align}

The spectrum of $L_{\xi}^0$ is thus contained in the imaginary axis. By \eqref{spec-symmetry}, each fixed simple eigenvalue of $L_{\xi}^0$ remain on the imaginary axis for small $\varepsilon>0$ and possible unstable spectrum bifurcates only from multiple eigenvalues of $L_{\xi}^0$ when $\varepsilon>0$ is sufficiently small. Therefore, we need to find the crossing spectrum of $L_{\xi}^0$ with $\xi \in\left(-\frac{1}{2}, \frac{1}{2}\right]$.

\begin{lemma}
\label{lem:crossing-origin}
The origin is an eigenvalue of $L^0_0$ of algebraic multiplicity four. More precisely,
\beq\label{crossing-0}
\lambda_{1}^{+}(0)=\lambda_{-1}^{+}(0)=\lambda_{0}^{\pm}(0)=0.
\eeq
Moreover, these are the only eigenvalue curves of $L^0_\xi$ crossing the origin for $\xi\in(-1/2,1/2]$.
\end{lemma}

\begin{proof}
The fact \eqref{crossing-0} is easy to verify and we just prove the second assertion.
Assume that $\zeta_j\neq0$ and $\lambda^\pm_j(\xi)=0$, then 
a further consideration of the dispersion relation \eqref{relation} gives
\beqs
\f{\zeta_j}{\tah(\zeta_j)}=\alpha+\gamma=\f{\kappa}{\tah(\kappa)}.
\eeqs
Since $\mu\mapsto \mu/\tah(\mu)$ is even and strictly increasing on $(0,\infty)$, it follows that
$\zeta_j=\pm\kappa$. Because $\zeta_j=\kappa(j+\xi)$ and $\xi\in(-1/2,1/2]$, the only possibilities are
$(j,\xi)=(\pm 1,0).$
It remains to consider $\zeta_j=0$, which occurs only for $(j,\xi)=(0,0)$. On the zero Fourier mode,
\beqs
L^0_0=
\begin{pmatrix}
0&0\\
-\alpha&0
\end{pmatrix}.
\eeqs
Hence $(0,1)^t$ is an eigenvector and $(1/\alpha,0)^t$ is a generalized eigenvector. Therefore the zero mode contributes algebraic multiplicity two. Together with the two simple crossings from $j=\pm1$, the total algebraic multiplicity at the origin is four. This proves the lemma.
\end{proof}

\section{Basis and dual basis for $L_{\xi}^{\vep}$}
In this section, we aim to construct a basis and its dual for the spectral subspace associated with the spectrum of $L_{\xi}^\vep$ near the origin. To this end, we first construct the basis and dual basis of the generalized kernel of $L_0^{\vep}$ and then extend them to be the basis and dual basis of $L_{\xi}^{\vep}.$

\subsection{Generalized kernel of $L_0^{\vep}$} 
We recall that the operator $L_0^{\varepsilon}= L_{\xi}^{\vep}\big|_{\xi=0}$ is defined as
$$
L_0^{\varepsilon}=\left(\begin{array}{cc}
k_{\varepsilon} \partial_x\left(d_{\varepsilon}{\cdot}\right) & G_0\left[\eta_{\varepsilon}\right] \\
-w_{\varepsilon} & k_{\varepsilon}d_{\varepsilon} \partial_x-\gamma H_0[\eta_{\varepsilon}]
\end{array}\right),
$$
where $G_0\left[\eta_{\varepsilon}\right]$ is defined in \eqref{2.4}. 
The main result is the following.
\begin{proposition}[Generalized kernel of $L_0^{\varepsilon}$]\label{prop-basis}
Zero is an eigenvalue of $L_0^{\varepsilon}$ with algebraic multiplicity four. Accordingly, there exist four (generalized) eigenfunctions $q_1^{\varepsilon}(0, \cdot), q_2^{\varepsilon}(0, \cdot), q_3^{\varepsilon}(0, \cdot), q_4^{\varepsilon}(0, \cdot)$ such that
$$
\begin{array}{ll}
L_0^{\varepsilon} q_1^{\varepsilon}(0, \cdot)=0, & L_0^{\varepsilon} q_2^{\varepsilon}(0, \cdot)=
d_{12}(\vep,0)q_1^{\varepsilon}(0, \cdot), \\
L_0^{\varepsilon} q_3^{\varepsilon}(0, \cdot)=0, & L_0^{\varepsilon} q_4^{\varepsilon}(0, \cdot)=-q_3^{\varepsilon}(0, \cdot),
\end{array}
$$
where 
$$
q_1^{\varepsilon}(0, \cdot)=\varepsilon^{-1}\binom{\partial_x \eta_{\varepsilon}}{\partial_x \varphi_{\varepsilon}-Z_{\varepsilon} \partial_x \eta_{\varepsilon}} 
$$
and $q_3^{\varepsilon}(0, \cdot), q_4^{\varepsilon}(0, \cdot)$ are smooth in  $\varepsilon$, with expansions
\begin{align}
& q_3^{\varepsilon}(0, \cdot)=\binom{0}{1}+ \Lambda_{\vep} q_1^{\vep}, 
\label{5.1}\\
& q_4^{\varepsilon}(0, \cdot)=\binom{\frac{1}{\alpha}}{0}+\varepsilon 
B_{\gamma}\binom{C\cos (x)}{-S\sin (x)}+\varepsilon^2 q_{42}^{\varepsilon}(x),\label{5.2}
\end{align}
where  $\Lambda_{\vep}=\vep\Lambda_0+\vep^2\Lambda_1(\vep)$, with
\beq\label{def-Lambda0}
\Lambda_0=\frac{\kappa\left(\alpha^2-\kappa^2\right)}{\alpha (2\alpha+\gamma)},  \qquad B_\gamma:=\frac{\alpha+\gamma}{\alpha(2\alpha+\gamma)}\left(\kappa-\frac{\gamma}{2}\tah(2\kappa)\right),
\eeq
and $\Lambda_1(\varepsilon)$ is smooth in $\varepsilon$.
Moreover, $q_{42}^{\varepsilon}(x)$ is smooth in $x$, with its first component being even and the second component being odd in $x$.

Finally, 
\begin{align}\label{def-d12ep0}
q_2^{\vep}(0,\cdot)=\hat{q}_2^{\varepsilon}(0,\cdot)+\mathfrak r_\vep q_4^\vep(0,\cdot), \qquad d_{12}(\vep,0)= 
\frac{\vep k_\vep\partial_\vep k_\vep}{c\partial_c k_\vep}
-\mathfrak r_\vep\Lambda_\vep=O(\vep^2),
\end{align}
where the constant $\mathfrak{r}_{\vep}$ is defined in \eqref{r-eps-def} and

$$
\hat{q}_2^{\varepsilon}(0, \cdot)=\left(\begin{array}{cc}
1 & 0 \\
-Z_{\varepsilon} & 1
\end{array}\right)\left(\partial_{\varepsilon}-\frac{\partial_{\varepsilon} k_{\varepsilon}}{\partial_c k_{\varepsilon}}\left(\begin{array}{cc}
\partial_c & 0 \\
0 & c^{-1} \partial_c(c\, \cdot)
\end{array}\right)\right)\binom{\eta_{\varepsilon}}{\varphi_{\varepsilon}}.
$$

\end{proposition}
\begin{rmk}\label{rem5.3}
Note that the background waves $\left(\eta_{\varepsilon}, \varphi_{\varepsilon}\right)$ depend on the wave speed $c$ through the parameters $\alpha=\frac{g h}{c^2}$ and $\gamma=\frac{\tilde{\gamma} h}{c}$ since the wave numbers $k_{\varepsilon}$ depend on these two parameters.
\end{rmk}
The kernel $q_1^{\varepsilon}(0, \cdot)$ can be found by using the translational invariance of the profile equation. Since $H_0[\eta_{\vep}](1)=G_0\left[\eta_{\varepsilon}\right](1)=0$, the profile equation is invariant under constant shifts of the velocity potential. This results in another kernel element $(0,1)^t$. Being the linear combination of $q_1^{\varepsilon}(0, \cdot)$ and $(0,1)^t$, we see that $q_3^{\varepsilon}(0, \cdot)$ also belongs to the kernel of $L_0^{\varepsilon}$. The generalized eigenfunction $q_2^{\varepsilon}(0, \cdot)$ can be found by differentiating the original profile equation before normalization in terms of the wave speed $c$ and the amplitude $\varepsilon$, while $q_4^{\varepsilon}(0, \cdot)$ is constructed using the Lyapunov-Schmidt method. The proof is given in Appendix \ref{app B}.

From the above proposition and expansions of the profiles in \eqref{ex}, we derive the following basis expansions.
\begin{corollary}[Expansions of the basis in $\varepsilon$]\label{cor-basisexpansion}
It holds that
$$
\begin{aligned}
&q_1^{\varepsilon}(0, \cdot)=\binom{-S \sin (x)}{C\cos (x)}+\varepsilon\binom{\partial_x \eta_2}{\partial_x \varphi_2-\kappa S\sin (x) \partial_x \eta_1}+\varepsilon^2\binom{\rm { odd }}{\rm { even }},\\
&\hat{q}_2^{\varepsilon}(0, \cdot)=\binom{S\cos (x)}{C\sin (x)}+\varepsilon\binom{2 \eta_2}{2 \varphi_2-\kappa S\sin (x) \eta_1}+\varepsilon^2\binom{\rm { even }}{\rm { odd }}, \\
& q_3^{\varepsilon}(0, \cdot)=\binom{0}{1}+\varepsilon \Lambda_0\binom{-S\sin (x)}{C\cos (x)}+\varepsilon^2\binom{\rm { odd }}{\rm { even }},
\end{aligned}
$$
and $q_4^{\varepsilon}(0, \cdot)$ expands as in \eqref{5.2}. We refer to \eqref{profile-expansion} for the definitions of $\eta_1, \eta_2, \varphi_2$.
Moreover, 
\beqs
q_2^\vep(0, \cdot)
=
\binom{S\cos x}{C\sin x}
+\vep
\binom{
2a_0^\sharp+2a_2\cos(2x)
}{
\Bigl(2W(0)-\frac{\kappa S^2}{2}\Bigr)\sin(2x)
}
+O(\vep^2), \quad \text{with} \;\;\; a_0^\sharp=-\frac{\kappa^2}{4\alpha}.
\eeqs

\begin{rmk}
    Let us remark that the coefficients in front of the second harmonics of $q_2^{\vep}$ do not enter the subsequent computations/analysis.
\end{rmk}

\end{corollary}
\begin{corollary}\label{cor:D0eps}
It holds that
$$
L_0^{\varepsilon}\left(q_1^{\varepsilon}, q_2^{\varepsilon}, q_3^{\varepsilon}, q_4^{\varepsilon}\right)(0, \cdot)=\left(q_1^{\varepsilon}, {q}_2^{\varepsilon}, q_3^{\varepsilon}, q_4^{\varepsilon}\right)(0, \cdot) \mathrm{D}_0^{\varepsilon},
$$
where the constant $4 \times 4$ matrix is
$$
\mathrm{D}_0^{\varepsilon}=\left(\begin{array}{cc}
\mathrm{A}_0^{\varepsilon} & 0 \\
0 & \mathrm{F}_0^{\varepsilon}
\end{array}\right) \quad \mathrm { with } \quad \mathrm{A}_0^{\varepsilon}=\left(\begin{array}{cc}
0 & d_{12}(\vep, 0) \\
0 & 0
\end{array}\right), \quad \mathrm{F}_0^{\varepsilon}=\left(\begin{array}{cc}
0 & -1 \\
0 & 0
\end{array}\right) .
$$
\end{corollary}

\subsection{Dual basis for $L_0^{\vep}$ }\label{sec-dualbasis-fixedconstant}
With a basis for the generalized kernel of $L_0^{\varepsilon}$ in hand, 
we shall construct the dual basis
\[
\left\{\tilde q_1^\vep,\tilde q_2^\vep,\tilde q_3^\vep,\tilde q_4^\vep\right\}
\subset \mathrm{gKer}\,(L_0^\vep)^*
\]
such that
\beq\label{dual-biorthogonality}
\left\langle \tilde q_j^\vep,q_k^\vep\right\rangle=\delta_{jk},
\qquad 1\leq j,k\leq4.
\eeq
Here the inner product is understood in
$\cY_{\per}$ and is conjugate-linear
in the first argument.
The above construction also gives the adjoint Jordan relations, namely,
\beq\label{adjoint-chain}
(L_0^\vep)^*\tilde q_2^\vep=0,\quad
(L_0^\vep)^*\tilde q_1^\vep=d_{12}(\vep,0)\tilde q_2^\vep, \quad (L_0^\vep)^*\tilde q_4^\vep=0,\quad
(L_0^\vep)^*\tilde q_3^\vep=-\tilde q_4^\vep.
\eeq
Recall that at $\vep=0$,
\beq\label{def-basis-00}
q_1^0=\binom{-S\sin x}{C\cos x},\qquad
q_2^0=\binom{S\cos x}{C\sin x},\qquad
q_3^0=\binom{0}{1},\qquad
q_4^0=\binom{1/\alpha}{0}.
\eeq
By \eqref{dual-00}, a direct normalization computation yields the following dual basis:
\beq\label{dualbasis-eps0}
\begin{aligned}
\tilde q_1^0
&=
\frac{1}{\pi P_\gamma}
\binom{-\dfrac{\alpha S}{\kappa}\sin x}{S\cos x}, \quad
\tilde q_2^0
=
\frac{1}{\pi P_\gamma}
\binom{\dfrac{\alpha S}{\kappa}\cos x}{S\sin x},
\quad
\tilde q_3^0
&=\dfrac{1}{2\pi}
\binom{0}{1},
\quad
\tilde q_4^0
=\dfrac{1}{2\pi}
\binom{\alpha}{0}.
\end{aligned}
\eeq
where $P_\gamma:=\frac{(2\alpha+\gamma)S^2}{\kappa}.$

\begin{prop}\label{prop-dualbasis}
For $|\vep|$ sufficiently small, there exists a unique smooth basis $\{\tilde q_j^\vep(0,\cdot)\}_{j=1}^4\subset \mathrm{gKer}\,(L_0^\vep)^*$ such that 
\eqref{dual-biorthogonality} holds and 
\begin{align}\label{paritysym}
    \tilde{q}_1^{\vep}(0,\cdot), \, \tilde{q}_3^{\vep}(0,\cdot) =\bpm \rm {odd}\\ \rm {even}\epm, \qquad    \tilde{q}_2^{\vep}(0,\cdot), \, \tilde{q}_4^{\vep}(0,\cdot) =\bpm  \rm {even}\\\rm {odd}\epm.
\end{align}
Moreover, the following expansions hold 
\begin{align*}
 \tilde q_j^{\vep}(0,\cdot)= \tilde q_j^{0}(0,\cdot) +\vep \tilde q_j^{1}(0,\cdot)+\cO(\vep^2),
\end{align*}
where $\tilde q_j^{0}(0,\cdot)$ is defined in
\eqref{dualbasis-eps0} and 
\beqs
\begin{aligned}
\tilde q_1^1(0,\cdot)
&=\bpm * \sin(2x)\\-\dfrac{\Lambda_0}{2\pi}+* \cos(2x)\epm, \qquad \tilde q_4^1(0,\cdot)
=\frac{\kappa^2}{2}\tilde q_2^0,
\\
\tilde q_2^1(0,\cdot)
&=\bpm * \cos(2x)\\ * \sin(2x)\epm+*\,\tilde q_4^0, \qquad \tilde q_3^1(0,\cdot)
=\frac{1}{2\pi}\bpm 
*\sin x\\
* \cos x \epm,
\end{aligned}
\eeqs
with
$*$ denoting some real constants uniquely determined by the second Fourier-mode parts of the adjoint equations. \footnote{Since the concrete expressions of these constants are irrelevant for the later use, we do not compute them.}
\end{prop}

\begin{proof}
Let $\Pi_0^\vep$ be the Riesz projector associated with the eigenvalue zero. Since $L_0^\vep$ depends smoothly on $\vep$ and zero has algebraic multiplicity four, $\Pi_0^\vep$ and $(\Pi_0^\vep)^*$ depend smoothly on $\vep$ and have rank four for $|\vep|$ sufficiently small. Set
\beqs
r_j^\vep:=(\Pi_0^\vep)^*\tilde q_j^0,\qquad \mathcal{G}_{jk}^\vep:=\ang{r_j^\vep,q_k^\vep}\,.
\eeqs
Since $\mathcal{G}^0=\Id$, the matrix $\mathcal{G}^\vep$ is invertible for small $\vep$. Therefore,
\beqs
\tilde q_j^\vep:=\sum_{\ell=1}^4\overline{\big((\mathcal{G}^\vep)^{-1}\big)_{j\ell}}\,r_\ell^\vep
\eeqs
defines a smooth basis of $\rm {gKer}(L_0^\vep)^*$ satisfying \eqref{dual-biorthogonality}. Let us point out that the constructed dual basis are real-valued. Indeed, 
since $(L_0^\vep)^*$ has real coefficients, and its Riesz
projector associated with zero preserves real-valuedness. 
Thus the  basis of $r_j^\vep$ and the matrix $\mathcal{G}^{\vep}$ are real,  and the same is true of $\{{\tilde q}_j^\vep\}_{j=1}^4$. For the 
uniqueness,  if $v\in {\rm {Ran}}(\Pi_0^\vep)^*={\rm {gKer}}(L_0^\vep)^*$ satisfies $\ang{v,u}=0$ for every $u\in {\rm {Ran}}\Pi_0^\vep={\rm {gKer}} (L_0^\vep)$, then, 
for any $w\in \cY_{\per}$
\beqs\ang{v,w}=\ang{(\Pi_0^\vep)^*v,w}=\ang{v,\Pi_0^\vep w}=0
\eeqs
thus $v=0.$

To prove the parity symmetries \eqref{paritysym}, we
define the reversibility operator $\mathcal{I}$ 
\beq\label{def-revisibility}
\mathcal{I} \bpm f \\ g \epm (x):={\bpm \overline {f(-x)}\\ -\overline {g(-x)}\epm }.
\eeq
From the parity of the wave profiles, we readily check that
\beqs
L_0^\vep\mathcal I=-\mathcal I L_0^\vep,\qquad (L_0^\vep)^*\mathcal I=-\mathcal I(L_0^\vep)^*.
\eeqs
Moreover, the basis constructed in Proposition~\ref{prop-basis} satisfies 
\beq\label{rev-basis}
\mathcal {I}q_j^\vep=(-1)^jq_j^\vep.
\eeq
To obtain the corresponding symmetry for the dual basis, set
$\bar q_j^\vep:=(-1)^j\mathcal I\tilde q_j^\vep. $
Then we have $\bar q_j^\vep\in\rm gKer(L_0^\vep)^*.$ 
Moreover, since $\mathcal I$ is anti-unitary,
\beqs
\begin{aligned}
\ang{\bar q_j^\vep,q_k^\vep}
&=(-1)^{j+k}\ang{\mathcal I\tilde q_j^\vep,\mathcal Iq_k^\vep}
=(-1)^{j+k}\overline{\ang{\tilde q_j^\vep,q_k^\vep}}
=\delta_{jk}.
\end{aligned}
\eeqs
By uniqueness of the dual basis, $\bar q_j^\vep=\tilde q_j^\vep$, and hence
\beq\label{rev-dualbasis}
\mathcal I\tilde q_j^\vep=(-1)^j\tilde q_j^\vep,
\eeq
which leads to the parity symmetries \eqref{paritysym}.

We now turn to the expansions of the dual basis. Write
\beqs
L_0^\vep=L_0^0+\vep L_0^1+O(\vep^2),\qquad \tilde q_j^\vep=\tilde q_j^0+\vep\tilde q_j^1+O(\vep^2).
\eeqs
By \eqref{adjoint-chain} and $d_{12}(\vep,0)=O(\vep^2)$,
\beqs
\begin{aligned}
(L_0^0)^*\tilde q_1^1&=-(L_0^1)^*\tilde q_1^0,&
(L_0^0)^*\tilde q_2^1&=-(L_0^1)^*\tilde q_2^0,\\
(L_0^0)^*\tilde q_4^1&=-(L_0^1)^*\tilde q_4^0,&
(L_0^0)^*\tilde q_3^1&=-\tilde q_4^1-(L_0^1)^*\tilde q_3^0.
\end{aligned}
\eeqs
Differentiating the duality relations \eqref{dual-biorthogonality} also gives
\beq
\ang{\tilde q_j^1,q_k^0}=-\ang{\tilde q_j^0,q_k^1}.
\eeq
The coefficients of $L_0^1$ contain only the first Fourier harmonics. Hence the first two adjoint equations, together with parity and the normalization conditions, give
\beqs
\tilde q_1^1=\binom{c_{11}\sin(2x)}{c_0+c_{12}\cos(2x)},\qquad
\tilde q_2^1=\binom{c_{21}\cos(2x)}{c_{22}\sin(2x)}+c_{20}\tilde q_4^0
\eeqs
for some real constants $c_{jk}$. Since $q_3^1=\Lambda_0q_1^0$, differentiating $\ang{\tilde q_1^\vep,q_3^\vep}=0$ yields
\beqs
2\pi c_0=\ang{\tilde q_1^1,q_3^0}=-\ang{\tilde q_1^0,q_3^1}=-\Lambda_0,
\eeqs
and therefore $c_0=-\Lambda_0/(2\pi)$. The constants multiplying the second Fourier harmonics are uniquely determined since the restriction of $(L_0^0)^*$ to the second Fourier mode is invertible.

We next determine $\tilde q_4^1$. Since $\tilde q_4^0=(\alpha/(2\pi),0)^{ t}$ and the first component of $L_0^1u$ has zero mean for every $u$, one has $(L_0^1)^*\tilde q_4^0=0$. Thus $\tilde q_4^1\in\ker(L_0^0)^*$. Its parity and the differentiated relations with $q_1^0,q_3^0,q_4^0$ show that 
$\tilde q_4^1=c\,\tilde q_2^0.$
Using $(q_2^1)_1=2a_0^\sharp+2a_2\cos(2x)$ and $a_0^\sharp=-\frac{\kappa^2}{4\alpha}$, we obtain
\beqs
c=\ang{\tilde q_4^1,q_2^0}=-\ang{\tilde q_4^0,q_2^1}=-2\alpha a_0^\sharp=\frac{\kappa^2}{2}.
\eeqs
Consequently, $\tilde q_4^1=\frac{\kappa^2}{2}\tilde q_2^0.$

Finally, the right-hand side of the equation for $\tilde q_3^1$ contains only the first Fourier harmonics. It then follows from parity that
$\tilde q_3^1=\binom{c_{31}\sin x}{c_{32}\cos x}$
where $c_{31},c_{32}\in\mathbb R$ are uniquely fixed by the adjoint equation and the differentiated duality conditions. This concludes the proof.
\end{proof}

\subsection{Basis and dual basis of $L_{\xi}^{\vep}, (\xi\neq 0)$ }
With the basis and dual basis of  $L_{0}^{\vep}$ available, we now construct a basis for the spectral subspace of $L_{\xi}^{\varepsilon}, (\xi\neq 0)$ corresponding to the four eigenvalues near zero. This can be done by means of the Kato extension operator with respect to the Floquet parameter $\xi$. 

By Lemma~\ref{lem:crossing-origin} and the compactness of
the resolvent of $L_0^0$, zero is an isolated eigenvalue
of $L_0^0$ of algebraic multiplicity four. We may
therefore choose $r_0>0$ sufficiently small such that
\beqs
\sigma(L_0^0)\cap\overline{B_{r_0}(0)}=\{0\}.
\eeqs
Since the operators $L_\xi^\vep$ have the common domain
$\cD_{\mathrm{per}}$ and depend continuously on $(\vep,\xi)$
in $\mathcal L(\cD_{\mathrm{per}},\cY_{\mathrm{per}})$,
we have
\beqs
\sup_{\lambda\in\partial B_{r_0}(0)}
\left\|
(L_\xi^\vep-L_0^0)(\lambda-L_0^0)^{-1}
\right\|_{\mathcal L(\cY_{\mathrm{per}})}
\rightarrow 0
\eeqs
as $|\vep|+|\xi|\to0$.
It follows from the resolvent identity
$$(\lambda I-L_\xi^\vep )^{-1}= (\lambda I-L_0^0 )^{-1}\big(\Id-(L_\xi^\vep-L_0^0)(\lambda I-L_0^0)^{-1}
\big)^{-1}$$
that $\partial B_{r_0}(0)\subset\rho(L_\xi^\vep)$
for $|\vep|\leq\vep_0$ and $|\xi|\leq\xi_0$, provided
$\vep_0,\xi_0>0$ are sufficiently small.
Moreover, the corresponding Riesz projector depends
continuously on $(\vep,\xi)$ and has constant rank four.

Let $|\xi| \leq \xi_0,|\varepsilon| \leq {\varepsilon}_0$ with $\xi_0>0,\vep_0>0$ sufficiently small.
Define the spectral projector
$$
\Pi_{\xi}^{\varepsilon}=\frac{1}{2\pi \,\i}\oint_{\partial B_{r_0}(0)}\left(\lambda I-L_{\xi}^{\varepsilon}\right)^{-1} \mathrm{~d} \lambda, \qquad (\Pi_{\xi}^{\varepsilon})^{*}=\frac{1}{2\pi \,\i}\oint_{\partial B_{r_0}(0)}\left(\lambda I-(L_{\xi}^{\varepsilon})^{*}\right)^{-1} \mathrm{~d} \lambda,
$$
We then define the Kato extension operator $\mathcal{U}^{\varepsilon}(\xi)$ and $\mathcal{V}^{\varepsilon}(\xi)$
by solving the ODEs
\begin{align}
&\partial_{\xi} \mathcal{U}^{\varepsilon}(\xi)=\left[\partial_{\xi} \Pi_{\xi}^{\varepsilon}, \Pi_{\xi}^{\varepsilon}\right] \mathcal{U}^{\varepsilon}(\xi), \qquad \qquad \,\,\,\mathcal{U}^{\varepsilon}(0)=\mathrm{Id}, \label{ODE-cU} \\
&\partial_{\xi} \mathcal{V}^{\varepsilon}(\xi)=\left[\partial_{\xi} (\Pi_{\xi}^{\varepsilon})^{*}, (\Pi_{\xi}^{\varepsilon})^{*}\right] \mathcal{V}^{\varepsilon}(\xi), \qquad \mathcal{V}^{\varepsilon}(0)=\mathrm{Id}. \label{ODE-cV}
\end{align}
Define 
\begin{equation}\label{basis-dualbasis}
q_j^{\varepsilon}(\xi, \cdot)=\mathcal{U}^{\varepsilon}(\xi) q_j^{\varepsilon}(0, \cdot), \qquad  \quad \tilde{q}_j^{\varepsilon}(\xi, \cdot)=\mathcal{V}^{\varepsilon}(\xi)\tilde{q}_j^{\varepsilon}(0, \cdot)  \qquad j=1,2,3,4,
\end{equation}
which, by Kato's perturbation theory \cite[Chapter 2, Section 4]{Kato}, 
form respectively a basis for the spectral subspace of $L_{\xi}^{\vep}$ and $(L_{\xi}^{\vep})^{*}$
associated to the spectrum inside of $B_{r_0}(0)$.  
We observe that
$$
\big[\p_\xi(\Pi_\xi^\vep)^*,
       (\Pi_\xi^\vep)^*\big]
=-\big[\p_\xi\Pi_\xi^\vep,
       \Pi_\xi^\vep\big]^*.
$$
It follows from the definitions \eqref{ODE-cU} and \eqref{ODE-cV} that
$
\frac{\d}{\d\xi}
\big((\cV^\vep(\xi))^*\cU^\vep(\xi)\big)=0.$
Since $\cU^\vep(0)=\cV^\vep(0)=\Id$, we obtain
$(\cV^\vep(\xi))^*\cU^\vep(\xi)=\Id.$
Therefore, the biorthogonality relations \eqref{dual-biorthogonality} are preserved for any $|\xi|\leq \xi_0$.

\section{Analysis of the representation matrix}
Having established the basis and its dual, we are now in a position to compute the representation matrix in this basis. Define the $4 \times 4$ matrix $\mathrm{D}_{\xi}^{\varepsilon}$ such that
$$
L_{\xi}^{\varepsilon}\left(q_1^{\varepsilon}, q_2^{\varepsilon}, q_3^{\varepsilon}, q_4^{\varepsilon}\right)(\xi, \cdot)=\left(q_1^{\varepsilon}, q_2^{\varepsilon}, q_3^{\varepsilon}, q_4^{\varepsilon}\right)(\xi, \cdot) \mathrm{D}_{\xi}^{\varepsilon} .
$$
Taking suitable inner products of the above identity with the dual basis $\left\{\tilde{q}_j^{\varepsilon}(\xi, \cdot)\right\}$, we find that
\begin{equation}\label{def_D}
\mathrm{D}_{\xi}^{\varepsilon}=\left(\left\langle\tilde{q}_j^{\varepsilon}(\xi, \cdot), L_{\xi}^{\varepsilon} q_{\ell}^{\varepsilon}(\xi, \cdot)\right\rangle\right)_{1 \leq j, \ell \leq 4}:=\left(d_{j \ell}(\varepsilon, \xi)\right)_{1 \leq j, \ell \leq 4}.
\end{equation}
We can, in fact, find the expansions of the entries $d_{j \ell}(\varepsilon, \xi)$ in $\xi$ and $\varepsilon$, using the definitions \eqref{basis-dualbasis}. However, it is not necessary to compute all of these expansions to determine whether the eigenvalues of $\mathrm{D}_{\xi}^{\varepsilon}$ are purely imaginary.

The first step is to get the the expression of $\mathrm{D}_{\xi}^{\varepsilon}$ when $\xi$ or $\varepsilon$ vanishes. We have already known $\mathrm{D}_0^{\varepsilon}$ from Corollary \ref{cor:D0eps}, so we now focus on determining $\mathrm{D}_{\xi}^0$.

\subsection{ Computation of $\mathrm{D}_{\xi}^0$}
In this subsection we assume $\varepsilon=0$ and $\xi \neq 0$. 
We infer from \eqref{eigenrelation} that
$$
L_{\xi}^0\left(v_1^{+}, v_{-1}^{+}
\right)=\left(v_1^{+}, v_{-1}^{+}
\right) \operatorname{diag}\left(\lambda_1^{+}, \lambda_{-1}^{+}\right)(\xi).
$$
To get $\mathrm{D}_{\xi}^0$, we need to relate $q_j^0(\xi, \cdot)(j=1,2)$ with $v_1^{+}(\xi, \cdot), v_{-1}^{+}(\xi, \cdot)$. 
\begin{lemma}\label{l5.7}
There is a $2 \times 2$  matrix 
$$
\begin{aligned}
T&=\operatorname{diag}\left(a_1(\xi), a_{-1}(\xi)\right) \frac{1}{2}\left(\begin{array}{cc}
1 & -\mathrm{i} \\
1 & \mathrm{i}
\end{array}\right)
\end{aligned}
$$
with
$a_1(\xi), a_{-1}(\xi)$ being some smooth coefficients,
such that
$$
\left(q_1^0, q_2^0\right)(\xi, \cdot)=\left(v_1^{+}, v_{-1}^{+}\right)(\xi, \cdot) T.
$$

\end{lemma}
\begin{proof}
Define
$$
{p}_{ \pm 1}(\xi, \cdot)=\mathcal{U}^0(\xi)\left(e^{ \pm \mathrm{i} x}\binom{ \pm \mathrm{i} \operatorname{sh}(\kappa)}{\operatorname{ch}(\kappa)}\right)=\pm\i\,\sh(\kappa)\,\mathcal{U}^0(\xi)\big(v_{\pm1}^{+}(0, \cdot)\big).
$$
Applying \eqref{extension-pm1}, we can find smooth functions  $a_{\pm}(\xi)$ such  that
${p}_{ \pm 1}(\xi, x)= a_{ \pm 1}(\xi)v_{\pm1}^{+}(\xi, x).$
The form of $T$ then follows from the fact that
$$
\begin{aligned}
\left(q_1^0, q_2^0\right)(\xi, \cdot)=&\mathcal{U}^0(\xi) \left(q_1^0(0, \cdot),q_2^0(0, \cdot)\right)\\
=&\mathcal{U}^0(\xi)\left(e^{ \mathrm{i} x}\binom{ \mathrm{i} \operatorname{sh}(\kappa)}{\operatorname{ch}(\kappa)}, e^{- \mathrm{i} x}\binom{ -\mathrm{i} \operatorname{sh}(\kappa)}{\operatorname{ch}(\kappa)} \right)\frac{1}{2}\left(\begin{array}{cc}
1 & -\mathrm{i} \\
1 & \mathrm{i}
\end{array}\right)\\
=&\left({p}_1, {p}_{-1}\right)(\xi, \cdot) \frac{1}{2}\left(\begin{array}{cc}
1 & -\mathrm{i} \\
1 & \mathrm{i}
\end{array}\right).
\end{aligned}
$$
\end{proof}
\begin{corollary}\label{cor_Dxi0} When $\varepsilon=0$, the representation matrix reads
$\mathrm{D}_{\xi}^0=\left(d_{j \ell}(0, \xi)\right)_{1 \leq j, \ell \leq 4}$,
where
$$d_{j \ell}(0, \xi)=0\quad \text{for } (j, \ell) \in(\{1,2\},\{3,4\}) \;\text{or } (\{3,4\},\{1,2\})$$
and
\begin{equation}\label{entry-ep=0}
\begin{aligned}
d_{11}(0, \xi)&=d_{22}(0, \xi)
=\mathrm{i} \kappa \xi\left(1-\nu_{*}\right)+\mathcal{O}\left(|\xi|^3\right), \\
d_{12}(0, \xi)&=-d_{21}(0, \xi)
=\delta_{*} (\kappa\xi)^2+\mathcal{O}\left(\xi^4\right), \\
d_{33}(0, \xi)&= \mathrm{i} (\kappa \xi -\gamma \,\tah(\kappa\xi)) \qquad 
d_{44}(0, \xi)=\mathrm{i} \kappa \xi, \\
d_{34}(0, \xi)&=-1, \quad d_{43}(0, \xi)=\alpha \kappa\xi \tah(\kappa\xi),
\end{aligned}
\end{equation}
with
\begin{equation}\label{signs}
\nu_{*}:=
\frac{\gamma}{2\ch^2(\kappa)}+ w^{\prime}(\kappa)<1, \quad 
\delta_{*}=\frac{1}{2} \left(\frac{\gamma  \tah(\kappa)}{\operatorname{ch}^2 (\kappa)}-w^{\prime\prime}(\kappa)\right)>0\,.
\end{equation}
\end{corollary}
\begin{proof}
We denote for short $\lambda_{\pm 1}=\lambda_{\pm 1}^{+}=\i\left(\zeta_j-\frac{\gamma\operatorname{th}(\zeta_j)} 2 - w(\zeta_j)\right),$ where $\zeta_j$ and $w(\zeta_j)$ are defined in \eqref{def-wzetaj}. It then follows from Lemma \ref{l5.7} 
that 
$$\mathrm{D}_{\xi}^0=\left(\begin{array}{cc} \rA_{\xi}^0 & 0\\ 0 & \rF_{\xi}^0 \end{array}\right), $$
 with
\begin{align*}
    \rA_{\xi}^0=T^{-1} \operatorname{diag}\left(\lambda_1, \lambda_{-1}\right) T =\frac 1 2 \left(\begin{array}{cc} \lambda_1+ \lambda_{-1}& \i\,\big(\lambda_{-1}-\lambda_{1}\big) \\
 \i\,\big(\lambda_{1}-\lambda_{-1}\big) & \lambda_1+ \lambda_{-1}\\   \end{array}\right).
\end{align*}
Therefore, 
\begin{align*}
d_{11}(0, \xi)&=d_{22}(0, \xi)=\frac{\mathrm{i}}{2}\bigg(2\kappa\xi - \frac{\gamma}{2}\big(\tah ( 
\zeta_1) + \tah (\zeta_{-1})\big)-w(\zeta_1)-w(\zeta_{-1})\bigg)\\
&\qquad\qquad\quad=\mathrm{i} \kappa \xi\left(1-\frac{\gamma }{2\operatorname{ch}^2 \kappa} -w^{\prime}(\kappa) \right)+\mathcal{O}\left(|\xi|^3\right), \\
d_{12}(0, \xi)&=-d_{21}(0, \xi)=\frac{1}{2}\left(2\kappa - \frac{\gamma}{2}\big( \tah \big( \zeta_1 \big) - \tah \big( \zeta_{-1} \big) \big)-w(\zeta_1)+w(\zeta_{-1})\right)\\
&\qquad\qquad\quad\;\;=\frac{\kappa^2\xi^2}{2} \left(\frac{\gamma \tah(\kappa)}{\operatorname{ch}^2 (\kappa)}-w^{\prime\prime}(\kappa)\right)+\mathcal{O}\left(\xi^4\right).
\end{align*}
Moreover, using the expression  \eqref{L0}
of $L_{\xi}^0,$ we compute readily that 
\begin{align*}
    L_{\xi}^0 q_3^0&= \binom{\mu\tah(\mu)}{\i\mu-\i\gamma\tah(\mu)}=\alpha\mu\tah(\mu) q_4^0+\i(\mu-\gamma\tah(\mu)) q_3^0, \\
     L_{\xi}^0 q_4^0&= \binom{\i\mu/\alpha}{-1}=\i\mu \,q_4^0-q_3^0,
\end{align*}
where $\mu=\kappa\xi.$
This yields that 
$$F_{\xi}^0=\bpm \i\big(\mu-\gamma\tah(\mu)\big) &-1 \\\alpha\mu\tah(\mu) & \i\mu \epm .$$
The proof for the sign of $\nu_{*},\delta_{*}$ in \eqref{signs} is
elementary but requires more computations. We present the details in Appendix \ref{Appendix-secondderi}.
\end{proof}

\subsection{Some useful properties of the entries.} We now gather some useful properties of the entries of $D_\xi^\vep,(\xi\neq0)$, defined in \eqref{def_D}.

\begin{lemma}
\label{lem:entries}
For $|\vep|$ and $|\xi|$ sufficiently small, the entries of $\rD_\xi^\vep$ have the following properties:\\
\textnormal{(1)}
$\quad
d_{j\ell}(\vep,\xi)\in \i\,\mathbb R,\quad \hbox{if } j+\ell \hbox{ is even},\quad \text{and }\quad
d_{j\ell}(\vep,\xi)\in \mathbb R,\quad \hbox{if } j+\ell \hbox{ is odd}.
$
\\
\textnormal{(2)}
There exist real-valued functions $m_{j\ell}(\vep,\xi)$, smooth in $(\vep,\xi)$, such that
\begin{align}
\tag{2a} d_{j\ell}&=\i\kappa\xi\,m_{j\ell},\qquad\quad  j,\ell\in \{1,3\} \,\,\mathrm{or}\,\, j,\ell\in \{2,4\},\\
\tag{2b} d_{j\ell}&=m_{j\ell},\qquad\qquad \; j\in \{1,3\},\, \ell\in \{2,4\},\\
\tag{2c} d_{j\ell}&=-\kappa^2\xi^2\,m_{j\ell},\quad \,\, \,\,j\in  \{2,4\},\ \ell\in \{1,3\},
\end{align}
and 
\begin{equation*}
\tag{2d}d_{14}(\vep,\xi)=O(\vep\kappa^2\xi^2),\qquad
d_{32}(\vep,\xi)=O(\vep\kappa^2\xi^2).
\end{equation*}
\textnormal{(3)}
$\quad d_{j\ell}(\vep,\xi)=O(\vep),\quad
(j,\ell)\in(\{1,2\},\{3,4\})\cup(\{3,4\},\{1,2\}).$ \\[3pt]
\textnormal{(4)}. 
For all $j,\ell,$ it holds that
\beqs
d_{j\ell}(\vep,-\xi)=\overline{d_{j\ell}(\vep,\xi)}
\eeqs
and $m_{j\ell}$ are even in $\xi.$

\end{lemma}

\begin{proof}
Let $\mathcal O:=\{1,3\}, \mathcal E:=\{2,4\}$ and $\cI$ be defined in \eqref{def-revisibility}.
In view of \eqref{rev-basis}, \eqref{rev-dualbasis}
we compute
\beqs
\begin{aligned}
\overline{d_{j\ell}(\vep,\xi)}
=\left\langle \mathcal I\tilde q_j^\vep(\xi,\cdot),\mathcal I L_\xi^\vep q_\ell^\vep(\xi,\cdot)\right\rangle 
&=-\left\langle (-1)^j\tilde q_j^\vep(\xi,\cdot),L_\xi^\vep (-1)^\ell q_\ell^\vep(\xi,\cdot)\right\rangle \\
&=(-1)^{j+\ell+1}d_{j\ell}(\vep,\xi).
\end{aligned}
\eeqs
This gives the property (1).

We next consider the property (2). If $j,\ell\in\mathcal O$ or $j,\ell\in\mathcal E$, then $d_{j\ell}$ is purely imaginary. Moreover, by Corollary \ref{cor:D0eps}, $d_{j\ell}(\vep,0)=0$ for these indices. Hence there is a real-valued smooth function $m_{j\ell}(\vep,\xi)$ such that (2a) holds.
If $j\in\mathcal O$ and $\ell\in\mathcal E$, then $d_{j\ell}$ is real-valued, and no additional power of $\xi$ is forced. Thus we get the property (2b).

For the property (2c), we consider $j\in\mathcal E$ and $\ell\in\mathcal O$. In this case $d_{j\ell}$ is real-valued and $d_{j\ell}(\vep,0)=0$. We prove also  that $\partial_\xi d_{j\ell}(\vep,0)=0$. Set
\beq\label{exp-operators}
L_1^\vep:=\frac{1}{\i \, \kappa}\partial_\xi L_\xi^\vep\big|_{\xi=0},\qquad \cU_1^\vep:=\frac{1}{\i \, \kappa} \p_{\xi}\cU^\vep(\xi)\big|_{\xi=0},\qquad \cV_1^\vep:=\frac{1}{\i \, \kappa}\partial_\xi \cV^\vep(\xi)\big|_{\xi=0},
\eeq
then these operators 
preserve real-valuedness. Differentiating
$d_{j\ell}(\vep,\xi)=\left\langle \tilde q_j^\vep(\xi,\cdot),L_\xi^\vep q_\ell^\vep(\xi,\cdot)\right\rangle$
at $\xi=0$, we obtain
\beqs
\frac{1}{\i \, \kappa}\partial_\xi d_{j\ell}(\vep,0)=-\left\langle \cV_1^\vep\tilde q_j^\vep(0,\cdot),L_0^\vep q_\ell^\vep(0,\cdot)\right\rangle+\left\langle \tilde q_j^\vep(0,\cdot),L_1^\vep q_\ell^\vep(0,\cdot)\right\rangle+\left\langle \tilde q_j^\vep(0,\cdot),L_0^\vep \cU_1^\vep q_\ell^\vep(0,\cdot)\right\rangle.
\eeqs
Since all the inner products are real, it holds that $
\partial_\xi d_{j\ell}(\vep,0)\in \i\mathbb R.$
Meanwhile, since $d_{j\ell}(\vep,\xi)$ is real-valued for $j\in\mathcal E$, $\ell\in\mathcal O$, we also have
$ \partial_\xi d_{j\ell}(\vep,0)\in \mathbb R.$
This gives the property (2c).

A similar argument shows that $\partial_\xi d_{14}(\vep,0)=\partial_\xi d_{32}(\vep,0)=0.$
On the other hand, Corollaries \ref{cor:D0eps} and \ref{cor_Dxi0} tell us that
\beqs
d_{14}(\vep,0)=d_{32}(\vep,0)=0,\qquad d_{14}(0,\xi)=d_{32}(0,\xi)=0.
\eeqs
Then the property (2d) holds.

Next, by Corollary \ref{cor_Dxi0}, the matrix $D_\xi^0$ is block diagonal with respect to the decomposition $\{1,2\}\cup\{3,4\}$. Hence all entries connecting these two blocks vanish when $\vep=0$. By smoothness in $\vep$,
the property (3) holds. 

 Finally, since $q_\ell^\vep(0,\cdot)$ and $\tilde q_j^\vep(0,\cdot)$ are real-valued, it follows from \eqref{ext-conj} in Proposition \ref{prop-extop} that $q_\ell^\vep(-\xi,\cdot)=\overline{q_\ell^\vep(\xi,\cdot)}$ and $\tilde q_j^\vep(-\xi,\cdot)=\overline{\tilde q_j^\vep(\xi,\cdot)}$. This, combined with the fact $L_{-\xi}=\overline{L_{\xi}},$  yields that $d_{j\ell}(\vep,-\xi)=\overline{d_{j\ell}(\vep,\xi)}$.
By (1), $d_{j\ell}(\vep,\cdot)$ is then even (resp.\ odd) when $j+\ell$ is odd (resp.\ even), and (2a)--(2c) show that every $m_{j\ell}(\vep,\cdot)$ is even.
This concludes the lemma.
\end{proof}

\subsection{Transformation of the representation matrix.}
In view of Lemma \ref{lem:entries}, we introduce the rescaled matrix
\beqs
\rM_{\xi}^{\vep}=(m_{j\ell}(\vep,\xi))_{1\leq j,\ell\leq4}:=\frac{1}{\i\kappa\xi}\mathcal P_\xi \rD_\xi^\vep\mathcal P_\xi^{-1},\qquad
\mathcal P_\xi:=\diag(\i\kappa\xi,1,\i\kappa\xi,1).
\eeqs
Then all entries $m_{j\ell}$ are real-valued and smooth in $(\vep,\xi)$ for $|\vep|$ and $|\xi|$ sufficiently small. Moreover, 
\beq \label{spec-relation}
\sigma(\rD_\xi^\vep)=\i\kappa\xi\,\sigma(\rM_{\xi}^{\vep}).
\eeq
With respect to the block decomposition $\{1,2\}\cup\{3,4\}$, 

we first collect the information on the diagonal blocks.
\begin{lemma}
\label{lem:diag-block-vorticity}
\, It holds that 

\textnormal{(1)} For the upper diagonal block,
\beqs
m_{11}=1-\nu_*+O(\vep^2+\xi^2),\qquad
m_{22}=1-\nu_*+O(\vep^2+\xi^2),
\eeqs
\beqs
m_{12}=\delta_*\kappa^2\xi^2+
\Gamma_{*}\vep^2
+O\big(\xi^2(\xi^2+|\vep|)+|\vep|^3\big),\quad m_{21}=\delta_*+O(\xi^2+|\vep|),
\eeqs
where $\Gamma_{*}=\frac{2\kappa k_2}{c\partial_c\kappa}-2\mathfrak b_2\Lambda_0$ and $\nu_*, \delta_{*}$ are defined in \eqref{signs}.\\
\textnormal{(2)}
For the lower diagonal block, 
\beqs
m_{33}=1-\gamma+O(|\vep|+\xi^2),\qquad
m_{44}=1+O(|\vep|+\xi^2),
\eeqs
\beqs
m_{34}=-1+O(|\vep|+\xi^2),\qquad
m_{43}=-\alpha+O(|\vep|+\xi^2).
\eeqs
\end{lemma}

\begin{proof} We first consider the upper diagonal block.
By Corollary \ref{cor_Dxi0}, \beqs
d_{12}(0,\xi)=-d_{21}(0,\xi)=\delta_*\kappa^2\xi^2+O(\xi^4).
\eeqs
On the other hand, Corollary \ref{cor:D0eps} and Proposition \ref{prop-basis} give
\beqs
d_{12}(\vep,0)=
\frac{\vep k_\vep\partial_\vep k_\vep}{c\partial_c k_\vep}
-\mathfrak r_\vep\Lambda_\vep=\left(\frac{2\kappa k_2}{c\partial_c\kappa}-2\frak{b}_2\Lambda_0\right)\vep^2+O(|\vep|^3)=\Gamma_{*}\vep^2+O(|\vep|^3),
\eeqs
where we used $k_\vep=\kappa+k_2\vep^2+O(\vep^4)$.
As in the proof of Lemma \ref{lem:entries} (2c), $\partial_\xi m_{12}(\varepsilon,0)$ is both real and purely imaginary, hence vanishes.
The above two identities then yield
\beqs
m_{12}=d_{12}
=\delta_*\kappa^2\xi^2+\Gamma_{*}\vep^2
+O\big(\xi^2(\xi^2+|\vep|)+|\vep|^3\big),\quad \text{and}\quad m_{21}=-\frac{d_{21}}{\kappa^2\xi^2}
=\delta_*+O(\xi^2+|\vep|).
\eeqs

Similarly, we get from Corollary \ref{cor_Dxi0} that
\beqs
d_{11}(0,\xi)=d_{22}(0,\xi)
=\i\kappa\xi\left(1-\frac{\gamma}{2\ch^{2}(\kappa)}-w'(\kappa)\right)+O(|\xi|^3),
\eeqs
and hence
\beqs
m_{11}(0,\xi)=m_{22}(0,\xi)=1-\nu_*+O(\xi^2).
\eeqs
It remains to prove the expansions of $m_{11}$ and $m_{22}$ when $\xi=0:$
\begin{equation}\label{m11m22}
m_{11}(\vep,0)=1-\nu_*+O(\vep^2),\qquad m_{22}(\vep,0)=1-\nu_*+O(\vep^2). 
\end{equation}
Since $d_{\ell\ell}=\i\kappa\xi m_{\ell\ell}$, one has
$m_{\ell\ell}(\vep,0)=\frac{1}{\i\kappa}\partial_{\xi}d_{\ell\ell}(\vep,0), \ell=1,2.$
Differentiating $d_{\ell\ell}(\vep,\xi)=\langle \tilde q_{\ell}^{\vep}(\xi,\cdot),L_{\xi}^{\vep}q_{\ell}^{\vep}(\xi,\cdot)\rangle$ at $\xi=0$ gives
\beqs
\begin{aligned}
m_{\ell\ell}(\vep,0)=&-\left\langle \cV_{1}^{\vep}\tilde q_{\ell}^{\vep}(0,\cdot),L_{0}^{\vep}q_{\ell}^{\vep}(0,\cdot)\right\rangle+\left\langle \tilde q_{\ell}^{\vep}(0,\cdot),L_{1}^{\vep}q_{\ell}^{\vep}(0,\cdot)\right\rangle+\left\langle \tilde q_{\ell}^{\vep}(0,\cdot),L_{0}^{\vep}\cU_{1}^{\vep}q_{\ell}^{\vep}(0,\cdot)\right\rangle,
\end{aligned}
\eeqs
where $L_1^{\vep}, \cV_1^{\vep}, \cU_1^{\vep}$ are defined in \eqref{exp-operators}.
By \eqref{adjoint-chain}, $(L_{0}^{\vep})^*\tilde q_{1}^{\vep}=\cO(\vep^2), (L_{0}^{\vep})^*\tilde q_{2}^{\vep}=0.$ Thus the two terms involving $\cU_{1}^{\vep}$ and $\cV_{1}^{\vep}$ are $O(\vep^2)$. Hence
\beqs
m_{\ell\ell}(\vep,0)=\left\langle \tilde q_{\ell}^{\vep}(0,\cdot),L_{1}^{\vep}q_{\ell}^{\vep}(0,\cdot)\right\rangle+O(\vep^2),\qquad \ell=1,2.
\eeqs
We expand the remaining scalar products to first order in $\vep$. For $\ell=1,2$, write
\beq\label{exp-L1}
L_1^\vep=L_1^0+\vep L_1^1+O(\vep^2),\qquad q_\ell^\vep=q_\ell^0+\vep q_\ell^1+O(\vep^2),\qquad \tilde q_\ell^\vep=\tilde q_\ell^0+\vep\tilde q_\ell^1+O(\vep^2).
\eeq
By  Corollary~\ref{cor-basisexpansion} and Proposition~\ref{prop-dualbasis},
the vectors $q_\ell^0$ and $\tilde q_\ell^0$ contain only the first Fourier mode, whereas $q_\ell^1$ and $\tilde q_\ell^1$ contain only the modes $0$ and $2$. Moreover, $L_1^0$ has constant coefficients, while the coefficients of $L_1^1$ contain only the first Fourier mode. It follows that
\beqs
\begin{aligned}
m_{\ell\ell}(\vep,0)=m_{\ell\ell}(0,0)+\vep\big(&\ang{\tilde q_\ell^1,L_1^0q_\ell^0}+\ang{\tilde q_\ell^0,L_1^1q_\ell^0}+\ang{\tilde q_\ell^0,L_1^0q_\ell^1}\big)+O(\vep^2).
\end{aligned}
\eeqs
Each scalar product in the coefficient of $\vep$ vanishes by Fourier orthogonality. 
Therefore, $m_{\ell\ell}(\vep,0)=m_{\ell\ell}(0,0)+O(\vep^2), \ell=1,2$.
By Corollary~\ref{cor_Dxi0}, $m_{11}(0,0)=m_{22}(0,0)=1-\nu_*$, and hence the equalities \eqref{m11m22} are valid. This leads to the stated expansions of $m_{11}$ and $m_{22}$.

For the lower diagonal block, Corollary \ref{cor_Dxi0} yields, with $\mu=\kappa\xi$,
\beqs
d_{33}(0,\xi)=\i\,(\mu-\gamma\,\tah\mu),\qquad
d_{44}(0,\xi)=\i\,\mu,\qquad
d_{34}(0,\xi)=-1,\qquad
d_{43}(0,\xi)=\alpha \mu\,\tah\mu.
\eeqs
Therefore
\beqs
m_{33}(0,\xi)=1-\gamma\frac{\tah(\kappa\xi)}{\kappa\xi}
=1-\gamma+O(\xi^2),\qquad
m_{44}(0,\xi)=1,
\eeqs
and
\beqs
m_{34}(0,\xi)=-1,\qquad
m_{43}(0,\xi)=-\alpha\frac{\tah(\kappa\xi)}{\kappa\xi}
=-\alpha+O(\xi^2).
\eeqs
The stated expansions then follow from smoothness in $\vep$. This completes the proof.
\end{proof}

We next collect the information on the off-diagonal blocks.
\begin{lemma}
\label{lem:offdiag-block-vorticity}
All entries of the off-diagonal blocks 
are of order $\cO(\vep)$. More precisely, there exist real constants $\mathfrak{m}_{j\ell},$ such that
$$m_{j\ell}(\vep,\xi)=\vep\big(\mathfrak m_{j\ell}+O(|\vep|+\xi^2)\big),\qquad
(j,\ell)\in(\{1,2\},\{3,4\})\cup(\{3,4\},\{1,2\}).$$ 
In addition,
\beqs
m_{14}(\vep,\xi)=\cO(\vep\kappa^2\xi^2),\qquad
m_{32}(\vep,\xi)=\cO(\vep\kappa^2\xi^2).
\eeqs
Moreover, 
\begin{align}
\mathfrak{m}_{13}&= -\kappa+\frac{\kappa(\alpha^2-\kappa^2)(\kappa^2-\alpha-\alpha^2)}
{\alpha(2\alpha+\gamma)^2}, \label{def-m13}\\
\mathfrak{m}_{42}&=-\frac{\alpha\kappa P_{\gamma}}{2}
-\frac{\kappa^2}{2}\nu_{*}.
\label{def-m42}
\end{align}
\end{lemma}

\begin{proof}
The fact that all entries in the off-diagonal blocks are of order $\cO(\vep)$ follows from property (3) of Lemma \ref{lem:entries} and the definition of $m_{j\ell}$.
Besides, the improved estimates for $m_{14}$ and $m_{32}$ follow from the property (2d) of Lemma \ref{lem:entries} and the definition of the rescaling $m_{14}=d_{14}, m_{32}=d_{32}.$

We now turn to the analysis of $m_{13}, m_{42}.$ Proposition \ref{prop-basis} and \eqref{adjoint-chain} imply that $$L_0^\vep q_3^\vep(0,\cdot)=0, \quad L_0^\vep q_2^\vep(0,\cdot)=O(\vep^2), \quad (L_0^\vep)^{*} \tilde{q}_1^\vep(0,\cdot)=\cO(\vep^2), \quad (L_0^\vep)^{*} \tilde{q}_4^\vep(0,\cdot)=0,$$ 
and thus the Kato extension operators contribute only $O(\vep^2)$. Hence
\beqs
m_{13}(\vep,0)=\ang{\tilde q_1^\vep,L_1^\vep q_3^\vep}+O(\vep^2),\qquad m_{42}(\vep,0)=\ang{\tilde q_4^\vep,L_1^\vep q_2^\vep}+O(\vep^2).
\eeqs
Note that we have suppressed the arguments $(0,\cdot)$ in $q_2^\vep, q_3^\vep, \tilde{q}_1^\vep, \tilde{q}_4^\vep$ for the sake of clarity.
Let us compute the  coefficients $\m_{13}$ and $\m_{42}$ which satisfy 
\beqs
\m_{13}
=\ang{\tilde q_1^1,L_1^0q_3^0}+\ang{\tilde q_1^0,L_1^1q_3^0}+\ang{\tilde q_1^0,L_1^0q_3^1},
\quad
\m_{42}
=\ang{\tilde q_4^1,L_1^0q_2^0}+\ang{\tilde q_4^0,L_1^1q_2^0}+\ang{\tilde q_4^0,L_1^0q_2^1}.
\eeqs
We begin with $\m_{13}$. The expressions of $\tilde q_1^0$ and $\tilde q_1^1$ in Proposition \ref{prop-dualbasis}, combined with the results for $L_1^0q_3^0$, $L_1^1q_3^0$ and $L_1^0 q_1^0$ in Corollary \ref{cor-algebraic}, and the identity $q_3^1=\Lambda_0 q_1^0,$ show that
\beqs
\begin{aligned}
\m_{13}
&=\ang{\tilde q_1^1,L_1^0q_3^0}+\ang{\tilde q_1^0,L_1^1q_3^0}+\ang{\tilde q_1^0,L_1^0q_3^1}\\
&=-(1-\gamma)\Lambda_0-\kappa+\frac{\Lambda_0}{P_\gamma}(SC-\kappa)\\
&=-\kappa+\frac{\Lambda_0}{P_\gamma}\bigl(SC-\kappa-(1-\gamma)P_\gamma\bigr).
\end{aligned}
\eeqs
Using $\alpha+\gamma=\kappa C/S$, $P_\gamma=(2\alpha+\gamma)S^2/\kappa$ and $\Lambda_0=\kappa(\alpha^2-\kappa^2)/(\alpha(2\alpha+\gamma))$, we obtain \eqref{def-m13}.

We next compute $\m_{42}$. By a similar argument to the above, combining the expressions of $\tilde q_4^0,$ $\tilde q_4^1$, and $\tilde q_2^0$ from Proposition \ref{prop-dualbasis} with the results for $L_1^0q_2^0$, $\int_0^{2\pi}(L_1^1q_2^0)_1\,\d x$ and $\int_0^{2\pi}(L_1^0q_2^1)_1\,\d x$ from Corollary \ref{cor-algebraic}, we find that
\beqs
\begin{aligned}
\m_{42}
&=\ang{\tilde q_4^1,L_1^0q_2^0}+\ang{\tilde q_4^0,L_1^1q_2^0}+\ang{\tilde q_4^0,L_1^0q_2^1}\\
&=\frac{\kappa^2}{2P_\gamma}(SC-\kappa)-\frac{\alpha\kappa P_\gamma}{2}+2\alpha a_0^\sharp\\
&=-\frac{\alpha\kappa P_\gamma}{2}-\frac{\kappa^2}{2}+\frac{\kappa^2}{2P_\gamma}(SC-\kappa).
\end{aligned}
\eeqs
Since $\nu_*:=1-\frac{SC-\kappa}{P_\gamma},$ this gives \eqref{def-m42}.
\end{proof}

\section{Eigenvalues of the transformed matrix $\rM_{\xi}^{\vep}$ and  proof of Theorem \ref{thm-main}}
Combining Lemmas \ref{lem:diag-block-vorticity} and \ref{lem:offdiag-block-vorticity}, the rescaled matrix has the form
\beqs
\rM_\xi^\vep
=\begin{pmatrix}
\rM_{11}&\vep\rM_{12}\\
\vep\rM_{21}&\rM_{22}
\end{pmatrix}(\vep,\xi),
\eeqs
where each block is $2\times2$ and real-valued and
\begin{equation} \label{def-M11}
\rM_{11}(\vep,\xi)=
\begin{pmatrix}
1-\nu_*+O(\vep^2+\xi^2)&
\delta_*\kappa^2\xi^2+\Gamma_{*}\vep^2
+O\big(\xi^2(\xi^2+|\vep|)+|\vep|^3\big)\\
\delta_*+O(|\vep|+\xi^2)&
1-\nu_*+O(\vep^2+\xi^2)
\end{pmatrix},    
\end{equation}
and
\begin{equation} \label{def-M22}
\rM_{22}(\vep,\xi)=
\begin{pmatrix}
1-\gamma&-1\\
-\alpha&1
\end{pmatrix}
+O(|\vep|+\xi^2).
\end{equation}

We now compute the leading eigenvalues of the diagonal blocks $\rM_{11}$ and $\rM_{22}$ and show that they are separated.
\begin{lemma}\label{lem:separation-diagonal-blocks}
Let $X_1^\pm(\vep,\xi)$ and $X_2^\pm(\vep,\xi)$ be the eigenvalues of $\rM_{11}(\vep,\xi)$ and $\rM_{22}(\vep,\xi)$, respectively. Then
\beqs
X_1^\pm(\vep,\xi)=1-\nu_*+O(|\vep|+|\xi|),
\eeqs
whereas
\beqs
X_2^\pm(\vep,\xi)=1-\frac{\gamma}{2}\pm\frac12\sqrt{\gamma^2+4\alpha}+O(|\vep|+\xi^2).
\eeqs
Moreover, there exist $\vep_0>0$, $\xi_0>0$ and $\iota_0>0$ such that, for $|\vep|\leq\vep_0$ and $|\xi|\leq\xi_0$,
\beqs
\left|X_1^\mu(\vep,\xi)-X_2^\sigma(\vep,\xi)\right|\geq \iota_0,\qquad \mu,\sigma\in\{+,-\}.
\eeqs
\end{lemma}

\begin{proof}
The estimates for the eigenvalues $X_1^\pm(\vep,\xi)$ follow directly from $\rM_{11}(\vep,\xi)$ in \eqref{def-M11}.
For the second diagonal block $\rM_{22}(\vep,\xi)$ in \eqref{def-M22}, 
the two eigenvalues of the leading matrix are
\beqs
X_2^\pm(\vep,\xi)=X_{2,0}^\pm+O(|\vep|+\xi^2),\qquad \text{with}\quad X_{2,0}^\pm:=1-\frac{\gamma}{2}\pm\frac12\sqrt{\gamma^2+4\alpha}.
\eeqs
Since $\alpha>0, \sqrt{\gamma^2+4\alpha}>|\gamma|,$ it holds that $X_{2,0}^+>1.$
Moreover, $X_{2,0}^+X_{2,0}^-=1-(\gamma+\alpha)=1-\frac{\kappa}{\tah(\kappa)}<0,$ thus $ X_{2,0}^-<0$.
These, together with the fact $0<1-\nu_*<1$, yield
\beq\label{X20pm}
X_{2,0}^-<0<1-\nu_*<1<X_{2,0}^+.
\eeq
Thus the limiting eigenvalues of $\rM_{11}$ lie strictly between the two limiting eigenvalues of $\rM_{22}$. By continuity of the eigenvalues with respect to $(\vep,\xi)$, the separation persists for $|\vep|$ and $|\xi|$ sufficiently small. This proves the lemma.
\end{proof}

We now continue the analysis of the real matrix $\rM_\xi^\vep$. 
Let $\rho_*:=\sqrt{\gamma^2+4\alpha}$, then $X_{2,0}^{\pm}=1-\frac{\gamma}{2}\pm\frac{\rho_*}{2}$.
The two eigenvalues of the leading part of $\rM_{22}$ are $X_{2,0}^{+}$ and $X_{2,0}^{-}$, and a diagonalizing matrix is
\beqs
\rG_0:=
\begin{pmatrix}
1&1\\
-\dfrac{\gamma+\rho_*}{2}&\dfrac{\rho_*-\gamma}{2}
\end{pmatrix},\qquad \text{with}\quad
\rG_0^{-1}
=\frac1{\rho_*}
\begin{pmatrix}
\dfrac{\rho_*-\gamma}{2}&-1\\
\dfrac{\gamma+\rho_*}{2}&1
\end{pmatrix}.
\eeqs
Therefore, for $|\vep|+|\xi|$ sufficiently small, there exists a real invertible matrix
$\rG(\vep,\xi)=\rG_0+O(|\vep|+\xi^2)$,
such that
\beqs
\rM_{22}=\rG
\begin{pmatrix}
X_2^+&0\\
0&X_2^-
\end{pmatrix}
\rG^{-1},
\qquad
X_2^\pm=X_{2,0}^\pm+O(|\vep|+\xi^2).
\eeqs

We conjugate $\rM_\xi^\vep$ by this matrix in the lower block and define
\beqs
\rN_\xi^\vep:=
\begin{pmatrix}
\Id&0\\
0&\rG^{-1}
\end{pmatrix}
\rM_\xi^\vep
\begin{pmatrix}
\Id&0\\
0&\rG
\end{pmatrix}
=
\begin{pmatrix}
\rN_{11}&\vep \rN_{12}\\
\vep \rN_{21}&\rN_{22}
\end{pmatrix},
\eeqs
where
\beq \label{def-N}
\rN_{11}=\rM_{11},\qquad
\rN_{12}=\rM_{12}\rG,\qquad
\rN_{21}=\rG^{-1}\rM_{21},\qquad
\rN_{22}=\begin{pmatrix}X_2^+&0\\0&X_2^-\end{pmatrix}.
\eeq
The next step is to eliminate the upper off-diagonal block.

\begin{lemma}
\label{lem:eliminate-upper}
There exists a real $2\times2$ matrix $\rQ=\rQ(\vep,\xi)$, smooth in $(\vep,\xi)$, such that
\beqs
\widetilde \rM_\xi^\vep:=
\begin{pmatrix}
\mathrm{Id} &-\vep \rQ\\
0& \mathrm{Id}
\end{pmatrix}
\rN_\xi^\vep
\begin{pmatrix}
\mathrm{Id}&\vep \rQ\\
0&\mathrm{Id}
\end{pmatrix}
=
\begin{pmatrix}
\rN_{11}-\vep^2 \rQ\rN_{21}&0\\
\vep \rN_{21}&\rN_{22}+\vep^2\rN_{21}\rQ
\end{pmatrix}.
\eeqs
\end{lemma}

\begin{proof}
A direct computation shows that it suffices to solve
\beqs
\rN_{11}\rQ-\rQ\rN_{22}+\rN_{12}-\vep^2 \rQ\rN_{21}\rQ=0.
\eeqs
Writing $\rQ=(\rV_+,\rV_-)$ and $\rN_{12}=(\rU_+,\rU_-)$ in columns, this equation is equivalent to
\beqs
\begin{cases}
(\rN_{11}-X_2^+\Id)\rV_+ + \rU_+ -\vep^2 \rQ\rN_{21}\rV_+=0,\\
(\rN_{11}-X_2^-\Id)\rV_- + \rU_- -\vep^2 \rQ\rN_{21}\rV_-=0.
\end{cases}
\eeqs
By Lemma \ref{lem:separation-diagonal-blocks}, the matrices $\rB_\pm:=\rN_{11}-X_2^\pm\Id$
are invertible for $|\vep|+|\xi|$ sufficiently small. Hence the implicit function theorem gives a smooth solution $\rQ$, with
\beq\label{def-V+V-}
\rV_+=-\rB_+^{-1}\rU_+ + O(\vep^2),\qquad
\rV_-=-\rB_-^{-1}\rU_- + O(\vep^2).
\eeq
This proves the lemma.
\end{proof}

The eigenvalues of $\rM_{\xi}^{\vep}$ are therefore the union of the eigenvalues of the two real $2\times2$ matrices
\beqs
\rN_{11}-\vep^2 \rQ\rN_{21},\qquad
\rN_{22}+\vep^2 \rN_{21}\rQ.
\eeqs
Since $\rN_{22}$ is diagonal and its two diagonal entries are separated, the second matrix has two real eigenvalues
\beqs
\widetilde X_2^\pm=X_2^\pm+O(\vep^2)=X_{2,0}^\pm+O(|\vep|+\xi^2).
\eeqs

It remains to compute the discriminant of the upper-left block. Denoting $-\rQ\rN_{21}=(n_{j\ell})_{1\leq j,\ell\leq2}$, we shall only need the leading term of $n_{12}$. By the definitions of $\rN_{12}, \rN_{21}$ in \eqref{def-N} and the information about the off-diagonal entries $\mathfrak m_{13}, \mathfrak m_{42}, m_{14}, m_{32}$ in Lemma \ref{lem:offdiag-block-vorticity}, we obtain, at leading order,
\beqs
\rN_{12}=
\begin{pmatrix}
\mathfrak m_{13}&\mathfrak m_{13}\\
*&*
\end{pmatrix}
+O(|\vep|+\xi^2),\qquad
\rN_{21}=
\frac1{\rho_*}
\begin{pmatrix}
*&-\mathfrak m_{42}\\
*&\mathfrak m_{42}
\end{pmatrix}
+O(|\vep|+\xi^2).
\eeqs
Let $\Delta_\pm:=1-\nu_*-X_{2,0}^\pm,$ then by \eqref{X20pm},
\beq\label{sign-D*}
\Delta_+\Delta_-=(1-\nu_*-X_{2,0}^+)(1-\nu_*-X_{2,0}^-)=\nu_*^2-\gamma\nu_*-\alpha<0.
\eeq
Using the leading approximation
\beqs
\rB_\pm^{-1}=
\begin{pmatrix}
\dfrac1{\Delta_\pm}&0\\
-\dfrac{\delta_*}{\Delta_\pm^2}&\dfrac1{\Delta_\pm}
\end{pmatrix}
+O(|\vep|+\xi^2),
\eeqs
together with the definitions $\rQ=(\rV_+,\rV_-), \rN_{12}=(\rU_+,\rU_-)$  and \eqref{def-V+V-}, we find
\beqs
n_{12}=\mathfrak n_{12}+O(|\vep|+\xi^2),\quad\text{with}\quad
\mathfrak n_{12}
=\frac{\mathfrak m_{13}\mathfrak m_{42}}{\rho_*}
\left(\frac1{\Delta_-}-\frac1{\Delta_+}\right)
=-\frac{\mathfrak m_{13}\mathfrak m_{42}}{\nu_*^2-\gamma\nu_*-\alpha}.
\eeqs
Then the two eigenvalues of 
\beqs
\rN_{11}-\vep^2 \rQ\rN_{21}=
\begin{pmatrix}
m_{11}+\vep^2n_{11}&m_{12}+\vep^2n_{12}\\
m_{21}+\vep^2n_{21}&m_{22}+\vep^2n_{22}
\end{pmatrix}
\eeqs
are
\beqs
\widetilde X_1^\pm
=\frac12\left(m_{11}+m_{22}+\vep^2(n_{11}+n_{22})
\pm \sqrt{\Delta}\right),
\eeqs
where
\beqs
\Delta
=\left(m_{11}-m_{22}+\vep^2(n_{11}-n_{22})\right)^2
+4\left(m_{21}+\vep^2n_{21}\right)
\left(m_{12}+\vep^2n_{12}\right).
\eeqs
Using the expansions of the entries of $\rN_{11}=\rM_{11}$ in Lemma \ref{lem:diag-block-vorticity} and the formula for $n_{12}$ above, we obtain
\beq\label{def-Delta-1}
\Delta
=4\delta_*
\left(\delta_*\kappa^2\xi^2+\mathcal C_*\vep^2\right)
+O\big((\xi^2+\vep^2)(\xi^2+|\vep|)\big),
\eeq
where $\delta_*=\frac12\left(\frac{\gamma\,\sh(\kappa)}{\ch^3(\kappa)}-w''(\kappa)\right)>0$ is defined in \eqref{signs} and
\beq\label{def-cC}
\mathcal C_{*}=\Gamma_{*}+\mathfrak n_{12}
=\frac{2\kappa k_2}{c\partial_c\kappa}-2\mathfrak b_2\Lambda_0
-\frac{\mathfrak m_{13}\mathfrak m_{42}}{\nu_*^2-\gamma\nu_*-\alpha}\,.
\eeq
Based on this, we are able to prove our Theorem \ref{thm-main}.

\begin{proof}[Proof of Theorem~\ref{thm-main}]
By \eqref{spec-relation}, the four eigenvalues of $L_{\xi}^{\vep}$ 
in the ball $B_{r_0}(0)$  are $$\lambda_{j}^{\pm}=\i \kappa \xi \widetilde {X}_j^\pm,\qquad (j=1,2),$$
where from the preceding block reduction,  the four eigenvalues $\widetilde {X}_j^\pm$ of $\rM_{\xi}^{\vep}$ take the form
\beqs
\widetilde {X}_1^\pm
=1-\nu_*+O(|\vep|+\xi^2)
\pm \frac12 \sqrt{\Delta}, \qquad \widetilde  {X}_2^\pm=X_{2,0}^\pm+O(|\vep|+|\xi|).
\eeqs
The eigenvalues $\widetilde X_2^\pm$ remain real for $|\vep|$ and $|\xi|$ sufficiently small. If $\mathcal C_*>0$, the expansion \eqref{def-Delta-1} implies
$\Delta\geq0$ for $|\vep|+|\xi|$ sufficiently small,
and hence all four eigenvalues $\lambda_{1}^{\pm},\lambda_{2}^{\pm}$ are purely imaginary.

Assume now $\mathcal C_*<0$. We first prove that
$\Delta(\vep,\xi)$ is even in $\xi$. By Lemma \ref{lem:entries} (4), the characteristic polynomial $\det(X-\rM_\xi^\vep)$ is even in $\xi$. Since $\widetilde X_1^\pm$ are separated from $\widetilde X_2^\pm$ by \eqref{X20pm}, the  factor $(X-\widetilde X_1^+)(X-\widetilde X_1^-)$ is even in $\xi$ and hence so is its discriminant $\Delta.$
By smoothness and \eqref{def-Delta-1},
$\partial_\xi^2\Delta(0,0)
=8\delta_*^2\kappa^2>0.$
Taking $\vep_0,\xi_0$ smaller if necessary, we conclude that
$\partial_\xi^2\Delta(\vep,\xi)>0$.
Together with
$\partial_\xi\Delta(\vep,0)=0$,
this shows that $\Delta(\vep,\xi)$ is strictly
increasing for $0<\xi\leq\xi_0$.
On the other hand, \eqref{def-Delta-1} implies
\beqs
\Delta(\vep,0)
=4\delta_*\mathcal C_*\vep^2+O(|\vep|^3)<0
\eeqs
when $\vep\neq0$ is sufficiently small, whereas
$\Delta(\vep,\xi_0)>0$.
Consequently, there exists exactly one
$\xi_*(\vep)\in(0,\xi_0)$ satisfying
$\Delta(\vep,\xi_*(\vep))=0$.
Moreover, for any fixed $0<a<\rho_0<b$, where
$\rho_0:=\frac{\sqrt{-\mathcal C_*}}{\kappa\sqrt{\delta_*}}, $
the expansion \eqref{def-Delta-1} shows that
$\Delta(\vep,a|\vep|)<0$ and
$\Delta(\vep,b|\vep|)>0$ for small $|\vep|$.
Hence $\xi_*(\vep)/|\vep|$ remains in a compact
subinterval of $(0,\infty)$. Substituting
$\xi=\xi_*(\vep)$ into \eqref{def-Delta-1}, we obtain
$\delta_*\kappa^2\frac{\xi_*(\vep)^2}{\vep^2}
+\mathcal C_*=O(|\vep|).$
Therefore,
$\xi_*(\vep)=\rho_0|\vep|+O(\vep^2)$.
The asserted spectral alternatives follow from the sign
of $\Delta$ and the reality of
$\widetilde X_2^\pm$.
\end{proof}

\section{Analysis of the modulational instability criterion}
\subsection{Comparison with the NLS coefficient in \cite{TKM2012}}
\label{app:TKM-comparison}

In this subsection we explain how the coefficient obtained from our spectral reduction is related to the nonlinear coefficient in the formal nonlinear Schr\"odinger (NLS) approximation derived by Thomas--Kharif--Manna \cite{TKM2012}. 

We now translate our notation into that of \cite{TKM2012}. Denote $\sigma:=\tah(\kappa ).$ Since in our convention the horizontal velocity is given by $U_1=\Phi_x-\gamma y$, while the sign convention for the shear parameter in \cite{TKM2012} is opposite, their dimensionless vorticity parameter $X,$ written in our normalization, is $X=-\frac{\gamma\sigma}{\kappa}.$
As a result, the parameters $\alpha$ and $\gamma$ in this paper are related to $X$ by 
\beq\label{ag-X}
\gamma=-\frac{X\kappa}{\sigma}, \qquad\alpha=\frac{\kappa}{\sigma}(1+X).
\eeq
Moreover, in our normalization the  phase speed and the fluid depth are normalized to be one. Thus, the notation of \cite{TKM2012} corresponds to our notation as follows 
\beqs
\mu=k=\kappa,\qquad c_p=1,\qquad \rho=\frac{c_g}{c_p}=\nu_*, \qquad \omega=\kappa,
 \qquad \overline{\Omega}=\frac{X}{\sigma}. 
\eeqs
Thomas--Kharif--Manna derive an NLS equation of the form
\beq\label{NLSeq}
\i \, \p_{\tau} a+L_{\rm TKM}a_{\xi\xi}=M_{\rm TKM}|a|^2a,
\eeq
where
\beqs
\begin{aligned}
L_{\rm TKM}
&=
\frac{1}{\kappa \sigma(2+X)} \bigg(\kappa (1-\sigma^2)
\bigl[1-\kappa\sigma+(1-\nu_*)X\bigr]-\sigma\nu_*^2\bigg),\\
M_{\rm TKM}&=\frac{\kappa^3}{8(1+X)(2+X)\sigma^4}\Big(U_{\rm TKM}+V_{\rm TKM}W_{\rm TKM}\Big),
\end{aligned}
\eeqs
with
\beq\label{def-UVW}
\begin{aligned}
U_{\rm TKM}&=9-12\sigma^2+13\sigma^4-2\sigma^6
+(27-18\sigma^2+15\sigma^4)X\\
&+(33-3\sigma^2+4\sigma^4)X^2+(21+5\sigma^2)X^3+(7+2\sigma^2)X^4+X^5, \\
V_{\rm TKM}
&=
(1+X)^2\left(1+\nu_*+\frac{\kappa X}{\sigma}\right)
+1+X-\nu_*\sigma^2-\kappa\sigma X,\\
W_{\rm TKM}&=2\sigma^3
\frac{(1+X)(2+X)+\nu_*(1-\sigma^2)}
{\sigma\nu_*^2+\kappa\nu_*X-\kappa(1+X)}, \quad \text{and} \quad \nu_{*}=\frac{1+X+\kappa(1-\sigma^2)/\sigma}{2+X}.
\end{aligned}
\eeq

In \cite[Section III]{TKM2012}, the authors linearize  \eqref{NLSeq} around a Stokes wave solution $\underline{a}_{\vep}=\vep S \exp(-\i M_{\rm TKM}S^2\vep^2 \tau)$
and find that the linearized operator has spectrum of form 
\beq\label{NLS-spectrum}
\lambda^2=-(\kappa\xi)^2L_{\rm TKM}\Big(2 S^2 M_{\rm TKM}\,\vep^2+(\kappa\xi)^2L_{\rm TKM}\Big)
\eeq
where we have chosen the characteristic frequency $\mu=\kappa\xi.$ As a result, the wave $\underline{a}_{\vep}$ is modulational unstable when 
$\vep$ and $\xi$ are sufficiently small
if and only if $L_{\rm TKM}M_{\rm TKM}<0. $

On the other hand, in our spectral analysis,  the eigenvalues of the modulational block satisfy
\beq\label{modu-spectrum}
\Big(\lambda-\i \kappa \xi c_{*}\Big)^2=-\delta_*(\kappa\xi)^2\Big(\delta_*(\kappa\xi)^2+\mathcal C_{*}\varepsilon^2\Big)+\cO\Big(|\xi|^2(|\xi|^2+\varepsilon^2)^{\f32}\Big), \quad 
\eeq
with $c_{*}=1-\nu_{*}+\cO(\vep^2+\xi^2).$
Now, we would like to find the relation between the coefficients 
$L_{\rm TKM}, M_{\rm TKM}$ and 
$\delta_{*}$ and $\cC_{*}.$
First, using the definition of $\delta_{*}$ in \eqref{signs} and the relation \eqref{ag-X}, we compute that 
\beq\label{delta-corespond}
{\delta_*}=\frac{1}
{2\kappa-\gamma\sigma} \bigg(\nu_*^2-(1-\sigma^2)\Big(\alpha(1-\kappa\sigma)+\gamma(\nu_*-\kappa\sigma)\Big)\bigg)\,=- L_{\rm TKM}\,.
\eeq

We next consider $\cC_{*}$ defined in \eqref{def-cC}. Using the definition of $k_2$ in \eqref{def-k2}, $\mathfrak{b}_2$ in \eqref{def-fb2}, $\Lambda_0$ in \eqref{def-Lambda0}, as well as the relation \eqref{ag-X}, 
we obtain, after some algebraic computations that
\beq\label{Gamma-corespond}
\Gamma_*
=-\frac{\kappa^3S^2}{4(1+X)(2+X)\sigma^4}\,U_{\rm TKM}.
\eeq

In fact, substituting \eqref{ag-X} and the relations
$\sh^2(\kappa)=\frac{\sigma^2}{1-\sigma^2}, \sh(2\kappa)=\frac{2\sigma}{1-\sigma^2}$
into the definition of $k_2$ in \eqref{def-k2}, we obtain
\beqs
k_2=\frac{\kappa^3P}{8\sigma (1-\sigma^2)\big(\sigma-\kappa (1-\sigma^2)\big)},
\eeqs
where
\beqs
\begin{aligned}
P={}&9-12\sigma^2+13\sigma^4-2\sigma^6+(18-12\sigma^2+10\sigma^4)X\\
&+(15-9\sigma^2+4\sigma^4)X^2+6(1-\sigma^2)X^3+(1-2\sigma^2)X^4.
\end{aligned}
\eeqs
We next compute $c\p_c\kappa$. Recall that $\alpha=gh/c^2$ and $\gamma=\tilde\gamma h/c$, and hence
$c\p_c\alpha=-2\alpha, c\p_c\gamma=-\gamma.$
Differentiating the dispersion relation $\alpha+\gamma=\kappa/\sigma$ with respect to $c$ and 
using $2\alpha+\gamma=\kappa(2+X)/\sigma$, we deduce that $c\p_c\kappa=-\frac{\kappa(2+X)\sigma}{\sigma-\kappa (1-\sigma^2)},$
and therefore
\beqs
\frac{2\kappa k_2}{c\p_c\kappa}=-\frac{\kappa^3P}{4(2+X)\sigma^2(1-\sigma^2)}.
\eeqs
On the other hand, using again \eqref{ag-X}, the formulas for $\mathfrak{b}_2$ in \eqref{def-fb2} and $\Lambda_0$ in \eqref{def-Lambda0} yield
\beqs
a_0=-\frac{\kappa\sigma Q}{4(1-\sigma^2)},\qquad \mathfrak b_2=-\alpha a_0-\frac{\kappa^2}{4}=\frac{\kappa^2XR}{4(1-\sigma^2)},\qquad \Lambda_0=\frac{\kappa Q}{(1+X)(2+X)}.
\eeqs
where
\beqs
Q:=(1+X)^2-\sigma^2,\qquad R:=X^2+3X+3-\sigma^2.
\eeqs
Therefore,
\beqs
-2\mathfrak b_2\Lambda_0=-\frac{\kappa^3XQR}{2(1+X)(2+X)(1-\sigma^2)},
\eeqs
and hence
\beqs
\Gamma_*=\frac{2\kappa k_2}{c\partial_c\kappa}-2\mathfrak b_2\Lambda_0=-\frac{\kappa^3}{4(1+X)(2+X)\sigma^2(1-\sigma^2)}\bigg((1+X)P+2X\sigma^2QR\bigg).
\eeqs
 Direct computations of powers of $X$ show that $(1+X)P+2X\sigma^2QR=U_{\rm TKM}$, which yields \eqref{Gamma-corespond}.

We now examine $\frak{n}_{12}.$ By \eqref{def-m13}, \eqref{def-m42},  and the relation \eqref{ag-X}, it holds that 
\beqs
\mathfrak m_{13}
=-\frac{\kappa}{(1+X)(2+X)}V_{\rm TKM}, \quad
\mathfrak m_{42}=-\frac{\kappa^2}{2(1-\sigma^2)} \big((1+X)(2+X)+\nu_*(1-\sigma^2)\big).
\eeqs
Moreover, from \eqref{ag-X} and \eqref{sign-D*}, we get that 
\begin{align*}
    \nu_*^2-\gamma\nu_*-\alpha = \frac{1} {\sigma}  \bigg( \sigma\nu_*^2+\kappa\nu_*X-\kappa(1+X)\bigg)<0.
\end{align*}
It then follows that
\beqs
\frak{n}_{12}=-\frac{\mathfrak m_{13}\mathfrak m_{42}}{\nu_*^2-\gamma\nu_*-\alpha}=-\frac{\kappa^3S^2}{4(1+X)(2+X)\sigma^4} V_{\rm TKM}W_{\rm TKM}.
\eeqs
Combining this with \eqref{Gamma-corespond}, we find that 
\beqs
\mathcal C_{*}=\Gamma_{*}+\mathfrak n_{12}
=-\frac{\kappa^3S^2}{4(1+X)(2+X)\sigma^4}
\Big(U_{\rm TKM}+V_{\rm TKM}W_{\rm TKM}\Big)=-2S^2M_{\rm TKM}.
\eeqs
This, together with \eqref{delta-corespond}, shows that the instability criterion for the Stokes waves in \eqref{modu-spectrum} coincides exactly with that obtained by linearizing the NLS equation about its standing-wave solution in \eqref{NLS-spectrum}.

\subsection{Proof of Proposition~\ref{prop:stability}}\label{sec:sign-index}
From the preceding subsection, we can write
\begin{align}\label{rewrite}
    \cC_{*}
    = \frac{\kappa^3 S^2}{4(1+X)(2+X)\sigma^4}  \frac{\widetilde{\cC}_{*}}{\gamma\nu_*+\alpha-\nu_*^2},
\end{align}
where 
\begin{align}\label{def-cctilde}
    \widetilde{\cC}_{*}(X,\kappa)=-(\gamma\nu_*+\alpha-\nu_*^2)U_{\rm TKM}+V_{\rm TKM}\cdot 2\sigma^2
\big((1+X)(2+X)+\nu_*(1-\sigma^2)\big),
\end{align}
with $U_{\rm TKM}, V_{\rm TKM}$ given in \eqref{def-UVW}.

\begin{proof}[Proof of Proposition~\ref{prop:stability}]
From \eqref{sign-D*},  \eqref{ag-X} and \eqref{rewrite}, we obtain that the sign of $\cC_{*}$ is determined by that of $\tilde{\cC}_{*}$ defined in \eqref{def-cctilde}.
We first prove $(2)$ and the second claim of $(1)$ by asymptotic expansion. In fact, using $\tah(\kappa)=1+O(e^{-2\kappa})$,
we compute, for fixed $X,$
\beqs
\widetilde{\mathcal C}_{*}(X,\kappa)
=-2\kappa(1+X)(3X+2)(X^2+2X+2)
+\cO(1),
\qquad \kappa\to+\infty.
\eeqs
The second claim of $(1)$ follows from the fact $1+X>0$ and $X^2+2X+2>0.$ For fixed $\kappa>0,$ 
\beqs
\widetilde{\mathcal C}_*(X,\kappa)
=\left(\sigma^2+(1-\sigma^2)\left(\frac{\kappa}{\sigma}-1\right)^2\right)X^5
+O_\kappa(X^4),\qquad X\to+\infty.
\eeqs
The leading coefficient is positive, proving $(2)$.

It remains to prove the first claim of $(1)$. Fix $\kappa>0$ and
$-1<X\leq-2/3.$ In view of the appearance of 
$1+X, \sigma^2, \frac{\kappa}{\sigma}$ in \eqref{def-UVW}, we set
\beqs
r=1+X\in(0,1/3],\qquad
y=\sinh^2\kappa>0,\qquad e=\frac{\sh (2\kappa)}{2\kappa}>1.
\eeqs
It then follows that 
$\sigma^2=\frac{y}{1+y},
\frac{\kappa}{\sigma}=\frac{1+y}{e},
\nu_*=\frac{r+1/e}{1+r}$ ($\nu_{*}$ is obtained from \eqref{def-UVW}).
Substituting these identities into the expression for
$\widetilde{\mathcal C}_*$, we obtain
\begin{equation}\label{eq:all-depth-normalization}
\widetilde{\mathcal C}_*(X,\kappa)
=\frac{ry}{(1+r)^2e^2(1+y)^2}\Phi_r(y,e),
\end{equation}
where
\beqs
\Phi_r(y,e)=A_r(y)+(e-1)B_r(y)+(e-1)^2C_r(y).
\eeqs
Here, with $p=r^2+r+1$, the three coefficient functions are
\beqs
\begin{aligned}
A_r(y)={}&(1+r)^2p^2
+(1-r^2)^2(2r^2+2r+1)y-(1+r)^3(1-r+5r^2-r^3)y^2,\\
B_r(y)={}&2(1+r)p(3r^2+3r+1)
+2r^2(1-r^2)(3r+1)y+2(1+r)^2(1-3r)(1+r^2)y(y-1),\\
C_r(y)={}&\frac{r^2p^2}{y}
+r(2r^5+10r^4+15r^3+16r^2+9r+4)+r^2(1+r)(r^3+7r^2+5r+7)y.
\end{aligned}
\eeqs
Note that these expressions also make sense at $r=0$ for $y>0.$
For $0\leq r\leq1/3$ and $y>0$, we have $C_r(y)\geq0$. Now we consider the term $B_r(y)$. Since
\beqs
0\leq(1+r)^2(1-3r)(1+r^2)
=1-r-4r^2-4r^3-5r^4-3r^5\leq 1
\eeqs
and $y(y-1)\geq-1/4$, we obtain
\beqs
B_r(y)\geq2-\frac12(1+r)^2(1-3r)(1+r^2)\geq\frac32.
\eeqs
Consequently,
\beqs
\partial_e\Phi_r(y,e)=B_r(y)+2(e-1)C_r(y)\geq\frac32
\qquad(e\geq1),
\eeqs
so $\Phi_r(y,e)$ is strictly increasing in $e$.

Suppose first that $0<y\leq1$. Note from the above estimates that $\Phi_r(y,e)\geq A_r(y)$ for $y>0$. Moreover, the function $A_r$ is  concave and quadratic
in $y$, and its endpoint values satisfy
\beqs
A_r(0)=(1+r)^2p^2>0,\qquad
A_r(1)=(1+r)^2\bigl(1+2r(1-r)(1-2r^2)\bigr)>0.
\eeqs
Thus $\Phi_r(y,e)>0$ for $0<y\leq1$.

Suppose next that $y\geq1$, and set $z=\sqrt y=\sinh\kappa\geq1$.
From the Taylor expansion of $\cosh\kappa$, we obtain
\beqs
\begin{aligned}
24\left(\cosh\kappa-\frac32\kappa\right)
&\geq\kappa^4+12\kappa^2-36\kappa+24>0.
\end{aligned}
\eeqs
Hence $e>3z/2$. 
By monotonicity of $\Phi_r$ in $e$,
$\Phi_r(y,e)>\Phi_r(z^2,3z/2)\geq G_r(z)$
with
\beqs
G_r(z):=\Phi_r(z^2,3z/2)
-\frac{r^2p^2}{z^2}(3z/2-1)^2.
\eeqs
It suffices to prove $G_r(z)>0$ for $0\leq r\leq1/3$ and $z\geq1$.
By substituting the expressions of $A_r$, $B_r$
and $C_r$ into the definition of $G_r$, and
collecting powers of $z-1$,  we obtain 
\beqs
G_r(z)-(1-3r)G_0(z)-27r^3G_{1/3}(z)
=r(1-3r)\sum_{j=0}^{5}p_j(r)(z-1)^j,
\eeqs
where
\begin{equation*}
\begin{aligned}
p_0(r)&=\frac{1728+7236r-436r^2-363r^3-171r^4}{108},\quad
p_1(r)=\frac{1296+5913r-269r^2-348r^3-234r^4}{54},\\
p_2(r)&=\frac{216+351r+341r^2-102r^3-216r^4}{36},\qquad
p_3(r)=\frac{486+81r+698r^2+204r^3-90r^4}{27},\\
p_4(r)&=\frac{2268+2349r+2027r^2+789r^3-117r^4}{108},\quad
p_5(r)=3(r^3+2r^2+2r+2).
\end{aligned}
\end{equation*}
Using $0\leq r\leq1/3$, we find
\beqs
\begin{aligned}
&436r^2+363r^3+171r^4\leq192r,\qquad 
269r^2+348r^3+234r^4\leq137r,\\
&102r^3+216r^4\leq58r^2,\qquad
90r^4\leq10r^2,\qquad 117r^4\leq39r^3.
\end{aligned}
\eeqs
Therefore
\beqs
p_0\geq16,\quad p_1\geq24,\quad p_2\geq6,\quad
p_3\geq18,\quad p_4\geq21,\quad p_5\geq6,
\eeqs
and hence
\beqs
G_r(z)\geq(1-3r)G_0(z)+27r^3G_{1/3}(z).
\eeqs
Both endpoint polynomials are positive for $z\geq1$. Indeed,
\beqs
\begin{aligned}
G_0(z)
&=3z^5-3z^4-3z^3+3z^2+3z-1\\
&=2+3(z-1)+6(z-1)^2+15(z-1)^3
+12(z-1)^4+3(z-1)^5>0,
\end{aligned}
\eeqs
while
\beqs
\begin{aligned}
729G_{1/3}(z)
&=256z^4-2640z^3+6999z^2+2928z-1548\\
&=\frac{z^2}{4}\bigl((32z-165)^2+771\bigr)
+2928z-1548>0.
\end{aligned}
\eeqs
Since  $1-3r$ and $27r^3$ do not vanish simultaneously on $[0,1/3]$,
we conclude that $G_r(z)>0$.
Therefore $\Phi_r(y,e)>0$ also for $y\geq1$. By
\eqref{eq:all-depth-normalization},
$\widetilde{\mathcal C}_*(X,\kappa)>0$ for $\kappa>0$ and
$-1<X\leq-2/3$. The proof is thus finished.
\end{proof}


\appendix
\section{Expansions of the small-amplitude periodic waves.}\label{appen-expansion}
We compute in this section the expansions of the steady periodic waves used in the main text, whose existence has been proven in several works; see for instance \cite{SimmenSaffman1985,SilvaPeregrine2018,Constaintin-Strauss-CPAM,CSV-acta,Erik-Weber,Kozlov-Lokharu}.

After the rescaling $X\mapsto x=kX$, the fluid domain is flattened onto $\cS:=\T\times[-1,0]$ by
\beqs
Y:=1+y,\qquad \rho(x,y):=y+Y\eta(x),\qquad \p_x^\rho:=\p_x-\f{Y\eta_x}{1+\eta}\p_y,\qquad \p_y^\rho:=\f1{1+\eta}\p_y.
\eeqs
The steady equations become
\beq\label{elliptic}
k^2(\p_x^\rho)^2\Phi+(\p_y^\rho)^2\Phi=0\quad\hbox{in }\cS,\qquad \Phi_y|_{y=-1}=0,
\eeq
and, at $y=0$,
\begin{align}
&k^2\eta_x\Phi_x-\f{1+k^2\eta_x^2}{1+\eta}\Phi_y=k(1+\gamma\eta)\eta_x,\\
&\f12k^2(\p_x^\rho\Phi)^2+\f12\f{\Phi_y^2}{(1+\eta)^2}-k(1+\gamma\eta)\p_x^\rho\Phi+(\alpha+\gamma)\eta+\f{\gamma^2}{2}\eta^2=0. \label{Bernouli}
\end{align}

\begin{prop}\label{prop-expan-profile}
Let $\alpha>0$ and $\gamma\in\mR$ be fixed, and assume that $\kappa>0$ satisfies
\beq\label{relation}
\alpha+\gamma=\frac{\kappa}{\tah(\kappa)}.
\eeq
The family of periodic solutions to
\eqref{elliptic}--\eqref{Bernouli},
parameterized by $\vep$ with wave number $k_\vep$,
admits the following expansions as $\vep\to0$:
\beq\label{ex}
\eta_\vep=\vep\eta_1+\vep^2\eta_2+\vep^3\eta_3+O(\vep^4),\quad \Phi_\vep=\vep\Phi_1+\vep^2\Phi_2+\vep^3\Phi_3+O(\vep^4),\quad k_\vep=\kappa+\vep^2k_2+O(\vep^4),
\eeq
where
\begin{align}
\eta_1&=S\cos x,\quad \eta_2=a_0+a_2\cos(2x),\quad \Phi_1=\ch(\kappa Y)\sin x, \quad
\Phi_2=W(y)\sin(2x), \label{profile-expansion}\\
k_2&=\f{\cN(\kappa,\gamma)}{2\kappa\sh(2\kappa)(\sh(2\kappa)-2\kappa)},\label{def-k2}
\end{align}
with
\beqs
\begin{aligned}
&a_0=-\f{\sh^2(\kappa)}{4(\alpha+\gamma)}(\alpha^2-\kappa^2),\qquad
a_2=\f{\sh^2(\kappa)}{4\kappa\tah\kappa}\left(\alpha^2+\frac{2\alpha\kappa}{\tah(2\kappa)}-3\kappa^2+\frac{2\kappa^2}{\tah\kappa\,\tah(2\kappa)}\right),\\
 &W(y)=\f{\kappa S}{2}Y\sh(\kappa Y)+\f{a_2-\frac{\sh^2(\kappa)}{4}(2\alpha+\gamma)}{\sh(2\kappa)}\ch(2\kappa Y),\\
&\mathcal N(\kappa,\gamma)=\gamma^4\sh^4(\kappa)\bigl(1-\sh^2(\kappa)\bigr)-6\gamma^3\kappa\,\ch(\kappa)\sh^3(\kappa)+\gamma^2\kappa^2\sh^2(\kappa)\bigl(10\sh^4(\kappa)+21\sh^2(\kappa)+15\bigr)\\
&-2\gamma\kappa^3CS\bigl(8\sh^4(\kappa)+12\sh^2(\kappa)+9\bigr)+\kappa^4\bigl(8\sh^6(\kappa)+16\sh^4(\kappa)+15\sh^2(\kappa)+9\bigr).
\end{aligned}
\eeqs
\end{prop}

\begin{proof}
Let us first write down the following identities that will be used throughout
\beqs
\begin{aligned}
(\p_x^\rho)^2\Phi={}&\Phi_{xx}-2\f{Y\eta_x}{1+\eta}\Phi_{xy}+\left(\f{Y\eta_x}{1+\eta}\right)^2\Phi_{yy}+\left(-\f{Y\eta_{xx}}{1+\eta}+\f{2Y\eta_x^2}{(1+\eta)^2}\right)\Phi_y,\\
(\p_y^\rho)^2\Phi={}&\f1{(1+\eta)^2}\Phi_{yy}.
\end{aligned}
\eeqs
We substitute the expansions \eqref{ex} into the flattened equations. 
At order $\vep$, one obtains
\beqs
\begin{aligned}
&\kappa^2\Phi_{1,xx}+\Phi_{1,yy}=0,\qquad \Phi_{1,y}|_{y=-1}=0,\\
&-\Phi_{1,y}=\kappa\eta_{1,x},\quad -\kappa\Phi_{1,x}+(\alpha+\gamma)\eta_1=0,\quad\hbox{at }y=0.
\end{aligned}
\eeqs
The stated expressions of $\eta_1$ and $\Phi_1$ satisfy the first three relations, while the last one is equivalent to the dispersion relation \eqref{relation}.

At order $\vep^2$, the interior equation is
\beqs
\kappa^2\Phi_{2,xx}+\Phi_{2,yy}=2\kappa^2Y\eta_{1,x}\Phi_{1,xy}+\kappa^2Y\eta_{1,xx}\Phi_{1,y}+2\eta_1\Phi_{1,yy},
\eeqs
while the kinematic and Bernoulli conditions become
\begin{align*}
    &\kappa^2\eta_{1,x}\Phi_{1,x}-\Phi_{2,y}+\eta_1\Phi_{1,y}=\kappa\eta_{2,x}+\kappa\gamma\eta_1\eta_{1,x},\\    &\f{\kappa^2}{2}\Phi_{1,x}^2+\f12\Phi_{1,y}^2-\kappa\Phi_{2,x}+\kappa\eta_{1,x}\Phi_{1,y}-\kappa\gamma\eta_1\Phi_{1,x}+(\alpha+\gamma)\eta_2+\f{\gamma^2}{2}\eta_1^2=0.
\end{align*}
Using the first-order relations and setting $\eta_2=a_0+a_2\cos(2x)$ and $\Phi_2=W(y)\sin(2x)$, these equations are reduced to
\beqs
\begin{aligned}
&W''-4\kappa^2W=\kappa^2S\left(\ch(\kappa Y)-\f32\kappa Y\sh(\kappa Y)\right),\qquad W'(-1)=0,\\
&W'(0)=2\kappa a_2-\f{\kappa}{2}(\alpha-1)S^2,\qquad 2\kappa W(0)=(\alpha+\gamma)a_2+\f{S^2}{4}(\alpha^2+\kappa^2).
\end{aligned}
\eeqs
The mean part of the Bernoulli condition gives $(\alpha+\gamma)a_0+\f{S^2}{4}(\alpha^2-\kappa^2)=0.$
A solution of the interior equation satisfying the bottom condition is
\beqs
W(y)=\f{\kappa S}{2}Y\sh(\kappa Y)+D\ch(2\kappa Y),
\eeqs
and therefore
\beqs
W(0)=\f{\kappa S^2}{2}+D\ch(2\kappa),\qquad W'(0)=\f{\kappa S}{2}(S+\kappa C)+2\kappa D\sh(2\kappa).
\eeqs
Substitution into the two surface relations gives the stated formulas for $a_2$ and $D$, and the preceding mean equation gives $a_0$.

It remains to determine $k_2$. Only the first harmonic of the order-$\vep^3$ equations is needed. Writing
\beqs
\eta_3=b_1\cos x+a_3\cos(3x),\qquad \Phi_3=Z_1(y)\sin x+Z_3(y)\sin(3x),
\eeqs
the cubic interior equation is
\beqs
\begin{aligned}
\kappa^2\Phi_{3,xx}+\Phi_{3,yy}={}&-2\kappa k_2\Phi_{1,xx}+2\eta_1\Phi_{2,yy}-(3\eta_1^2-2\eta_2)\Phi_{1,yy}+2\kappa^2Y\eta_{1,x}\Phi_{2,xy}\\
&+2\kappa^2Y(\eta_{2,x}-\eta_1\eta_{1,x})\Phi_{1,xy}-\kappa^2Y^2\eta_{1,x}^2\Phi_{1,yy}+\kappa^2Y\eta_{1,xx}\Phi_{2,y}\\
&+\kappa^2Y(\eta_{2,xx}-\eta_1\eta_{1,xx}-2\eta_{1,x}^2)\Phi_{1,y}.
\end{aligned}
\eeqs
The first-harmonic projections of the nonlinear terms are
\beqs
\begin{aligned}
2\eta_1\Phi_{2,yy}&\mapsto SW'',\\
-(3\eta_1^2-2\eta_2)\Phi_{1,yy}&\mapsto\kappa^2\ch(\kappa Y)\left(2a_0-a_2-\f34S^2\right),\\
2\kappa^2Y\eta_{1,x}\Phi_{2,xy}+\kappa^2Y\eta_{1,xx}\Phi_{2,y}&\mapsto\f32\kappa^2YSW',\\
2\kappa^2Y(\eta_{2,x}-\eta_1\eta_{1,x})\Phi_{1,xy}&+\kappa^2Y(\eta_{2,xx}-\eta_1\eta_{1,xx}-2\eta_{1,x}^2)\Phi_{1,y}\mapsto-\f34\kappa^3YS^2\sh(\kappa Y),\\
-\kappa^2Y^2\eta_{1,x}^2\Phi_{1,yy}&\mapsto-\f34\kappa^4Y^2S^2\ch(\kappa Y).
\end{aligned}
\eeqs
Hence
\beqs
\begin{aligned}
Z_1''-\kappa^2Z_1={}&2\kappa k_2\ch(\kappa Y)+SW''+\f32\kappa^2YSW'+\kappa^2\ch(\kappa Y)\left(2a_0-a_2-\f34S^2\right)\\
&-\f34\kappa^3YS^2\sh(\kappa Y)-\f34\kappa^4Y^2S^2\ch(\kappa Y).
\end{aligned}
\eeqs
Substituting the expression of $W$, the terms containing $Y^2\ch(\kappa Y)$ cancel and we obtain
\beqs
Z_1''-\kappa^2Z_1=2\kappa k_2\ch(\kappa Y)+F_1(Y),
\eeqs
where
\beqs
\begin{aligned}
F_1(Y)={}&\f{\kappa^2}{4}\left((8a_0-4a_2+S^2)\ch(\kappa Y)+2\kappa YS^2\sh(\kappa Y)\right)\\
&+4\kappa^2DS\ch(2\kappa Y)+3\kappa^3YDS\sh(2\kappa Y).
\end{aligned}
\eeqs
A convenient solution satisfying the bottom condition is
\beqs
Z_1(y)=E\ch(\kappa Y)+k_2Y\sh(\kappa Y)+\f{\kappa}{2}(2a_0-a_2)Y\sh(\kappa Y)+\f{\kappa^2S^2}{8}Y^2\ch(\kappa Y)+\kappa DSY\sh(2\kappa Y).
\eeqs
Thus
\beqs
Z_1(0)=EC+k_2S+P_1,\qquad Z_1'(0)=E\kappa S+k_2(S+\kappa C)+P_1',
\eeqs
where
\beqs
\begin{aligned}
P_1&=\f{\kappa S}{2}(2a_0-a_2)+\f{\kappa^2S^2}{8}C+\kappa DS\sh(2\kappa),\\
P_1'&=\f{\kappa}{2}(2a_0-a_2)(S+\kappa C)+\f{\kappa^2S^2}{8}(2C+\kappa S)+\kappa DS\big(\sh(2\kappa)+2\kappa\ch(2\kappa)\big).
\end{aligned}
\eeqs
At order $\vep^3$, the kinematic condition reads
\beqs
\begin{aligned}
&\kappa^2(\eta_{1,x}\Phi_{2,x}+\eta_{2,x}\Phi_{1,x})-\Phi_{3,y}+\eta_1\Phi_{2,y}-(\eta_1^2-\eta_2+\kappa^2\eta_{1,x}^2)\Phi_{1,y}\\
&\hspace{35mm}=\kappa\eta_{3,x}+k_2\eta_{1,x}+\kappa\gamma(\eta_1\eta_{2,x}+\eta_2\eta_{1,x}).
\end{aligned}
\eeqs
Setting $A_1:=\Phi_{1,x}$, $A_2:=\Phi_{2,x}-\eta_{1,x}\Phi_{1,y}$ and $A_3:=\Phi_{3,x}-\eta_{1,x}\Phi_{2,y}-(\eta_{2,x}-\eta_1\eta_{1,x})\Phi_{1,y}$, the Bernoulli condition becomes
\beqs
\kappa^2A_1A_2+\Phi_{1,y}\Phi_{2,y}-\eta_1\Phi_{1,y}^2-\kappa A_3-k_2A_1-\kappa\gamma(\eta_1A_2+\eta_2A_1)+(\alpha+\gamma)\eta_3+\gamma^2\eta_1\eta_2=0.
\eeqs
Their first harmonics are
\beqs
Z_1'(0)=\kappa b_1+Sk_2+K_1,\qquad \kappa Z_1(0)=(\alpha+\gamma)b_1-Ck_2+M_1,
\eeqs
where
\beqs
\begin{aligned}
K_1&=\kappa^2SW(0)+\kappa S\left((\gamma+1)a_0-\f{\gamma+2\alpha-1}{2}a_2\right)-\f{\kappa S^3}{4}(\alpha+3\kappa^2),\\
M_1&=-\gamma\alpha S\left(a_0+\f{a_2}{2}\right)+\alpha\kappa SW(0)+\f{\alpha\kappa^2S^3}{4}-\kappa^2Sa_2.
\end{aligned}
\eeqs
The first boundary relation gives
\beqs
E=\f{\kappa b_1+K_1-P_1'-\kappa Ck_2}{\kappa S}.
\eeqs
Substituting this into the second one and using \eqref{relation}, the terms containing $b_1$ cancel and we find
\beqs
\left(C-\f{\kappa}{S}\right)k_2=M_1-\coth(\kappa)(K_1-P_1')-\kappa P_1.
\eeqs
This identity already determines $k_2$ independently of the choice of amplitude parametrization. Finally, substituting the formulas for $a_0$, $a_2$, $D$ and $W(0)$, together with \eqref{relation} and the identity $2CS=\sh(2\kappa)$ yields \eqref{def-k2}.
\end{proof}

\begin{rmk}
When $\gamma=0$, the formulas in \eqref{profile-expansion} and \eqref{def-k2} reduce to
\beqs
a_0=-\f{\kappa\tah\kappa}{4},\qquad a_2=\f{\kappa(2+\ch(2\kappa))}{4\tah\kappa},
\qquad
k_2=\f{\kappa^3(8\sh^6\kappa+16\sh^4\kappa+15\sh^2\kappa+9)}{2\,\sh(2\kappa)(\sh(2\kappa)-2\kappa)}.
\eeqs
which is consistent with the computation of 
small amplitude irrotational waves \cite{JRSY-cmp}. 
\end{rmk}

\begin{rmk} 
Proposition \ref{prop-expan-profile} gives the expansion of the surface velocity potential $\varphi_{\vep}$ as
\beq\label{exp-vep}
\vp_{\vep}=\varepsilon\vp_1+\varepsilon^2\vp_2+O(|\varepsilon|^3),\quad \text{with} \quad \vp_1=C\sin x, \;\;\vp_2=W(0)\sin (2x).
\eeq
\end{rmk}

\section{Discussion on the Hamiltonian structure for $L_0^{\vep}.$}\label{appen-Hamiltonian-0}
In this section, we discuss the Hamiltonian structure for $L_0^{\vep}.$
Let
$$
\langle f\rangle
:=
\frac{1}{2\pi}\int_{\T_{2\pi}}f(x)\,\d x,
\qquad
\widetilde\Pi_0f:=f-\langle f\rangle,
$$
and define
\beq\label{def-J0ep}
J_0^{\vep}
:=
\left(\begin{array}{cc}
0&1\\
-1&\gamma(k_{\vep}\p_x)_0^{-1}\widetilde\Pi_0
\end{array}\right),
\eeq
where $(k_{\vep}\p_x)_0^{-1}$ denotes the zero-mean Fourier inverse:
$$
(k_{\vep}\p_x)_0^{-1}\widetilde\Pi_0f
=
\sum_{j\neq 0}\frac{f_j}{\i k_{\vep}j}e^{\i jx}.
$$
In particular,
$$
(k_{\vep}\p_x)_0^{-1}\widetilde\Pi_0
k_{\vep}\p_xf
=
f-\langle f\rangle.
$$
On the other hand, the inverse derivative appearing in the physical operator $L_0^{\vep}$ is determined by the bottom-normalized stream function. More precisely,
$
-k_{\vep}\p_xH[\eta_{\vep}]g
=
G_0[\eta_{\vep}]g.$
Since
$
\int_{\T_{2\pi}}G_0[\eta_{\vep}]g\,\d x=0,
$
we obtain
$$
(k_{\vep}\p_x)_0^{-1}
G_0[\eta_{\vep}]g
=
-H[\eta_{\vep}]g
+
\left\langle H[\eta_{\vep}]g\right\rangle.
$$
Thus, the zero-mean primitive appearing in $J_0^{\vep}$ and the physical primitive $-H[\eta_{\vep}]g$ differ by an additive constant.
Let
$$
A_0^{\vep}
:=
\left(\begin{array}{cc}
w_{\vep}+\gamma d_{\vep}&-k_{\vep}d_{\vep}\p_x\\
k_{\vep}\p_x(d_{\vep}\cdot)&G_0[\eta_{\vep}]
\end{array}\right).
$$
It follows from the direct computation that
\beqs
J_0^{\vep}A_0^{\vep}
\binom{f}{g}
=
L_0^{\vep}\binom{f}{g}
+
\gamma
\binom{0}{
-\langle d_{\vep}f\rangle
+\left\langle H[\eta_{\vep}]g\right\rangle
}.
\eeqs
Equivalently,
\beq \label{Almost-hamiltonian}
L_0^{\vep}=J_0^{\vep}A_0^{\vep}+R^{\vep},
\eeq
where
$$
R^{\vep}\binom{f}{g}
:=
\gamma
\left(
\langle d_{\vep}f\rangle
-
\left\langle H[\eta_{\vep}]g\right\rangle
\right)
\binom{0}{1}.
$$
Therefore, the physical fixed-constant operator $L_0^{\vep}$ is not, in general, equal to the  Hamiltonian operator $J_0^{\vep}A_0^{\vep}$. This distinction is already visible when $\vep=0$, in which case, $d_0=1$ and $\langle H[0]g\rangle=0$, so that
$$
R^0\binom{f}{g}
=\gamma\langle f\rangle\binom{0}{1}.
$$
Consequently, on the constant Fourier mode,
\beqs
L_0^0
=
\left(\begin{array}{cc}
0&0\\
-\alpha&0
\end{array}\right),
\qquad
J_0^0A_0^0
=
\left(\begin{array}{cc}
0&0\\
-(\alpha+\gamma)&0
\end{array}\right).
\eeqs

\section{Construction of the basis for the generalized kernel of $L_0^{\varepsilon}$}\label{app B}
The purpose of this section is to prove Proposition \ref{prop-basis}, concerning the construction of the basis for the generalized kernel of $L_0^{\varepsilon}$. 

\subsection{Determination of $q_1^{\vep}(0,\cdot)$ and the preliminary vector $\hat q_2^{\vep}(0,\cdot)$}

To find $q_1^{\vep}(0,\cdot)$ and $\hat q_2^{\vep}(0,\cdot)$, we first work with the original linearized operator, before introducing the good unknown,
\beqs
\mathcal L_0^{\vep}:=
\begin{pmatrix}
1&0\\
-Z_\vep&1
\end{pmatrix}^{-1}
L_0^{\vep}
\begin{pmatrix}
1&0\\
-Z_\vep&1
\end{pmatrix}.
\eeqs
The steady profiles solve the profile equations with a possibly nonzero Bernoulli constant:
\beq\label{profile-with-Bernoulli-constant}
\left\{\begin{array}{l}
k_\vep\partial_x\eta_\vep+G_0[\eta_\vep]\varphi_\vep+\dfrac{\gamma}{2}k_\vep\partial_x(\eta_\vep^2)=0,\\[0.5em]
k_\vep\partial_x\varphi_\vep-\dfrac12(k_\vep\partial_x\varphi_\vep)^2+\dfrac12\dfrac{\bigl(G_0[\eta_\vep]\varphi_\vep+k_\vep^2\partial_x\varphi_\vep\,\partial_x\eta_\vep\bigr)^2}{1+|k_\vep\partial_x\eta_\vep|^2}-\alpha\eta_\vep+\gamma k_\vep\eta_\vep\partial_x\varphi_\vep-\gamma 
H_0[\eta_\vep]\varphi_\vep=\mathfrak b_\vep,
\end{array}\right.
\eeq
where, from the expansions of the waves in Proposition \ref{prop-expan-profile}, 
\beq\label{def-fb2}
\mathfrak b_\vep=\mathfrak b_2\vep^2+O(|\vep|^3), \quad \text{with} \quad \mathfrak b_2=-\big( a_0 \alpha + \frac{\kappa^2}{4} \big).
\eeq
In fact, taking the average of the second profile equation gives
\[
\begin{aligned}
\mathfrak b_\varepsilon={}&-\frac12\left\langle (k_\varepsilon\partial_x\vp_{\vep})^2\right\rangle+\frac12\left\langle \frac{(G_0[\eta_\varepsilon]\vp_{\vep}+k_\varepsilon^2\partial_x\vp_{\vep}\partial_x\eta_\varepsilon)^2}{1+|k_\varepsilon\partial_x\eta_\varepsilon|^2}\right\rangle-\alpha\langle\eta_\varepsilon\rangle\\
&+\gamma\left\langle k_\varepsilon\eta_\varepsilon\partial_x\vp_{\vep}\right\rangle-\gamma\left\langle H[\eta_\varepsilon]\vp_{\vep}\right\rangle .
\end{aligned}
\]
With the wave expansions in \eqref{ex} and \eqref{exp-vep}, the first quadratic term gives
\[
-\frac12\left\langle (k_\varepsilon\partial_x\vp_{\vep})^2\right\rangle=-\frac{\kappa^2C^2}{4}\varepsilon^2+O(|\vep|^3).
\]
and the positive quadratic term gives
\[
\frac12\left\langle \frac{(G_0[\eta_\varepsilon]\vp_{\vep}+k_\varepsilon^2\partial_x\vp_{\vep}\partial_x\eta_\varepsilon)^2}{1+|k_\varepsilon\partial_x\eta_\varepsilon|^2}\right\rangle=\frac{\kappa^2S^2}{4}\varepsilon^2+O(|\vep|^3),
\]
where we have used $G_0[0]\vp_1=\kappa S\sin x$.  Moreover,
\[
-\alpha\langle\eta_\varepsilon\rangle=-\alpha a_0\varepsilon^2+O(|\vep|^3),\qquad \gamma\left\langle k_\varepsilon\eta_\varepsilon\partial_x\vp_{\vep}\right\rangle=\frac{\gamma\kappa SC}{2}\varepsilon^2+O(|\vep|^3).
\]
It remains to compute the mean of $-\gamma H[\eta_\varepsilon]\vp_{\vep}$. From the  definition \eqref{def-Heta},
we obtain
\[
-\gamma\left\langle H[\eta_\varepsilon]\vp_{\vep}\right\rangle=-\gamma\varepsilon^2\left\langle \eta_1(x)\kappa\partial_x\Phi_1(x,0)\right\rangle+O(\varepsilon^3)=-\gamma\frac{\kappa SC}{2}\varepsilon^2+O(|\vep|^3).
\]
Combining the above estimates yields \eqref{def-fb2}.

Differentiating \eqref{profile-with-Bernoulli-constant} with respect to $x$ gives
$\mathcal L_0^\vep
\binom{\partial_x\eta_\vep}{\partial_x\varphi_\vep}=0.$
Thus we define
\beq\label{p1-q1-def}
p_1^\vep:=\vep^{-1}\binom{\partial_x\eta_\vep}{\partial_x\varphi_\vep},
\qquad
q_1^\vep(0,\cdot):=
\begin{pmatrix}
1&0\\
-Z_\vep&1
\end{pmatrix}p_1^\vep
=
\vep^{-1}\binom{\partial_x\eta_\vep}{\partial_x\varphi_\vep-Z_\vep\partial_x\eta_\vep},
\eeq
with $L_0^\vep q_1^{\varepsilon}(0, \cdot)=0.$

We next define the preliminary generalized vector $\hat q_2^\vep$. Since the nondimensional parameters $\alpha, \gamma$ depend on the physical wave speed $c$ through $\alpha=\frac{gh}{c^2}, \gamma=\frac{\tilde\gamma h}{c},$
we set
\beq\label{preliminary-p2-def}
p_{2,*}^\vep:=
\partial_\vep\binom{\eta_\vep}{\varphi_\vep}
-\frac{\partial_\vep k_\vep}{\partial_c k_\vep}
\binom{\partial_c\eta_\vep}{c^{-1}\partial_c(c\varphi_\vep)}
\eeq
and
\beq\label{q2hat-def}
\hat q_2^\vep(0,\cdot):=
\begin{pmatrix}
1&0\\
-Z_\vep&1
\end{pmatrix}p_{2,*}^\vep .
\eeq
Differentiating \eqref{profile-with-Bernoulli-constant} with respect to $\vep$ and $c$, and combining the two differentiated equations as in \eqref{preliminary-p2-def}, gives
\beq\label{preliminary-q2-identity}
\mathcal L_0^\vep p_{2,*}^\vep
=
\frac{\vep k_\vep\partial_\vep k_\vep}{c\partial_c k_\vep}\,p_1^\vep
+\mathfrak r_\vep e_2,
\qquad
e_2:=\binom{0}{1},
\eeq
with
\beq\label{r-eps-def}
\mathfrak r_\vep:=
\partial_\vep\mathfrak b_\vep
-\frac{\partial_\vep k_\vep}{\partial_c k_\vep}\, c^{-2}\partial_c (c^2 \mathfrak{b}_\vep)=2\mathfrak b_2\vep+O(\vep^2), 
\eeq
where we have used \eqref{def-fb2} and $\partial_\vep k_\vep=O(\vep)$. Equivalently, after passing to the good unknown,
\beq\label{q2hat-identity-good}
L_0^\vep\hat q_2^\vep(0,\cdot)
=
\frac{\vep k_\vep\partial_\vep k_\vep}{c\partial_c k_\vep}\,q_1^\vep(0,\cdot)
+\mathfrak r_\vep e_2.
\eeq
Now we record the expansion of $\hat q_2^\vep$. 
Since $Z_\vep=\vep Z_1+O(\vep^2)$ with $Z_1=\kappa S\sin x$, one obtains
\beq\label{q2hat-expansion}
\hat q_2^\vep
=
\binom{S\cos x}{C\sin x}
+\vep\binom{2\eta_2}{2\varphi_2-Z_1\eta_1}
+O(\vep^2).
\eeq
The residual term $\mathfrak r_\vep e_2$ in \eqref{q2hat-identity-good} will be compensated after the construction of $q_4^\vep$ in the next subsection.

\subsection{Construction of $q_3^{\varepsilon}(0, \cdot), q_4^{\varepsilon}(0, \cdot)$}
We will prove the existence of $q_3^{\varepsilon}(0, \cdot), q_4^{\varepsilon}(0, \cdot)$ in the form of \eqref{5.1}-\eqref{5.2} by the Lyapunov-Schmidt method. It suffices to find a constant $\Lambda_0$, a smooth function $\Lambda_1(\varepsilon)$ and a function $u^{\varepsilon}(x)$ in the space $\tilde{\cY}_{\text {per}}:=\left\{f \in \cY_{\text{per }}: f=(\text{even, odd})^t, f \perp_{L^2} \operatorname{ker}\left(L_0^0\right)\right\}$, such that
\begin{equation}\label{b.4}
L_0^{\varepsilon} u^{\varepsilon}=\binom{0}{1}+\left(\varepsilon \Lambda_0+\varepsilon^2 \Lambda_1(\varepsilon)\right) q_1^{\varepsilon}(0, \cdot) ,
\end{equation}
so that 
by defining
\[
q_3^\vep:=\begin{pmatrix}0\\1\end{pmatrix}
+\bigl(\vep\Lambda_0+\vep^2\Lambda_1(\vep)\bigr)q_1^\vep,
\qquad 
q_4^\vep:=-u^\vep.
\]
Then $L_0^\vep q_3^\vep(0, \cdot)=0$ and $L_0^\vep q_4^\vep(0, \cdot)=-q_3^\vep$.

We search for $u^{\varepsilon}$ of the form
\begin{equation}\label{b.5}
u^{\varepsilon}(x)=u^0(x)+\varepsilon u^1(x)+\varepsilon^2 u^2(\varepsilon, x), \quad\left(u^j \in \tilde{\cY}_{\text{per }}\right).
\end{equation}
Let $L_0^{\varepsilon}=L_0^0+\varepsilon L_0^1+\varepsilon^2 \frac{L_0^{\varepsilon}-L_0^0-\varepsilon L_0^1}{\varepsilon^2}$. Substituting \eqref{b.5} into \eqref{b.4} and collecting terms of order $\varepsilon^0, \varepsilon^1, \varepsilon^2$, we find
\begin{align}
 L_0^0 u^0&=\binom{0}{1}, \quad \text { so that } u^0=-\binom{1 / \alpha}{0}, \nonumber\\
L_0^0 u^1&=-L_0^1 u^0+\Lambda_0 q_1^0, \label{b.6}\\ 
L_0^{\varepsilon} u^2&=-\varepsilon^{-2}\left(L_0^{\varepsilon}-L_0^0-\varepsilon L_0^1\right) u^0-\varepsilon^{-1}\left(L_0^{\varepsilon}-L_0^0\right) u^1+\Lambda_0 \varepsilon^{-1}\left(q_1^{\varepsilon}-q_1^0\right)+\Lambda_1(\varepsilon) q_1^{\varepsilon}
:=\mathcal{H}^{\varepsilon}.\label{b.7}
\end{align}
Therefore, a basis for the generalized kernel of $L_0^0$ reads
\beqs
q_1^0=\binom{-S\sin x}{C\cos x},\qquad q_2^0=\binom{S\cos x}{C\sin x},\qquad
q_3^0=\binom{0}{1},\qquad q_4^0=\binom{1/\alpha}{0}\,.
\eeqs
We first claim that 
\beq\label{dual-00}
\begin{aligned}
&(J_0^0)^{-1} q_1^0(0, \cdot)=-
\begin{pmatrix}
\dfrac{\alpha\,S}{\kappa}\cos x\\
S\sin x
\end{pmatrix},
\quad
(J_0^0)^{-1}q_2^0(0,\cdot)
=
\begin{pmatrix}
-\dfrac{\alpha\,S}{\kappa}\sin x\\
S\cos x
\end{pmatrix},\\
 &(J_0^0)^{-1} q_3^0(0, \cdot)=\binom{-1}{0}, \qquad\qquad\quad (J_0^0)^{-1} q_4^0= \binom{0}{1/\alpha}
\end{aligned}
\eeq
form a basis for the generalized kernel of $(L_0^0)^{*},$
where $(J_0^0)^{-1} q_j\in \ker \left(L_0^0\right)^*\, (j=1,2,3)$ and
$(J_0^0)^{-1} q_4$ is a generalized eigenfunction of $\left(L_0^0\right)^*$. Indeed, by \eqref{Almost-hamiltonian} and since 
\begin{align*}
R^0 q_j^0(0, \cdot)=0, \qquad A_0^0 q_j^0(0, \cdot)=0, \qquad (R^0)^{*}(J_0^0)^{-1}q_j^0(0, \cdot)=0,\qquad j=1,2,3,
\end{align*}
it holds that
 \beq\label{pro=0}
 (L_0^0)^{*}(J_0^0)^{-1} q_j^0(0, \cdot)=-A_0^0 q_j^0(0, \cdot)+(R^0)^{*}(J_0^0)^{-1}q_j=0.
 \eeq
Moreover, it follows from direct computation that 
\begin{align*}
     (L_0^0)^{*}(J_0^0)^{-1} q_4= (J_0^0)^{-1} q_3\,.
\end{align*}
From the Fredholm alternative, the solvability of \eqref{b.6} is equivalent to the orthogonality condition
\begin{equation}\label{b.8}
\left\langle (J_0^0)^{-1} q_j^0(0, \cdot),-L_0^1 u^0+\Lambda_0 q_1^0(0, \cdot)\right\rangle=0, \quad j=1,2,3.
\end{equation}
 However, one can readily check that
$$
\begin{aligned}
L_0^1 u^0 & =-\left.\frac{\kappa}{\alpha}\binom{\partial_x (-v_{\varepsilon}+\gamma \eta_{\varepsilon})}{\partial_x Z_{\varepsilon}}/\vep\right|_{\varepsilon=0}=-\frac{\kappa}{\alpha}\binom{(\kappa C-\gamma\operatorname{sh}(\kappa))\sin (x)}{\kappa S\cos (x)}=\binom{\text { odd }}{\text { even }}, \\
q_1^0(0, \cdot) & =\binom{-S\sin (x)}{C\cos (x)}=\binom{\text { odd }}{\text { even }},
\end{aligned}
$$
and thus \eqref{b.8} holds for $j=1,3$ automatically. On the other hand, one can determine $\Lambda_0$ by the condition for $j=2$ :
$$
\Lambda_0=\frac{\left\langle (J_0^0)^{-1} q_2^0(0, \cdot), L_0^1 u^0\right\rangle}{\left\langle (J_0^0)^{-1} q_2^0(0, \cdot), q_1^0(0, \cdot)\right\rangle}=
\frac{\kappa\left(\alpha^2-\kappa^2\right)}{\alpha (2\alpha+\gamma)}.
$$
With this in hand, and using the fact $G_0[0]\left(e^{ \pm \mathrm{i} x}\right)=\kappa|D| \operatorname{th}(\kappa|D|)\left(e^{ \pm \mathrm{i} x}\right)=\kappa \operatorname{th}(\kappa)\left(e^{ \pm \mathrm{i} x}\right)$, and
$$
-L_0^1 u^0+\Lambda_0 q_1^0=\frac{\kappa\left(\kappa\operatorname{ch}(2\kappa)-\frac {\gamma}{2}\operatorname{sh}(2\kappa)\right)}{\left(2\alpha+ \gamma\right)\operatorname{sh}(\kappa)}\binom{\frac{\kappa}{\alpha}\sin (x)}{\cos (x)},
$$
we can solve $u^1$ explicitly by \eqref{b.6}:
$$
u^1=-\frac{\alpha + \gamma}{\alpha(2\alpha + \gamma)} \left( \kappa - \frac{\gamma}{2} \tanh(2\kappa) \right)\binom{C\cos (x)}{-S\sin (x)}.
$$

Let us now proceed to prove the existence of $\left(\Lambda_1(\varepsilon), u^2(\varepsilon, \cdot)\right)$ that satisfies \eqref{b.7}, which will be provided by the Implicit Function Theorem. Denote by $\mathbb{Q}_0$ the orthogonal projection onto the orthogonal complement in $\cY_{\text {per }}$ of $\operatorname{ker}\left(L_0^0\right)$ and by $\mathbb{P}_0$ the orthogonal projection onto the orthogonal complement in $\cY_{\text {per }}$ of $\operatorname{ker}\left(L_0^0\right)^*$. Then \eqref{b.7} is equivalent to
\begin{align}
& \tilde{L}_0^{\varepsilon} u^2=\mathbb{P}_0 \mathcal{H}^{\varepsilon}, \quad \text { with } \quad \tilde{L}_0^{\varepsilon}=\mathbb{P}_0 L_0^{\varepsilon} \mathbb{Q}_0, \label{b.9}\\
& F_j\left(\varepsilon, \Lambda_1\right)=\left\langle (J_0^0)^{-1} q_j^0(0, \cdot),\left(L_0^{\varepsilon}-L_0^0\right) u^2-\mathcal{H}^{\varepsilon}\right\rangle=0, \quad j=1,2,3,\label{b.10}
\end{align}
where we have used \eqref{pro=0} and, abusing notation slightly, regarded $\Lambda_1$ as a parameter rather than a function of $\varepsilon$.

On the one hand, the formula $\left(\tilde{L}_0^{\varepsilon}\right)^{-1}=\left(\tilde{L}_0^0\right)^{-1}\left(\operatorname{Id}+\left(\tilde{L}_0^{\varepsilon}-\tilde{L}_0^0\right)\left(\tilde{L}_0^0\right)^{-1}\right)^{-1}$ defines an inverse $\left(\tilde{L}_0^{\varepsilon}\right)^{-1}: \tilde{\tilde{\cY}}_{\text {per }} \rightarrow \tilde{\cY}_{\text {per }} \cap \operatorname{Dom}\left(L_0^0\right)$ to $\tilde{L}_0^{\varepsilon}$ with
\[
\left\|\left(\tilde{L}_0^{\varepsilon}\right)^{-1}-\left(\tilde{L}_0^0\right)^{-1}\right\|_{B\left(\tilde{\tilde{\cY}}_{\mathrm{per}}, \tilde{\cY}_{\mathrm{per}}\right)} \lesssim \varepsilon,
\]
where $\tilde{\tilde{\cY}}_{\text {per }}:=\left\{f \in \cY_{\text{per}}: f=(\text {odd, even })^t, f\perp_{L^2}\ker((L_0^0)^*)\right\}$. The existence of the function $u^2(\varepsilon, \cdot) \in \tilde{\cY}_{\text {per }}$ is thus ensured, once $\Lambda_1(\varepsilon)$ is known. On the other hand, it is observed that
$$
\mathcal{H}^{\varepsilon}=\binom{\text { odd }}{\text { even }} \quad \text { and } \quad\left(L_0^{\varepsilon}-L_0^0\right) u^2=\binom{\text { odd }}{\text { even }}.
$$
Thus \eqref{b.10} holds directly for $j=1,3$ and it suffices to find a smooth function $\Lambda_1(\varepsilon)$ such that
$$
F_2\left(\varepsilon, \Lambda_1\right)=\left\langle (J_0^0)^{-1} q_2^0(0, \cdot),\left(L_0^{\varepsilon}-L_0^0\right) u^2-\mathcal{H}^{\varepsilon}\right\rangle=0.
$$
Setting $\varepsilon=0$, one finds a constant $\Lambda_{10}$ such that $F_2\left(0, \Lambda_{10}\right)=0$. Moreover,
$$
\left(\partial_{\Lambda_1} F_2\right)\left(0, \Lambda_1\right)=-\left\langle (J_0^0)^{-1} q_2^0(0, \cdot), q_1^0(0, \cdot)\right\rangle= 
-\frac{\pi(2\alpha+\gamma)\sh^2(\kappa)}{\kappa} 
\neq 0.
$$
Consequently, by the Implicit Function Theorem, there exists $\varepsilon_3>0$ and a smooth curve $\Lambda_1(\varepsilon)$ : $\left[-\varepsilon_3, \varepsilon_3\right] \rightarrow \mathbb{R}$, such that $\Lambda_1(0)=\Lambda_{10}$ and $F_2\left(\varepsilon, \Lambda_1(\varepsilon)\right)=0$, for $\varepsilon \in\left[-\varepsilon_3,\varepsilon_3\right].$

\subsection{Construction of $q_2^{\vep}$}
We now return to the preliminary vector $\hat q_2^\vep$ constructed in \eqref{q2hat-def}. Since $q_4^\vep$ satisfies
\beqs
L_0^\vep q_4^\vep=-e_2-\Lambda_\vep q_1^\vep,
\qquad
\Lambda_\vep=\vep\Lambda_0+\vep^2\Lambda_1(\vep),
\eeqs
we define the corrected generalized vector
\beq\label{q2-corrected-def}
q_2^\vep(0,\cdot):=\hat q_2^\vep(0,\cdot)+\mathfrak r_\vep q_4^\vep(0,\cdot).
\eeq
Using \eqref{q2hat-identity-good}, we obtain
\beqs
L_0^\vep q_2^\vep
=d_{12}(\vep,0) q_1^\vep,
\qquad \text{with} \quad
d_{12}(\vep,0)=
\frac{\vep k_\vep\partial_\vep k_\vep}{c\partial_c k_\vep}
-\mathfrak r_\vep\Lambda_\vep=\left(
\frac{2\kappa k_2}{c\partial_c\kappa}
-2\mathfrak b_2\Lambda_0
\right)\vep^2+O(\vep^3).
\eeqs
Moreover, since $q_4^\vep=q_4^0+O(\vep)$ and $q_4^0=(1/\alpha,0)^t$, the corrected expansion of $q_2^\vep$ is
\beqs
q_2^\vep
=\binom{S\cos x}{C\sin x}
+\vep\left[
\binom{2\eta_2}{2\varphi_2-Z_1\eta_1}
+2\mathfrak b_2\binom{1/\alpha}{0}
\right]
+O(\vep^2).
\eeqs
Equivalently, the zero Fourier mode of the first component of $q_2^1$ is $2a_0^\sharp$, where
$a_0^\sharp=a_0+\frac{\mathfrak b_2}{\alpha}=-\frac{\kappa^2}{4\alpha}.$

\section{Proof of \eqref{signs}}\label{Appendix-secondderi}

We use the variables
\beqs
\sigma=\tah\kappa, \qquad r=1+X=\frac{\alpha}{\alpha+\gamma}>0,\qquad
y=\sh^2\kappa>0,\qquad
e=\frac{\sh(2\kappa)}{2\kappa}>1,
\eeqs
introduced in the proof of Proposition~\ref{prop:stability} in 
Section \ref{sec:sign-index}. 
 By the dispersion
relation \eqref{relation}, we have
$$
\alpha=\frac{\kappa r}{\sigma},\quad
\gamma=\frac{\kappa(1-r)}{\sigma},\quad \frac{\kappa(1-\sigma^2)}{\sigma}
=\frac1e,\quad
\kappa\sigma=\frac ye,\quad
\kappa^2(1-\sigma^2)=\frac{y}{e^2}.
$$
Recall from \eqref{def-wzetaj} that for $z>0,$ $w(z)=\frac12\sqrt{F(z)}$ with $F(z)=4\alpha z\tah z+\gamma^2\tah^2z$,
and thus 
$F(\kappa)=\kappa^2(1+r)^2$.
Using the preceding identities, we compute further that
\beqs
\begin{aligned}
F'(\kappa)
&=2\kappa\left(2r+\frac{1+r^2}{e}\right), \qquad 
F''(\kappa)=\frac{2}{e^2}
\left(4re+(1-r)^2-2(1+r^2)y\right).
\end{aligned}
\eeqs
Since
\beqs
w'(\kappa)=\frac{F'(\kappa)}{4\sqrt{F(\kappa)}},
\qquad
w''(\kappa)
=\frac{2F(\kappa)F''(\kappa)-F'(\kappa)^2}
{8F(\kappa)^{3/2}},
\eeqs
it follows from \eqref{signs} that
$$
\nu_*=\frac{r+e^{-1}}{1+r}<1,\qquad \delta_*
=\frac{r^2(e-1)^2+(1+r)^2y}
{\kappa e^2(1+r)^3}>0.
$$
This completes the proof.

\section{Some properties of the extension operators}

In this section, we collect and prove some useful properties of the extension operators
$\cU^{\vep}(\xi), \cV^{\vep}(\xi)$ defined in \eqref{ODE-cU} and
\eqref{ODE-cV}.
\begin{prop}\label{prop-extop}
Let $\cK g:=\overline{g}$ denote complex conjugation and define the reversibility operator: 
\beqs 
\mathcal{I} g =\mathrm{diag}\, (1,-1) \cK({g(-x)}).
\eeqs
The following properties hold. 

\textnormal{(1)}
\beq\label{ext-conj}
\cU^{\vep}(-\xi)=\cK\,\cU^{\vep}(\xi)\,\cK,\qquad
\cV^{\vep}(-\xi)=\cK\,\cV^{\vep}(\xi)\,\cK.
\eeq

 \textnormal{(2)} The extension operators $\cU^{\vep}(\xi), \cV^{\vep}(\xi)$ commute with the reversibility operator: 
\beq\label{reversibilitycom}
\mathcal{U}^{\vep}(\xi)\mathcal{I}=\mathcal{I}\,\mathcal{U}^{\vep}(\xi), \qquad \mathcal{V}^{\vep}(\xi)\mathcal{I}=\mathcal{I}\,\mathcal{V}^{\vep}(\xi),
\eeq

 \textnormal{(3)}
 Let
 \beqs
 \cU_1^\vep:=\frac{1}{\i\kappa}\p_{\xi} \cU^\vep(\xi)\big|_{\xi=0},\qquad \cV_1^\vep:=\frac{1}{\i\kappa}\partial_\xi \cV^\vep(\xi)\big|_{\xi=0}.
\eeqs
then $\cU_1^\vep$ and $\cV_1^\vep$ preserve real-valuedness. 
\end{prop}

\begin{proof}
Proof of  \textnormal{(1).}
Since $L^\vep$ has real coefficients, $\cK L_\xi^\vep\cK=L_{-\xi}^\vep$ and $\cK(L_\xi^\vep)^*\cK=(L_{-\xi}^\vep)^*$; as $\partial B_{r_0}(0)$ is invariant under conjugation, the same holds for $\Pi_\xi^\vep$ and $(\Pi_\xi^\vep)^*$.
Hence $\cK\,\cU^\vep(-\xi)\,\cK$ and $\cK\,\cV^\vep(-\xi)\,\cK$ solve \eqref{ODE-cU} and \eqref{ODE-cV} with initial value $\mathrm{Id}$, and \eqref{ext-conj} follows by uniqueness.

Proof of  \textnormal{(2).} 
From the parity of the background waves, we readily find that $$L_{\xi}^{\vep}\mathcal{I}=-\mathcal{I}\,L_{\xi}^{\vep}, \quad (L_{\xi}^{\vep})^{*}\mathcal{I}=-\mathcal{I}\,(L_{\xi}^{\vep})^{*}, $$
which implies that $$\Pi_{\xi}^{\vep}\mathcal{I}=\mathcal{I}\,\Pi_{\xi}^{\vep}, \qquad (\Pi_{\xi}^{\vep})^{*}\mathcal{I}=\mathcal{I}\,(\Pi_{\xi}^{\vep})^{*}.$$ 
The identity 
\eqref{reversibilitycom} then follows from the 
uniqueness of the ODEs system \eqref{ODE-cU}, \eqref{ODE-cV} satisfied by $\cU^{\vep}(\xi), \cV^{\vep}(\xi).$

Proof of  \textnormal{(3).} Let $L_1^\vep:=\frac{1}{\i\kappa}\partial_\xi L_\xi^\vep\big|_{\xi=0}, 
\Pi_1^\vep:=\frac{1}{\i\kappa}\partial_\xi\Pi_\xi^\vep\big|_{\xi=0}.$
Differentiating the Riesz formula for $\Pi_{\xi}^\vep$ gives
\beqs
\Pi_1^\vep=\frac{1}{2\pi \i}\int_{\partial B_{r_0}(0)}(\lambda I-L_0^\vep)^{-1}L_1^\vep(\lambda I-L_0^\vep)^{-1}\,\d\lambda .
\eeqs
Since $L_0^\vep$ and $L_1^\vep$ are real operators, and since the contour $\partial B_{r_0}(0)$ is invariant under complex conjugation, the operator $\Pi_1^\vep$ also preserves real-valuedness. The same is true for ${\Pi}_0^\vep$. Therefore
\beqs
\cU_1^\vep=\frac{1}{\i\kappa}\partial_\xi\cU^\vep(\xi)\big|_{\xi=0}=[\Pi_1^\vep,{\Pi}_0^\vep]
\eeqs
preserves real-valuedness. The similar argument applies $\cV_1^{\vep}.$
\end{proof}

\begin{lemma}\label{lem-extension-zero}
Let $V\in\mathbb C^2$ be constant. Then $\cU^0(\xi)V=V$. Moreover, there exist smooth non-vanishing functions $b_{\pm1}(\xi)$ such that
\beq\label{extension-pm1}
\cU^0(\xi)\big(v_{\pm1}(0,\cdot)\big)=b_{\pm1}(\xi)v_{\pm1}(\xi,\cdot).
\eeq
where $v_{\pm1}(\xi,\cdot)=v_{\pm1}^{+}(\xi,\cdot)$ are defined in \eqref{lam}.
\end{lemma}

\begin{proof}
For $j\in\mathbb Z$, set $E_j:=\{e^{\i jx}V:\ V\in\mathbb C^2\}$. Since $L_\xi^0$ has constant coefficients, each $E_j$ is invariant under $L_\xi^0$, its resolvent, and the Riesz projector $\Pi_\xi^0$. The same is therefore true for $\partial_\xi\Pi_\xi^0$ and $[\partial_\xi\Pi_\xi^0,\Pi_\xi^0]$, and uniqueness in the Kato equation implies
\beqs
\cU^0(\xi)E_j\subset E_j.
\eeqs
On $E_0$, both eigenvalues belong to the spectral region enclosed by the contour, so that $\Pi_\xi^0|_{E_0}=\Id$. Hence the Kato generator vanishes on $E_0$ and $\cU^0(\xi)V=V$ for every constant vector $V$.

We now consider $E_{\pm1}$. For $|\xi|$ sufficiently small, the restriction of $\Pi_\xi^0$ to $E_{\pm1}$ has rank one and
\beqs
E_{\pm1}\cap{\rm Ran}\,\Pi_\xi^0={\rm Span}\,\{v_{\pm1}(\xi,\cdot)\}.
\eeqs
Since the Kato extension intertwines the spectral subspaces,
\beqs
\cU^0(\xi)\rm{Ran}\,\Pi_0^0=\rm{Ran}\,\Pi_\xi^0.
\eeqs
For $|\xi|$ sufficiently small, $\lambda_{\pm1}^+(\xi)$ are the only eigenvalues of $L_\xi^0|_{E_{\pm1}}$ inside the contour. Consequently,
\beqs
E_{\pm1}\cap{\rm Ran}\Pi_\xi^0= {\rm Span}\{v_{\pm1}(\xi,\cdot)\}.
\eeqs
Since $v_{\pm1}(0,\cdot)\in E_{\pm1}\cap\rm{Ran}\Pi_0^0$, there exist smooth functions $b_{\pm1}(\xi)$ such that
\beq\label{relationI}
\cU^0(\xi)v_{\pm1}(0,\cdot)=b_{\pm1}(\xi)v_{\pm1}(\xi,\cdot).
\eeq
The coefficients $b_{\pm1}$ may also be characterized explicitly. In fact, let $\widetilde v_{\pm1}(\xi,\cdot)$ be the corresponding adjoint eigenfunctions, normalized by $\ang{\widetilde v_{\pm1}(\xi,\cdot),v_{\pm1}(\xi,\cdot)}=1.$
Then the restriction of the projector to $E_{\pm1}$ is
\beqs
\Pi_\xi^0f=\ang{\widetilde v_{\pm1}(\xi,\cdot),f}\,v_{\pm1}(\xi,\cdot),\qquad f\in E_{\pm1}.
\eeqs
Differentiating the  relation \eqref{relationI} and using the definition of $\cU^0(\xi),$  the identities 
\beqs
(\p_{\xi}\Pi_\xi^0) ( v_{\pm 1})=(I-\Pi_\xi^0) (\p_{\xi} v_{\pm 1}), \qquad \Pi_\xi^0 (\p_{\xi}\Pi_\xi^0)\Pi_\xi^0 =0
\eeqs
we find
\beqs
b_{\pm 1}'(\xi)v_{\pm1}=-b_{\pm1}(\xi)\Pi_\xi^0 (\partial_\xi v_{\pm1}).
\eeqs
Taking the scalar product with $\widetilde v_{\pm1}$ yields
\beqs
b_{\pm1}'(\xi)=-b_{\pm1}(\xi)\ang{\widetilde v_{\pm1}(\xi,\cdot),\partial_\xi v_{\pm1}(\xi,\cdot)},
\eeqs
and consequently
\beqs
b_{\pm1}(\xi)=b_{\pm1}(0)\exp\left(-\int_0^\xi\ang{\widetilde v_{\pm1}(r,\cdot),\partial_rv_{\pm1}(r,\cdot)}\,\d r\right).
\eeqs 
\end{proof}

\section{Linearized operator expansions and some elementary computations. }

Recall that the 
physical Bloch stream function in the linearized operator $L_{\xi}^{\vep}$ given in \eqref{def-G-xi-esp} is defined by
\beq\label{def-Hxi-physical}
H_\xi[\eta_{\vep}]f(x)
:=
\int_{-1}^{\eta_{\vep}(x)}
k_{\vep}(\p_x+\i\xi)\Phi_{\xi}^{\vep}(x,y)\,\d y,
\eeq 
where $\Phi_{\xi}^{\vep}=\Phi_{\xi}[\eta_\vep](f)$ solves the elliptic problem
\beqs
\left(\p_y^2+k_{\vep}^2(\p_x+\i\xi)^2\right)\Phi_{\xi}^{\vep}=0
\quad\hbox{in}\quad
\Omega_{\vep}:=\{(x,y): -1<y<\eta_{\vep}(x)\},
\eeqs
with boundary conditions
\[
\Phi_{\xi}^{\vep}(x,\eta_{\vep}(x))=f(x),\qquad
\p_y\Phi_{\xi}^{\vep}(x,-1)=0.
\]
Also recall the notations
\beqs
L_0^{\vep}= L_{\xi}^{\vep}\big|_{\xi=0}, \qquad 
L_{1}^{\vep}=\frac{1}{\i\kappa}\partial_{\xi}L_{\xi}^{\vep}\big|_{\xi=0},\qquad
L_{k}^{\vep}=L_{k}^{0}+\vep L_{k}^{1}+O(\vep^2), \qquad k=0,1.
\eeqs

\begin{lem}\label{lem-Hxi-expansion}
For every $s>\frac12$, the maps $
(\vep,\xi)\mapsto H_\xi[\eta_\vep],
\,(\vep,\xi)\mapsto G_\xi[\eta_\vep]$
are analytic near $(0,0)$ with values respectively in
\beqs
B\big(H^{s+\frac12}(\T),H^{s+\frac12}(\T)\big)
\quad\text{and}\quad
B\big(H^{s+\frac12}(\T),H^{s-\frac12}(\T)\big).
\eeqs
Define
\begin{align*}
    G_{\ell}^k= \big(\f{1}{\i\kpa}\p_{\xi}\big)^{\ell} \p_{\vep}^k\big(G_{\xi}[\eta_{\vep}]\big)\big|_{\xi=\vep=0}, \qquad  H_{\ell}^k= \big(\f{1}{\i\kpa}\p_{\xi}\big)^{\ell} \p_{\vep}^k\big(H_{\xi}[\eta_{\vep}]\big)\big|_{\xi=\vep=0}.
\end{align*}
The following identities hold:

\textnormal{(1)}
\beq\label{expan-stream}
\begin{aligned}
 H_0^0=\i\tah(\kappa D),\qquad &
H_1^0=\sech^2(\kappa D), \\ 
 H_0^1=\eta_1\kappa\p_x 
-
H_0^0\left(\eta_1G_0^0\right), \qquad & H_1^1=\eta_1 
-
H_1^0\left(\eta_1G_0^0\right)
-
H_0^0\left(\eta_1G_1^0\right).
\end{aligned}
\eeq

\textnormal{(2)}
\beq\label{expan-DN}
\begin{aligned}
G_0^0=|\kpa D|\tah |\kpa D|\,, & \qquad G_1^0=-\i\, \tah(\kpa D)-\f{\kpa\p_x}{\ch^2(\kpa D)}, \\ 
 G_0^1=-G_0^0\,\eta_1G_0^0-\kpa \p_x(\eta_1 \kpa\p_x), &\qquad G_1^1=-\big(G_0^0\,\eta_1G_1^0+G_1^0\,\eta_1 G_0^0\big)-\kpa \big(\p_x(\eta_1\cdot)+\eta_1\p_x\big).
\end{aligned}
\eeq
\end{lem}

\begin{proof}
We first prove the analyticity. 
Introduce the change of variables
\beqs
y=(1+\eta_\vep)(z+1)-1,
\eeqs
and let $\phi_{\xi}^{\vep}
=\phi_\xi[\eta_{\vep}](f)$ solve the problem on the fixed strip $\cS=\T\times[-1,0]$:
\beq\label{elliptic-strip}
\Div_\xi\big(\mathbb A_\vep\nabla_\xi\phi_{\xi}^{\vep}\big)=0\quad\hbox{in }\mathcal S,\qquad \phi_{\xi}^{\vep}\big|_{z=0}=f,\qquad \p_z\phi_{\xi}^{\vep}\big|_{z=-1}=0.
\eeq
where
\beqs
\begin{aligned}
\nabla_\xi:=&\binom{\p_x+\i\xi}{\p_z},\qquad \Div_\xi(F_1,F_2):=(\p_x+\i\xi)F_1+\p_zF_2,\\
\mathbb A_\vep:=&
\begin{pmatrix}
k_\vep^2(1+\eta_\vep)&-k_\vep^2(z+1)\p_x\eta_\vep\\
-k_\vep^2(z+1)\p_x\eta_\vep&\dfrac{1+k_\vep^2(z+1)^2(\p_x\eta_\vep)^2}{1+\eta_\vep}
\end{pmatrix}.
\end{aligned}
\eeqs
For $s>\frac12$, consider the operator 
\beqs
\mathfrak A(\vep,\xi):
H^{s+1}(\mathcal S)\rightarrow 
H^{s-1}(\mathcal S)\times H^{s+\frac12}(\T)\times H^{s-\frac12}(\T):=X^s
\eeqs
defined by
\beqs
\mathfrak A(\vep,\xi)\phi:=\left(\Div_\xi(\mathbb A_\vep\nabla_\xi\phi),\,\phi\big|_{z=0},\,\p_z\phi\big|_{z=-1}\right).
\eeqs
Since $\vep\mapsto(k_\vep,\eta_\vep)$ is analytic and $1+\eta_\vep$ stays bounded away from zero for sufficiently small $|\vep|$, the map $(\vep,\xi)\mapsto\mathfrak{A} (\vep,\xi)$ is analytic between a neighborhood of $(0,0)$ and $B(H^{s+1}(\mathcal S), X^s).$
At $(\vep,\xi)=(0,0)$,
\beqs
\mathfrak A(0,0)\phi=\left((\p_z^2+\kappa^2\p_x^2)\phi,\phi\big|_{z=0},\p_z\phi\big|_{z=-1}\right)
\eeqs
is invertible. Hence, for $(\vep,\xi)$ sufficiently small,
\beqs
\mathfrak A(\vep,\xi)^{-1}=\mathfrak A(0,0)^{-1}\left(\Id+\big(\mathfrak A(\vep,\xi)-\mathfrak A(0,0)\big)\mathfrak A(0,0)^{-1}\right)^{-1}
\eeqs
is analytic, and therefore
\beqs
\phi_{\xi}^{\vep}[f]=\mathfrak A(\vep,\xi)^{-1}(0,f,0)
\eeqs
depends analytically on $(\vep,\xi)$ and linearly on $f$.

After the same change of variables, the definition \eqref{def-Hxi-physical} becomes
\beqs
H_\xi[\eta_\vep]f=k_\vep\int_{-1}^0(1+\eta_\vep)(\p_x+\i\xi)\phi_{\xi}^{\vep}-(z+1)\p_x\eta_\vep\,\p_z\phi_{\xi}^{\vep}\,\d z:=\int_{-1}^0 F_{1,\xi}^{\vep}(x,z)\, \d z.
\eeqs
By the elliptic equation \eqref{elliptic-strip}, 
\begin{align*}
    k_{\vep}(\p_x+\i \,\xi)F_{1,\xi}^{\vep}=-\p_z F_{2,\xi}^{\vep},
\end{align*}
where 
\begin{align*}
F_{2,\xi}^{\vep}=-k_\vep^2(z+1)\p_x\eta_\vep (\p_x+\i \xi)\phi_{\xi}^{\vep}+\dfrac{1+k_\vep^2(z+1)^2(\p_x\eta_\vep)^2}{1+\eta_\vep}\p_z \phi_{\xi}^{\vep}
\end{align*}
belongs to $ H^{s}(\cS)$ and vanishes on $z=-1.$ Consequently, $H_\xi[\eta_\vep]\in B(H^{s+\f12}(\T),H^{s+\f12}(\T)).$ Moreover, the transformed Dirichlet--Neumann operator takes the form
\beqs
G_\xi[\eta_\vep]f=\left.-k_\vep^2\p_x\eta_\vep(\p_x+\i\xi)\phi_{\xi}^{\vep}+\frac{1+k_\vep^2(\p_x\eta_\vep)^2}{(1+\eta_{\vep})}\p_z\phi_{\xi}^{\vep}\right|_{z=0}.
\eeqs
From the analyticity of  $\phi_{\xi}^{\vep}[f]$, we derive the analyticity of
\begin{align*}
(\vep,\xi)\longmapsto H_\xi[\eta_\vep]\in B\big(H^{s+\frac12}(\T),H^{s+\frac12}(\T)\big),\quad \text{and }
(\vep,\xi)\longmapsto G_\xi[\eta_\vep]\in B\big(H^{s+\frac12}(\T),H^{s-\frac12}(\T)\big).
\end{align*}

We next compute the coefficients in \eqref{expan-stream} and \eqref{expan-DN}. It is simpler to go back to the original definition \eqref{def-Hxi-physical}.
For $\eta_{\vep}=0$ and $f=e^{\i jx}$,
\beqs
\Phi_{\xi}^0(x,y)=e^{\i jx}\frac{\ch\big((y+1)\zeta_j\big)}{\ch\big(\zeta_j\big)}, \quad \hbox{where} \quad \zeta_j=\kappa(j+\xi).
\eeqs
Consequently,
\beqs
H_\xi[0]e^{\i jx}=\i\,\tah\big(\zeta_j\big)e^{\i jx},\qquad G_\xi[0]e^{\i jx}=\zeta_j\tah\big(\zeta_j\big)e^{\i jx}.
\eeqs
With $D=-\i\p_x$, differentiation at $\xi=0$ gives the corresponding first-order coefficients $H_0^0, H_1^0, G_0^0, G_1^0$ in \eqref{expan-stream} and \eqref{expan-DN}.
It remains to compute the first shape derivatives. Let $h$ be a variation of the free surface.
We compute from \eqref{def-Hxi-physical} that
\beqs
\f{\d}{\d \vep} H_\xi[\vep h]\bigg|_{\vep=0}(f)=-H_\xi^0\big(h\,G_\xi^0f\big)+h\,\kappa(\p_x+\i\xi)f.
\eeqs
Moreover, the shape derivative for the Dirichlet-Neumann operator (see \cite[(3.40)]{Lannes}) gives
\beqs
\f{\d}{\d \vep} G_\xi[\vep h]\bigg|_{\vep=0}(f)=-G_\xi^0\big(hG_\xi^0f\big)-\kappa^2(\p_x+\i\xi)\big(h(\p_x+\i\xi)f\big).
\eeqs
These lead to  \eqref{expan-stream} and \eqref{expan-DN}.
\end{proof}

\begin{corollary}[Some algebraic computations]\label{cor-algebraic}

Let $q_j^0(0,\cdot)$ be defined in \eqref{def-basis-00}.
Suppressing the dependence on the arguments $(0,\cdot)$ on the basis, we have that
\beqs
\begin{aligned}
L_{1}^{0}q_{1}^{0}&=
\begin{pmatrix}
\dfrac{\kappa}{C}\sin x\\
\left(C-\dfrac{\gamma}{C}\right)\cos x
\end{pmatrix},\quad
L_{1}^{0}q_{2}^{0}=
\begin{pmatrix}
-\dfrac{\kappa}{C}\cos x\\
\left(C-\dfrac{\gamma}{C}\right)\sin x
\end{pmatrix}, \quad  L_{1}^{0}q_{3}^{0}=
\begin{pmatrix}0\\ 1-\gamma\end{pmatrix},  \\
(L_1^0q_2^1)_1&=
2a_0^{\sharp}+* \cos(2x), \quad 
L_{1}^{1}q_{3}^{0}
=\kappa\begin{pmatrix}
S\sin x\\
-C\cos x
\end{pmatrix}, \quad (L_{1}^{1}q_{2}^{0})_1=d_1S\cos x+G_1^1(C\sin x).
\end{aligned}
\eeqs
In particular,
\begin{align*}
\int_0^{2\pi}(L_1^0q_2^1)_1\,\d x=4\pi a_0^{\sharp}=-\frac{\pi \kappa^2}{\alpha},\qquad
\int_0^{2\pi}(L_{1}^{1}q_{2}^{0})_1\,dx
=\pi(\gamma S^2-2\kappa CS)
=-\pi\kappa P_\gamma.
\end{align*}
\end{corollary}
\begin{proof}
Recalling  \eqref{L0} and definition of the operator $L_{1}^{0}$ from \eqref{exp-operators} and \eqref{exp-L1}, we have
\beqs
L_{1}^{0}=
\begin{pmatrix}
1 & -\i \left(\tah(\kappa D)+\kappa D\,{\sech^2(\kappa D)}\right)\\
0 & 1-\gamma {\sech^2((\kappa D))}
\end{pmatrix}, \qquad (D=-\i\p_x).
\eeqs
Therefore, on the Fourier mode $e^{\i jx}$, with $\mu_j:=\kappa j$, the operator $L_{1}^{0}$ has multiplier
\beqs
L_{1}^{0}(j)=
\begin{pmatrix}
1 & -\i\big(\tah\mu_j+\frac{\mu_j}{\ch^2\mu_j}\big)\\
0 & 1-\frac{\gamma}{\ch^2\mu_j}
\end{pmatrix}.
\eeqs
The identities involving $L_1^0$ can then be readily computed by using the definitions 
\eqref{def-basis-00}.

We now compute $L_{1}^{1}$. As $Z_{\vep}=\vep G_0[0]\varphi_1+O(\vep^2)=\vep\kappa S\sin x+O(\vep^2),$
\beqs
\begin{aligned}
v_{\vep}&=\vep v_1+O(\vep^2),\qquad \quad \text{with } v_1=\kappa\partial_x\varphi_1=\kappa C\cos x,\\
d_{\vep}&=1+\vep d_1+O(\vep^2),\quad\text{with }
d_1=-v_1+\gamma\eta_1=-(\kappa C-\gamma S)\cos x=-\alpha S\cos x,
\end{aligned}
\eeqs
we obtain
\beqs
L_{1}^{1}=
\begin{pmatrix}
d_1 & G_{1}^{1}\\
0 & d_1-\gamma H_{1}^{1}
\end{pmatrix},
\eeqs
where $H_1^1, G^1_1$ are defined in \eqref{expan-stream}, \eqref{expan-DN}.
Therefore
\beqs
G_{1}^{1}(1)=-\kappa\partial_x\eta_1=\kappa S\sin x,\qquad
H_{1}^{1}(1)=\eta_1=S\cos x.
\eeqs
The properties for $L_{1}^{1}q_{3}^{0}, (L_{1}^{1}q_{2}^{0})_1$ follow from direct computations, and we omit the details.

\end{proof}

\vspace{0.5cm}
\noindent {\bf Acknowledgements}
The work of L. Fan is partially supported by Tianyuan Fund for Mathematics of the NSFC (Grant No. 12426607) and the NSF of Henan Province of China (Grant No. 252300421218), the work of C. Sun is partially supported by the ANR project ANR-24-CE40-3260. C. Sun is grateful to Henan Normal University for its kind hospitality during his visits.

\vspace{0.5cm}
\noindent {\bf Conflict of interest}
The authors declare that there are no conflicts of interest.
\vspace{0.5cm}

\noindent {\bf Data Availability}
No external datasets were used in this study. The numerical evaluations and figures were generated
from the explicit formulas provided in the manuscript.
The corresponding code is available from the authors
upon reasonable request.

\vspace{0.5cm}

\noindent {\bf Declaration of AI use.}
An initial version of this manuscript was completed on 10 July 2026, which contains proof of the main result, Theorem \ref{thm-main}, and the verification of agreement between the resulting stability criterion and that predicted by \cite{TKM2012}. The proof strategies of these results were developed by the authors. During subsequent revisions, one of the authors used OpenAI's GPT-6.1 Sol 
to develop a rigorous but elementary proof that $\mathcal C_*>0$ for $X\in(-1,-2/3]$, as presented in Section \ref{sec:sign-index}. 
ChatGPT was also used to assist with proofreading and language editing. 
The AI-assisted arguments in Section 6.2 were  checked and refined by the authors, who take responsibility for their correctness.

\bibliographystyle{siam}
\bibliography{ww-vorticity}
\end{document}